\documentclass[11pt,a4paper,english,reqno,a4paper]{amsart}
\usepackage{amsmath,amssymb,amsthm, comment,graphicx} \usepackage{tikz}
\usepackage{fancyhdr} \usepackage{epsfig} \usepackage{mathrsfs}

\usepackage{pifont}
\usepackage{amsfonts}
\usepackage{geometry}
\usepackage{graphicx}
\usepackage[compress,sort]{cite}
\usepackage[title]{appendix}
\usepackage[colorlinks,linkcolor=blue,anchorcolor=blue,citecolor=blue]{hyperref}
\allowdisplaybreaks
\newcommand{\RNum}[1]{\uppercase\expandafter{\romannumeral #1\relax}}
\numberwithin{equation}{section}
\numberwithin{figure}{section}

\newcounter{En}[section]
\renewcommand{\theEn}{\thesection.\arabic{En}}
\newcommand{\Eqn}[1][]{\refstepcounter{En}(\theEn)}

\newcounter{As}[section]
\renewcommand{\theAs}{\textbf{A\arabic{As}}}
\newcommand{\Assumption}[1][]{\noindent\refstepcounter{As}(\theAs)}

\newcounter{TH}[section]
\renewcommand{\theTH}{\thesection.\arabic{TH}}
\newcommand{\Theorem}[1][]{\noindent\textbf{Theorem}\ \refstepcounter{TH}\textbf{\theTH.}\ }

\newcounter{CA}
\renewcommand{\theCA}{\arabic{CA}}
\newcommand{\Case}[1][]{\noindent\textbf{Case}\ \refstepcounter{CA}\textbf{\theCA.}\ }

\newcounter{LM}[section]
\renewcommand{\theLM}{\thesection.\arabic{LM}}
\newcommand{\Lemma}[1][]{\noindent\textbf{Lemma}\ \refstepcounter{LM}\textbf{\theLM.}\ }

\newcounter{PS}[section]
\renewcommand{\thePS}{\thesection.\arabic{PS}}
\newcommand{\Proposition}[1][]{\noindent\textbf{Proposition}\ \refstepcounter{PS}\textbf{\thePS.}\ }
\newcounter{DF}[section]
\renewcommand{\theDF}{\thesection.\arabic{DF}}
\newcommand{\Definition}[1][]{\noindent\textbf{Definition}\ \refstepcounter{DF}\textbf{\theDF.}\ }

\newcounter{Rm}[section]
\renewcommand{\theRm}{\thesection.\arabic{Rm}}
\newcommand{\Remark}[1][]{\noindent\textbf{Remark}\ \refstepcounter{Rm}\textbf{\theRm.}\ }

\newcounter{Co}[section]
\renewcommand{\theCo}{\thesection.\arabic{Co}}
\newcommand{\Corollary}[1][]{\noindent\textbf{Corollary}\ \refstepcounter{Co}\textbf{\theCo.}\ }

\usepackage{enumitem}
\allowdisplaybreaks[4]

\newcounter{math}
\addtocounter{math}{1}
\newcommand{\dd}{{\rm d}}

\begin{document}
	\title[Stability of Inverse Problems for Steady Supersonic Flows Past Perturbed Cones]{Stability of Inverse Problems for Steady Supersonic Flows Past Lipschitz Perturbed Cones}
\date{\today}
\author{Gui-Qiang G. Chen}
\address{ Gui-Qiang G. Chen: \, Mathematical Institute, University of Oxford,
Radcliffe Observatory Quarter, Woodstock Road, Oxford, OX2 6GG, UK;
School of Mathematical Sciences, Fudan University, Shanghai 200433, China}
\email{\tt  chengq@maths.ox.ac.hk}

\author{Yun Pu}
\address{ Yun Pu: \, Academy of Mathematics and Systems Science, Chinese Academy of Sciences, Beijing 100190, China;
School of Mathematical Sciences, Fudan University, Shanghai 200433, China}
\email{\tt ypu@amss.ac.cn}

\author{Yongqian Zhang}
\address{Yongqian Zhang: \, School of Mathematical Sciences, Fudan University, Shanghai 200433, China}
\email{\tt  yongqianz@fudan.edu.cn}


\begin{abstract}
We are concerned with inverse problems for supersonic potential flows past infinite axisymmetric Lipschitz cones.
The supersonic flows under consideration are governed by the steady isentropic Euler equations
for axisymmetric potential flows, which give rise to a singular geometric source term.
We first study the inverse problem for the stability of an oblique conical shock
as an initial-boundary value problem with both the generating curve of the cone surface
and the leading conical shock front as free boundaries.
We then establish the existence and asymptotic behavior of global entropy solutions
with bounded BV norm
of the inverse problem, under the condition that the Mach number of the incoming flow is sufficiently large
and the total variation of the pressure distribution on the cone is sufficiently small.
To this end, we first develop a modified Glimm-type scheme to construct approximate solutions
by self-similar solutions as building blocks to balance the influence of the geometric source term.
Then we define a Glimm-type functional, based on the local interaction estimates between weak waves,
the strong leading conical shock, and self-similar solutions.
Meanwhile, the approximate generating curves of the cone surface are also constructed.
Next, when the Mach number of the incoming flow is sufficiently large,
by asymptotic analysis of the reflection
coefficients in those interaction estimates,
we prove that appropriate weights can be chosen
so that the corresponding Glimm-type functional decreases in the flow direction.
Finally, we determine the generating curves of the cone surface and establish the existence
of global entropy solutions containing a strong leading conical shock, besides weak waves.
Moreover, the entropy solution is proved to approach asymptotically the self-similar solution
determined by the incoming flow and the asymptotic pressure on the cone surface at infinity. 	
\end{abstract}
\maketitle

\tableofcontents

\section{Introduction}\label{Section-Intro}
\begin{figure}
\begin{center}
\includegraphics{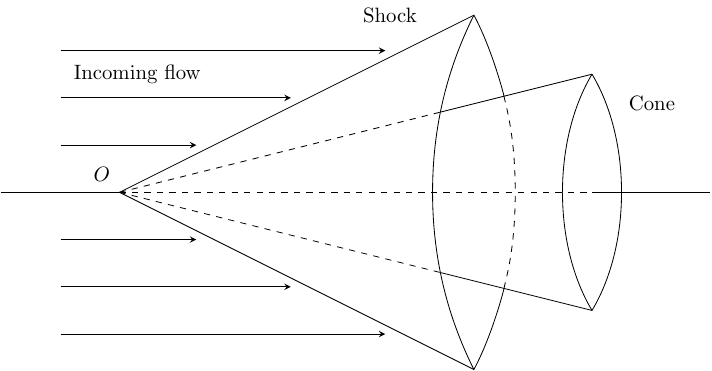}
\end{center}
\caption{The strong straight-sided conical shock}
\label{fig:1.0}
\end{figure}
\begin{figure}
\begin{center}
\includegraphics{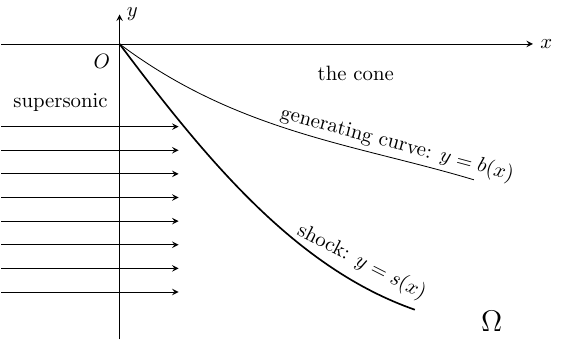}
\end{center}
\caption{Supersonic flow past an axisymmetric cone}
\label{fig:1.1}
\end{figure}

We are interested in the structural stability of inverse problems for the three-dimensional (3-D) steady
supersonic potential flows past a Lipschitz perturbed cone with given states of the incoming flow
together with Lipschitz perturbed pressure distributions on its surface.
The shock stability problem of steady supersonic flows past Lipschitz cones is fundamental for
the mathematical theory of the multidimensional (M-D) hyperbolic systems of conservation laws,
since its solutions are time-asymptotic states and
global attractors of general entropy solutions of time-dependent initial-boundary value problems (IBVP)
with abundant nonlinear phenomena, besides its significance to many fields
of applications including aerodynamics; see \cite{Anderson2019,Chen&Feldman2018,Hu2019,Courant1948}
and references cited therein.
Meanwhile, the corresponding inverse problems play essential roles in airfoil design;
see \cite{Abbott1959,Abzalilov2005,Caramia2019,Goldsworthy1952,Golubkin1988,Golubkin1994,Maddalena2020,
Mohammadi2001,Robinson1956,Vorobev1998}.
As indicated in \cite{Courant1948}, when a uniform supersonic flow of constant speed from
the far-field (negative infinity) hits a straight cone,
given a constant pressure distribution that is less than a critical value on the cone surface,
the vertex angle of the cone can be determined such that
there is a supersonic straight-sided conical shock attached to the cone vertex,
and the state between the conical shock-front and the cone can be obtained by the shooting method,
which is a self-similar solution; see Fig. \ref{fig:1.0}.
In this paper, we focus our analysis on the stability of an inverse problem,
along with the background self-similar solutions, in the steady potential flows
that are axisymmetric with respect to the $x$--axis,
given the pressure distributions of gas
on the cones, whose boundary surfaces in $\mathbb{R}^3$,
formed by the rotation of generating curves of form $\Gamma:=
\{(x,b(x))\, : \, x\geq0\}$
around the $x$--axis, are to be determined; see Fig. \ref{fig:1.1}.
\par

To be precise, the governing $3$-D Euler equations for steady potential conical flows are of the form:
\begin{equation}
\left\{\begin{aligned}
&(\rho u)_x+(\rho v)_y=-\frac{\rho v}{y}, \\
&v_x-u_y=0,
\end{aligned}
\right.
\label{eqn:orig}
\end{equation}
together with the Bernoulli law:
\begin{equation}
\frac{u^2+v^2}{2}+\frac{c^2}{\gamma-1}= \frac{u_{\infty}^2}{2}+\frac{c_{\infty}^2}{\gamma-1},
\label{eqn:Bernoulli}
\end{equation}
where $U:=(u,v)^\top$ is the velocity in the $(x,y)$--coordinates,
$\rho$ is the flow density, and $U_{\infty}=(u_{\infty},0)^\top$ and $\rho_{\infty}$
are the velocity and the density
of the incoming flow, respectively.
The Bernoulli law in (\ref{eqn:Bernoulli}) is obtained
via the constitutive relation between pressure $p$ and density $\rho$: $p=\rho^{\gamma}$,
with $\gamma>1$ for the polytropic isentropic
gas and $\gamma=1$ for the isothermal flow, under scaling.
In particular, $c=:\sqrt{\frac{\gamma p}{\rho}}$ is called the sonic speed, and
$M:=\sqrt{\frac{u^2+v^2}{c^2}}$ is called the Mach number.

The Bernoulli law (\ref{eqn:Bernoulli}) can be written as
\begin{equation}
\frac{u^2+v^2}{2}+\frac{(u^2+v^2)M^{-2}}{\gamma-1}= \frac{u_\infty^2}{2}+\frac{M_{\infty}^{-2}}{\gamma-1}.
\label{eqn:BerHyper}
\end{equation}
Without loss of generality, we may choose $u_\infty=1$ by scaling; otherwise,
we can simply scale: $U\rightarrow u_\infty^{-1}U$, in system (\ref{eqn:orig}) and (\ref{eqn:BerHyper}).
With the fixed $u_\infty=1$, $M_{\infty}\to\infty$ is equivalent to $p_{\infty}\to0$, or $c_{\infty}\to0$.

System (\ref{eqn:orig}) can be written in the form:
\begin{equation}
\partial_xW(U)+\partial_yH(U)=E(U,y)
\label{eqn:Conseve1}
\end{equation}
with $U=(u,v)^\top$, where
\begin{align*}
W(U)=(\rho u,v)^\top,\quad\,\, H(U)=(\rho v,-u)^\top,\quad\,\,
E(U,y)=(-\frac{\rho v}{y},0)^\top,
\end{align*}
and $\rho$ is a function of $U$ determined by the Bernoulli law (\ref{eqn:Bernoulli}).\par
When $\rho>0$ and $u>c$, $U$ can also be presented by $W(U)=(\rho u,v)^\top$, \textit{i.e.}, $U=U(W)$, by the implicit function theorem, since the Jacobian:
\begin{align*}
\det\big(\nabla_UW(U)\big)=-\frac{\rho}{c^2}(u^2-c^2)<0.
\end{align*}
Regarding \textit{x} as the \textit{time} variable, (\ref{eqn:Conseve1}) can be written as
\begin{equation}
\partial_xW+\partial_yH(U(W))=E(U(W),y).
\label{eqn:Conseve2}
\end{equation}
Therefore, system (\ref{eqn:orig})--(\ref{eqn:Bernoulli}) becomes a hyperbolic system of conservation laws
with source terms of form (\ref{eqn:Conseve2}).
Such nonhomogeneous hyperbolic systems of conservation laws also arise naturally in other problems
from many important applications, which exhibit rich phenomena;
for example, see \cite{Chen&Feldman2018,Chen&Liu1994,Chen&Wagner2003,Chen&Wang2022,Courant1948,Dafermos2016}
and the references cited therein.\par

Throughout this paper, the following conditions are assumed:

\medskip
\Assumption\label{Asum:pres} $\,p^b(x)>0$ for $x>0$,
\begin{align*}
p^{b}(x)=p_{0}\qquad \mbox{for $x\in[0,x_{0}]$},
\end{align*}
where
$x_{0}>0$, $p_0\in (0, p^*)$ for some $p^*>0$ to be determined by $\gamma>1$,
and
\begin{align*}
(p^{b})'_{+}(x)=\lim_{t\rightarrow x+}\frac{p^{b}(t)-p^b(x)}{t-x}\in\text{BV}([0,\infty)).
\end{align*}

\Assumption\label{Asum:super} The velocity of the incoming flow $U_{\infty}=(1,0)^\top$ is supersonic: $M_{\infty}>1$.\\

Given a perturbed pressure distribution $p^b(x)$ on the cone surface,
the problem is axisymmetric with respect to the $x-$axis.
Thus, it suffices to analyze the problem in the half-space $\{y\leq0\}$.
Then the inverse problem is to find the generating curve $y=b(x)\le 0$ of the cone surface
and a global solution in the domain:
\begin{equation}\label{domain}
\Omega=\big\{(x,y)\,:\,x\geq0,y<b(x)\big\}
\end{equation}
with its upper boundary:
\begin{equation}\label{bdry}	
\Gamma=\big\{(x,y)\,:\,x\geq0,y= b(x)\big\}
\end{equation}
such that
\begin{equation}
U\cdot\textbf{n}|_{\Gamma}=0,
\label{eqn:Boundary}
\end{equation}
where
$\textbf{n}=\textbf{n}\big(x,b(x)\big)
=\frac{(-b'(x),1)^\top}{\sqrt{1+(b'(x))^2}}
$
is the corresponding outer normal vector to $\Gamma$ at a differentiable point $(x,b(x))\in \Gamma$.

\smallskip
With this setup, the inverse stability problem can be formulated into
the following initial-boundary value problem (IBVP)
for system (\ref{eqn:Conseve1}):

\medskip
\textbf{Cauchy Condition:}
\begin{equation}
U|_{x=0}=U_{\infty}:=
(1,0)^\top,
\label{eqn:Cauchyda}
\end{equation}	

\textbf{Boundary Condition:}
\begin{equation}
p\big(x,b(x)\big)=p^b(x).
\label{eqn:Boundarycon}
\end{equation}

\medskip
We first introduce the notion of entropy solutions for
problem (\ref{eqn:Conseve2})--(\ref{eqn:Boundarycon}).

\medskip
\Definition\label{def:entrosol}(Entropy Solutions).
Consider the inverse problem (\ref{eqn:Conseve2})--(\ref{eqn:Boundarycon}).
A function $b(x)\in\text{Lip}([0,\infty))$ is called a generating curve $\Gamma$ of the cone surface
as defined in \eqref{bdry},
and a vector function $U=(u,v)^\top\in(\text{BV}_{\text{loc}}\cap\text{L}^{\infty})(\Omega)$
with $\Omega$ defined in \eqref{domain}
is called an entropy solution of (\ref{eqn:Conseve2})--(\ref{eqn:Boundarycon})
if they satisfy the following conditions{\rm :}
\begin{itemize}
\item[(\expandafter{\romannumeral1})] For any test function $\phi\in \text{C}_{0}^1(\mathbb{R}^2;\mathbb{R})$
and $\psi\in\text{C}_{0}^1(\Omega;\mathbb{R})$,
\begin{equation}
\begin{aligned}
&\iint_{\Omega}\big(\rho u\phi_x+\rho v\phi_y-\dfrac{\rho v}{y}\phi\big)\,\dd x \dd y
+\int_{-\infty}^{0}\phi(0,y)\rho_{\infty}u_{\infty}\,\dd y=0,\\
&\iint_{\Omega}\left(v\psi_x-u\psi_y\right) \dd x \dd y=0,
\end{aligned}	
\label{eqn:DEFofwea}
\end{equation}

\item[(\expandafter{\romannumeral2})] For any convex entropy pair $(\mathcal{E},\mathcal{Q})$
with respect to $W$ of (\ref{eqn:Conseve2}), \textit{i.e.}, $\nabla^2\mathcal{E}(W)\geq0$ and $\nabla\mathcal{Q}(W)=\nabla\mathcal{E}(W)\nabla H(U(W))$,
\begin{equation}	
\begin{aligned}
&\iint_{\Omega}\big(\mathcal{E}(W(U))\varphi_x+\mathcal{Q}(W(U))\varphi_y
+\nabla_W\mathcal{E}(W(U))E(U,y)\varphi\big)\,\dd x \dd y\\
&+\int_{-\infty}^{0}\mathcal{E}(W(U_{\infty}))\varphi(0,y)\,\dd y\geq0
\end{aligned}	
\label{eqn:DEFofentro}
\end{equation}
for any $\varphi\in\text{C}_{0}^1(\Omega;\mathbb{R})$ with
$\varphi\geq0$.
\end{itemize}

\medskip
\Remark For the potential flow, the Bernoulli law (\ref{eqn:Bernoulli}) gives
\begin{align*}
\frac{M^2}{2}+\frac{1}{\gamma-1}=\frac{B_\infty}{c^2}
\end{align*}
for $B_\infty=\frac{u_\infty^2}{2}+\frac{c_\infty^2}{\gamma-1}$, $c^2=\gamma\rho^{\gamma-1}$, and $p=\rho^\gamma$.
Then the assumptions on pressure $p^b$ can be reduced to the equivalent ones on the Mach number $M_b$ on the unknown boundary $\Gamma$.

\medskip
\medskip
\noindent
\textbf{Main Theorem} (Existence and stability). Let (\ref{Asum:pres})--(\ref{Asum:super}) hold,
and let $1<\gamma<3$ and
$$
0<p_{0}<p^*:=\Big((\sqrt{\gamma+7}-\sqrt{\gamma-1})
\sqrt{\frac{\gamma-1}{16\gamma}} \Big)^\frac{2\gamma}{\gamma-1}.
$$
Assume that $M_{\infty}$ is sufficiently large and $\epsilon_0$ is sufficiently small such that
\begin{equation}	
\begin{aligned}
\text{T.V.}\,p^{b}<\epsilon_0,
\end{aligned}	
\label{eqn:Condpre}
\end{equation}
then the following statements hold{\rm :}
\newcounter{Main}
\begin{list}{(\roman{Main})}{\usecounter{Main}}
\item \textit{Global existence}:
IBVP (\ref{eqn:Conseve1})--(\ref{eqn:Boundarycon}) determines
a boundary $y=b(x)=\int_{0}^{x}b'(t)\,\dd t$ with $b'(x)\in $BV$(\mathbb{R}_{+})$,
which is a small perturbation of the generating curve of the straight-sided cone surface: $y=b_{0}x$,
and admits a global entropy solution $U(x,y)$ with bounded total variation:
\begin{equation}	
\sup_{x>0}\text{T.V.}\{U(x,y)\, :\, -\infty<y<b(x)\}<\infty
\label{eqn:bdTV}
\end{equation}
in the sense of Definition \ref{def:entrosol},
containing a strong leading shock-front $y=\chi(x)=\int_{0}^xs(t)\dd t$
with $s(x)\in$BV$(\mathbb{R}_{+})$ which is a small perturbation of
the strong straight-sided conical shock-front $y=s_{0}x$,
so that the solution between the leading shock-front and the cone surface
is a small perturbation of the background self-similar solution of the straight-sided cone case,
where $s_0$ denotes the slope of the corresponding straight-sided shock-front,
and $b_{0}$ is the slope of the generating curve of the straight-sided cone surface.

\item \textit{Asymptotic behavior}:
For the entropy solution $U(x,y)$,
\begin{equation}	
\lim_{x\rightarrow\infty}\sup\left\{\big|U_{\vartheta}(x,y)-\tilde{U}
(\sigma;s_{\infty},G(s_{\infty}))\big|\,:\,\, \chi_{\vartheta}(x)<y<b_{\vartheta}(x)\right\}=0
\label{eqn:asymU}
\end{equation}
with $\tilde{U}(\sigma;s_{\infty},G(s_{\infty}))$
satisfying $\tilde{U}(s_{\infty};s_{\infty},G(s_{\infty}))=G(s_{\infty})$,
\begin{equation}	
\begin{aligned}
&\tilde{U}(b_{\infty}';s_{\infty},G(s_{\infty}))\cdot(-b_{\infty}',1)=0,\\
&\dfrac{1}{2}\big|\tilde{U}(b_{\infty}';s_{\infty},G(s_{\infty}))\big|^2
+\dfrac{\gamma (p^{b}_{\infty})^{\frac{\gamma-1}{\gamma}}}{\gamma-1}
=\dfrac{1}{2}+\dfrac{\gamma p_{\infty}^{\frac{\gamma-1}{\gamma}}}{\gamma-1},
\end{aligned}	
\label{eqn:asymnota}
\end{equation}
where
\begin{equation}	
\begin{aligned}
p^{b}_{\infty}=\lim_{x\rightarrow \infty}p^{b}(x),\quad\,\,
s_{\infty}=\lim_{x\rightarrow \infty}s_{\vartheta}(x),\quad\,\,
b'_{\infty}=\lim_{x\rightarrow \infty}(b_{\vartheta})'_{+}(x),
\end{aligned}	
\label{eqn:asymlim}
\end{equation}
$\tilde{U}(\sigma;s,G(s))$ is the state of the self-similar solution,
and $G(s)$ denotes the state connected to state $U_{\infty}$ by the strong leading shock-front
of speed $s$.
\end{list}

\smallskip
During the last forty years,
the shock stability problem has been studied
for the perturbed cones with small perturbations of the straight-sided cone.
For polytropic potential flow near the cone vertex,
the local existence of piecewise smooth solutions was established in \cite{Chen2000,Chen2001}
for both symmetrically perturbed cone and pointed body, respectively.
Lien-Liu in \cite{Liu1999} first analyzed the global existence of weak solutions
via a modified Glimm scheme for the uniform supersonic isentropic Euler flow past
over a piecewise straight-side cone, provided that the cone has a small opening
angle (the initial strength of the shock-front is relatively weak)
and the Mach number of the incoming flow is sufficiently large.
Later on, Wang-Zhang considered in \cite{Wang2009} for supersonic potential flow
for the adiabatic exponent $\gamma\in(1,3)$ over a symmetric Lipschitz cone
with an arbitrary opening angle less than the critical angle
and constructed
global weak solutions
that are small perturbations of the self-similar solution,
given that the total variation of the slopes of the perturbed generating curves
of the cone is sufficiently small and the Mach number of the incoming flow
is sufficiently large.
In addition, for the isothermal flows (\textit{i.e.}, $\gamma=1$),
Chen-Kuang-Zhang in \cite{Kuang2021} made full use of delicate expansions up to second-order as the Mach number
of the incoming flow goes to infinity and provided a complete proof of the global existence
and asymptotic behavior of conical shock-front solutions in BV
when the isothermal flows past Lipschitz perturbed cones that are small perturbations
of the straight-sided one.

When the surface of the perturbed cone is smooth, using the weighted energy methods,
Chen-Xin-Yin established the global existence of piecewise smooth solutions in \cite{Xin2002}.
They considered a 3-D axisymmetric potential flow past a symmetrically perturbed cone
under the assumption that the attached angle is sufficiently small and the Mach number of the incoming flow is sufficiently large.
This result was also extended to the M-D potential flow case;  see \cite{Li2014} for more details.
Under a certain boundary condition on the cone surface, the global existence of the M-D conical shock
solutions was obtained in \cite{Xin2006} when the uniform supersonic incoming flow with large Mach number
past a generally curved sharp cone.
Meanwhile, using a delicate expansion of the background solution, Cui-Yin established the global
existence and stability of a steady conical shock wave in \cite{Yin2007,Yin2009}
for the symmetrically perturbed supersonic flow past an infinitely long cone whose vertex angle
is less than the critical angle. More recently, by constructing new background solutions
that allow the speeds of the incoming flows to approach the limit speed,
the global existence of steady symmetrically conical shock solutions was established
in Hu-Zhang \cite{Hu2019} when a supersonic incoming potential flow hits a symmetrically
perturbed cone with an opening angle less than the critical angle.
We also remark that some pivotal results have been obtained on the stability of M-D transonic shocks
under symmetric perturbations of the straight-sided cones or the straight-sided wedges,
as well as on Radon measure solutions for steady compressible Euler equations of hypersonic-limit conical flows;
see \cite{Xiang2021,Fang2009,Fang2017,Qu2020,Xu2009} and the references cited therein.

Corresponding to these shock stability problems, two types of inverse problems have been considered.
One type is for the problems of determining the shape of the wedge in the planar steady supersonic flow
for the given location of the leading shock front.
This kind of inverse problems and the related inverse piston problems have been considered
by Li-Wang
in \cite{Wang2007,Li2009,Li2006,Li2007,Wang2014,Wang2019}, where the leading shock-front
is assumed to be smooth and the characteristic method is applied to find the piecewise
smooth solution with the leading shock as its only discontinuity; see also \cite{Li2022}.
The other one is for the problems of determining the shape of the wedge or the cone with given
pressure distribution on it in the planar steady supersonic flow ({\it cf}. \cite{Pu2023}) or axisymmetric conical
steady supersonic flow.
Though various numerical methods and the linearized method have been proposed to deal
with this type of problems, there seems no rigorous result on the existence of solutions
to such inverse problems for steady supersonic flow past a cone.

\vspace{2pt}
In this paper, we develop a modified Glimm scheme to establish the global existence
and the asymptotic behavior of conical shock-front solutions
of the inverse problem in $BV$ in the flow direction,
when the isentropic flows past cones with given pressure distributions on their surfaces,
which are small perturbations of a constant pressure less than the critical value.
Mathematically, our problem can be formulated as a free boundary problem governed
by 2-D steady isentropic irrotational Euler flows with geometric structure.

There are two main difficulties in solving this problem: One of them is the singularity generated
by the geometric source term, and the other is that, compared to the shock stability problem
for supersonic flows past a cone, the generating curve of the cone is unknown.
For supersonic flows past an axisymmetric cone with the given generating curve,
a modified Glimm scheme developed by Lien-Liu in \cite{Liu1999} is used to
construct approximate solutions (see also \cite{Wang2009,Kuang2021}).
In the previous construction, in order to incorporate with the geometric source term
and the boundary condition on the approximate generating curve,
the center $(x_0,0)$ of the self-similar variable $\sigma=\frac{x-x_0}{y}$
is defined to be the intersection of the $x$-axis and the line
on which the current approximate generating curve (a line segment of a polyline)
lies, and the center is changed according to the random choice at each step
when the ordinary differential equations (\ref{eqn:selfsimilar}) are solved.
As a result, the approximate solution on the approximate generating curve
is a piecewise constant vector-valued function that satisfies the boundary condition everywhere.
However, in the inverse problem under consideration in this paper,
the generating curve of the cone is to be determined, {\it apriori} unknown,
so that the approach in Lien-Liu ({\it c.f.}\cite{Liu1999}) could not apply directly.\par

To overcome the new difficulties, we first fix the center of the self-similar
variable to be the origin when solving the differential equations (\ref{eqn:selfsimilar})
and then develop a modified Glimm scheme to construct approximate solutions
$U_{\Delta x,\vartheta}(x,y)$ via the self-similar solutions as building blocks
in order to incorporate the geometric source term.
In our construction, the grid points are fixed at the beginning,
which are the intersections of lines $x=x_h$, $h\in\mathbb{N}$, and the rays
issuing from the origin (the vertex point of the cone).
Consequently, this construction allows us to find new terms $\theta_{b}(h)$ to control
the increasing part of the Glimm type functional near the approximate
boundary (see Lemma \ref{Lem:Estchabdry}),
while it brings us an extra error so that the boundary conditions on the approximate boundary
are no longer satisfied everywhere, but are satisfied at the initial point of each approximate
boundary at each step.
Nevertheless, in Proposition \ref{Prop:Conver3}, we are able to prove that this error goes
to zero as the grid size $\Delta x$ tends to zero. \par

Furthermore, we make careful asymptotic expansions of the self-similar solutions with respect to
$M_{\infty}^{-1}$. We then make full use of the asymptotic expansion analysis of the
background solutions with respect to $M_{\infty}^{-1}$ to calculate the reflection
coefficients $K_{r,1}$, $K_{w,2}$, $K_s$, and $\mu_{w,2}$ of the weak waves reflected
from both the boundary and the strong leading shock,
and the self-similar solutions reflected from the strong leading shock to prove that
\smallskip
\begin{equation*}
\lim_{M_{\infty}\rightarrow\infty}\big(|K_{r,1}||K_{w,2}|+|K_{r,1}||K_s||\mu_{w,2}|\big)<1.
\end{equation*}
Based on this, we choose some appropriate weights, independent of $M_{\infty}$,
in the construction of the Glimm-type functional and show that the functional is
monotonically decreasing.
Then the convergence of the approximate solutions is followed by the standard approach
for the Glimm-type scheme as in \cite{Glimm1965,Lax1973};
see also \cite{Chen2006,Chen2004,Bressan2000,Dafermos2016,Smoller1983}.
Finally, the existence of entropy solutions and the asymptotic behavior of the entropy solutions
are also proved.

\smallskip
The remaining part of this paper is organized as follows:
In \S \ref{Section-Pre}, we give some preliminaries of the homogeneous system (\ref{eqn:orig})
and then study Riemann-type problems in several cases and self-similar solutions
of the unperturbed conic flow.
Also, we calculate the limit states of related quantities as $M_{\infty} \to \infty$.
In \S \ref{Section-Glimm-Appro}, we construct a family of approximate solutions
via a modified Glimm scheme.
In \S \ref{Section-RIemann-Estimate}, we establish some essential interaction estimates
in a small neighborhood in the limit state.
Then, in \S \ref{Section-Decrease}, we define the Glimm-type functional and
show the monotonicity of the Glimm-type functional and, in \S\ref{Section-final},
we prove that there exists a subsequence of approximate solutions converging to the entropy
solution.
Finally, in \S\ref{Section-asymp},
we give the asymptotic behavior of the entropy solution
which, together with the existence theory, leads to our main theorem.	

\bigskip
\section{Riemann Problems and Self-Similar Solutions of the Unperturbed Conic Flow}\label{Section-Pre}

\smallskip
Regarding $x$ as the \emph{time} variable, the simplified system of (\ref{eqn:orig}):
\begin{equation}
\left\{\begin{aligned}
&(\rho u)_x+(\rho v)_y=0, \\
&v_x-u_y=0,
\end{aligned}
\right.
\label{eqn:simorig}
\end{equation}
is strictly hyperbolic with two distinctive eigenvalues:
\begin{align*}
\lambda_{i}=\frac{uv+(-1)^ic\sqrt{u^2+v^2-c^2}}{u^2-c^2},\qquad  i=1,\ 2,
\end{align*}
for $u>c_*$ and $u^2+v^2<q_*^2$, where
\begin{align*}
c_*=\sqrt{\frac{\gamma-1}{\gamma+1}+\frac{2c_{\infty}^2}{\gamma+1}},
\quad\, q_*=\sqrt{1+\frac{2c_{\infty}^2}{\gamma+1}}.
\end{align*}
Denote $q:=\sqrt{u^2+v^2}$ and $\theta:=\arctan\frac{v}{u}$. Then
\begin{align*}
\lambda_{i}=\tan(\theta+(-1)^i\theta_{m}),\qquad i=1, 2,
\end{align*}
where
\begin{align*}
\theta_{m}:=\arctan\frac{c}{\sqrt{q^2-c^2}}
\end{align*}
is the Mach angel. A direct computation indicates $\theta_{m}\in(0,\frac{\pi}{2})$.\\

Next, we introduce the following lemma, whose proof can be found in \cite{Wang2009}.

\smallskip
\Lemma\label{Lem:eigenvalue} For $u>c_*$ and $q<q_*$,
\begin{align*}
(-\lambda_{i},1)\cdot(\frac{\partial\lambda_{i}}{\partial u},\frac{\partial\lambda_{i}}{\partial v})
= \frac{\gamma+1}{2\sqrt{q^2-c^2}}\sec^3(\theta+(-1)^i\theta_{m}),\quad\, i=1, 2.
\end{align*}

Then, setting
\begin{align*}
r_{i}(U)=\frac{(-\lambda_{i}(U),1)}{(-\lambda_{i}(U),1)\cdot\nabla\lambda_{i}(U)},\quad\, i=1, 2,
\end{align*}
we see that $r_{i}(U)\cdot\nabla\lambda_{i}(U)=1$ for $i=1, 2$.

Denote the supersonic part of the shock polar by
$$
S((u_{\infty},0))=\big\{(\bar{u},\bar{v})\, :\, \bar{c}^2< \bar{u}^2+\bar{v}^2\leq 1\big\},
$$
where $(\bar{u},\bar{v})$ satisfies the Rankine–Hugoniot condition:
\begin{equation}
\left\{\begin{aligned}
&\bar{\rho}(\bar{u}s-\bar{v})=\rho_{\infty}s,\\
&\bar{u}+\bar{v}s=1,
\end{aligned}
\right.
\label{eqn:RHCondition}
\end{equation}
with 	
$\frac{\bar{u}^2+\bar{v}^2}{2}+\frac{\gamma\bar{\rho}^{\gamma-1}}{\gamma-1}
=\frac{1}{2}+\frac{c_{\infty}^2}{\gamma-1}$.
Let
$$
S_{1}^-((u_{\infty},0))=\big\{(\bar{u},\bar{v})\,:\, (\bar{u},\bar{v})\in S((u_{\infty},0)),\,\bar{v}<0\big\}
$$
be the part of shock polar corresponding to the $\lambda_{1}$--characteristic field.
Similar to \cite{Keyfitz1991,Zhang1999,Zhang2003}, we can parameterize the shock polar
$S_{1}^-((u_{\infty},0))$ by a C$^2$--function:
$$
G:s\mapsto G(s;U_{\infty}) \qquad\,\,\mbox{with $U_{\infty}=(1,0)$},
$$
where $G(s;U_{\infty})$ is a supersonic state connected with $U_{\infty}$ by a shock of speed $s$.
For simplicity, we write $G(s;U_{\infty})$ as $G(s)$ and use $\bar{u}(s)$ and $\bar{v}(s)$ to
denote the components of $G(s)$,
that is, $G(s)=(\bar{u}(s),\bar{v}(s))^\top$.
Then we have the following property for $S_{1}^-((u_{\infty},0))$ ({\it cf}. \cite{Chen2020}).

\smallskip
\Lemma\label{Lem:monoshockpole} For $s<\lambda_{1}(U_{\infty})$, $\,\bar{\rho}(s)$ is a strictly monotonically
decreasing function of $s$, and  $\bar{u}(s)$ is a strictly monotonically increasing function of $s$.

\begin{figure}
\begin{center}
\includegraphics{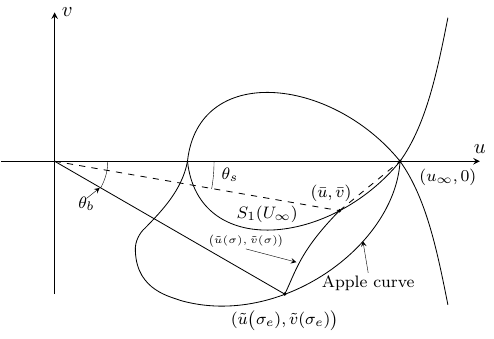}
\end{center}
\caption{Apple curve}
\label{fig:Appcur}
\end{figure}

\smallskip
As in \cite{Courant1948}, let $\sigma=\frac{y}{x}$. Then the equations in (\ref{eqn:orig}) become
\begin{equation}
\left\{\begin{aligned}
&\big(1-\frac{u^2}{c^2}\big)\sigma^2 u_{\sigma} -\frac{2uv}{c^2}\sigma^2 v_{\sigma}
-\big(1-\frac{v^2}{c^2}\big)\sigma v_{\sigma}- v=0, \\
&u_{\sigma}+\sigma v_{\sigma}=0.
\end{aligned}
\right.
\label{eqn:selfsimilar}
\end{equation}
or, equivalently,
\begin{equation}
\left\{\begin{aligned}
&u_\sigma=\frac{c^2v}{(1+\sigma^2)c^2-(v-\sigma u)^2},\\
&v_\sigma=-\frac{c^2v}{\sigma\big((1+\sigma^2)c^2-(v-\sigma u)^2\big)},\\
&\rho_\sigma=\frac{\rho v(v-\sigma u)}{\sigma\big((1+\sigma^2)c^2-(v-\sigma u)^2\big)}.
\end{aligned}
\right.
\label{eqn:eqselfsimilar}
\end{equation}
Given a constant state $(\bar{u},\bar{v})=G(s)$ on $S_{1}^-((u_{\infty},0))$,
there exists a local solution $\tilde{U}(\sigma;s,G(s))=(\tilde{u}(\sigma;s),\tilde{v}(\sigma;s))$
of system (\ref{eqn:selfsimilar}) with initial data:
\begin{align*}
(\tilde{u}(s;s),\,\tilde{v}(s;s))=(\bar{u},\,\bar{v}).
\end{align*}
This solution can be extended to an end-point $(\tilde{u}(\sigma_e;s),\,\tilde{v}(\sigma_e;s))$
with $\frac{\tilde{v}(\sigma_e;s)}{\tilde{u}(\sigma_e;s)}=\sigma_e$.
As $(\bar{u},\bar{v})$ varies on $S_{1}^-((u_{\infty},0))$,
the collection of these end-states forms an \textit{apple curve} (Fig. \ref{fig:Appcur})
through $U_{\infty}$; see \cite{Courant1948}.
For these solutions, we have following properties, whose proof can be found in \cite{Wang2009}.

\medskip
\Lemma\label{Lem:posinume}
For $\tilde{u}(s;s)>\tilde{c}(s;s)$ and $\sigma\in(s,\sigma_e)$,
then
$\tilde{u}(\sigma;s)\sigma-\tilde{v}(\sigma;s)<0$,
\begin{align*}
&\frac{\partial \tilde{u}}{\partial \sigma}<0,\qquad \frac{\partial \tilde{v}}{\partial \sigma}<0,\\
&\tilde{c}(\sigma;s)-\frac{\tilde{v}(\sigma;s)-\sigma \tilde{u}(\sigma;s)}{\sqrt{1+\sigma^2}}
>\tilde{c}(s;s)-\frac{\tilde{v}(s;s)-s\tilde{u}(s;s)}{\sqrt{1+s^2}}>0,
\end{align*}
with
$\frac{\tilde{u}^2+\tilde{v}^2}{2}+\frac{\tilde{c}^2}{\gamma-1}= \frac{1}{2}+\frac{c_{\infty}^2}{\gamma-1}$.

\smallskip
Thus, we obtain the following estimate of the self-similar solution $(\tilde{u}(\sigma;s),\,\tilde{v}(\sigma;s))$.

\smallskip
\Lemma\label{Lem:proficomp} For $\tilde{u}(s;s)>\tilde{c}(s;s)$ and $\sigma\in(s,\sigma_e)$,
\begin{align*}
\frac{1}{1+s^2}< \tilde{u}(\sigma;s)<\tilde{u}(s;s),\qquad 	
\tilde{c}(s;s)<\tilde{c}(\sigma;s)<\tilde{c}(\sigma_e;s)<\sqrt{\frac{(\gamma-1)s^2}{2(1+s^2)}+\frac{1}{M_{\infty}^2}}.
\end{align*}

To obtain the asymptotic expansion of the self-similar solution,
we need the following properties of the shock polar. Set
\begin{equation}
p^*=\left(\big(\sqrt{\gamma+7}-\sqrt{\gamma-1}\big)
\sqrt{\frac{\gamma-1}{16\gamma}} \right)^\frac{2\gamma}{\gamma-1}.
\label{eqn:critp}
\end{equation}	

\Lemma \label{Lem:limshockpol}
Let $1<\gamma<3$ and $p_{0}\in(0,p^*)$. For $M_{\infty}$ large enough, the equations:
\begin{equation}
\left\{\begin{aligned}
&\rho_{0}(u^{\sharp}s^{\sharp}-v^{\sharp})=\rho_{\infty}s^{\sharp},\\
&u^{\sharp}+v^{\sharp} s^{\sharp}=1,\\
&\frac{(u^{\sharp})^2+(v^{\sharp})^2}{2}+\frac{c_{0}^2}{\gamma-1}=\frac{1}{2}+\frac{c_{\infty}^2}{\gamma-1}
\end{aligned}
\right.
\label{eqn:shockpol}
\end{equation}
have a unique solution $(u^{\sharp},v^{\sharp},s^{\sharp})$
with $s^{\sharp}<0$, where
$\rho_{\infty}=p_{\infty}^{\frac{1}{\gamma}}$,  $\rho_{0}=p_{0}^{\frac{1}{\gamma}}$, and $c_{0}=\sqrt{\gamma} p_{0}^{\frac{\gamma-1}{2\gamma}}$, such that
\begin{equation}
\begin{aligned}
\lim_{M_{\infty}\rightarrow \infty}u^{\sharp}&=\lim_{M_{\infty}\rightarrow\infty}u_{a}=1-\frac{2c_{0}^2}{\gamma-1},\\
\lim_{M_{\infty}\rightarrow \infty}v^{\sharp}&=\lim_{M_{\infty}\rightarrow\infty}v_{a}
=-\sqrt{\frac{2c_{0}^2}{\gamma-1-2c_{0}^2}}\Big(1-\frac{2c_{0}^2}{\gamma-1}\Big),\\
\lim_{M_{\infty}\rightarrow\infty} c_{a}^2&= c_{0}^2.
\end{aligned}
\label{eqn:shockpolim}
\end{equation}
where
\begin{equation}
u_{a}:=\frac{1}{1+(s^{\sharp})^2},\quad  v_{a}:=\frac{s^{\sharp}}{1+(s^{\sharp})^2},
\quad c_{a}:=\sqrt{\frac{(\gamma-1)(s^{\sharp})^2}{2(1+(s^{\sharp})^2)}+\frac{1}{M_{\infty}^2}},
\label{eqn:shockpolimset}
\end{equation}

\noindent\textbf{Proof.} From the first two equations of (\ref{eqn:shockpol}), we have
\begin{equation}
u^{\sharp}=\frac{\rho_{0}+\rho_{\infty}(s^{\sharp})^2}{\rho_{0}(1+s^2)},
\qquad v^{\sharp}=\frac{(\rho_{0}-\rho_{\infty})s^{\sharp}}{\rho_{0}\big(1+(s^{\sharp})^2\big)}.
\label{eqn:shockpolvelocity}
\end{equation}
With the help of the third equation of (\ref{eqn:shockpol}), we have
\begin{align*}
\bigg( \frac{\rho_{0}+\rho_{\infty}(s^{\sharp})^2}{\rho_{0}\big(1+(s^{\sharp})^2\big)}\bigg)^2
+\bigg(\frac{(\rho_{0}-\rho_{\infty})s^{\sharp}}{\rho_{0}\big(1+(s^{\sharp})^2\big)} \bigg)^2
=1+\frac{2(c_{\infty}^2-c_{0}^2)}{\gamma-1},
\end{align*}
which gives
\begin{align*}
(s^{\sharp})^2&=\frac{2(c_{0}^2-M_{\infty}^{-2})\rho_{0}^2}
{(\gamma-1)\big(\rho_{0}^2-(\gamma M_{\infty}^2)^{-\frac{2}{\gamma-1}}\big)-2(c_{0}^2-M_{\infty}^{-2})\rho_{0}^2}\\
&=\frac{2c_{0}^2}{\gamma-1-2c_{0}^2}-\frac{2(\gamma-1)}{(\gamma-1-2c_{0}^2)^2}\, M_{\infty}^{-2}
+O(M_{\infty}^{-\frac{4}{\gamma-1}})\qquad \text{as $M_{\infty}\rightarrow \infty$}.
\end{align*}
Therefore, noting that $s^{\sharp}<0$, we obtain
\begin{align*}
s^{\sharp}=-\sqrt{\frac{2c_{0}^2}{\gamma-1-2c_{0}^2}}
\bigg(1-\frac{\gamma -1}{2c_{0}^2\left(\gamma-1-2c_{0}^2\right)}\, M_{\infty}^{-2}\bigg)
+O(M_{\infty}^{-\frac{4}{\gamma-1}})\qquad \text{as $M_{\infty}\rightarrow \infty$}.
\end{align*}
Substituting the above expansion into (\ref{eqn:shockpolimset})--(\ref{eqn:shockpolvelocity}),
yields (\ref{eqn:shockpolim}).
\hfill{$\square$}

\medskip
\Lemma \label{Lem:limshockpolctrl}
Let $1<\gamma<3$ and $p_{0}\in(0,p^*)$.
For $M_{\infty}$ large enough, there exists $\rho_{d}$ such that
\begin{equation}
\rho_{0}^{\gamma-1}=\frac{\rho_{d}^{\gamma+1}-\rho_{\infty}^{\gamma+1}}{\rho_{d}^2-\rho_{\infty}^2},
\label{eqn:shockpolrho}
\end{equation}
and the equations:
\begin{equation}
\left\{\begin{aligned}
&\rho_{d}(u_{d}s_{d}-v_{d})=\rho_{\infty}s_{d},\\
&u_{d}+v_{d}s_{d}=1,\\
&\frac{u_{d}^2+v_{d}^2}{2}+\frac{c_{d}^2}{\gamma-1}=\frac{1}{2}+\frac{c_{\infty}^2}{\gamma-1},
\end{aligned}
\right.
\label{eqn:shockpolctrl}
\end{equation}
have a unique solution $(u_{d},v_{d},s_{d})$ with $s_{d}<0$ and $c_{d}=\sqrt{\gamma\rho_{d}^{\gamma-1}}$.
Moreover, for
\begin{equation}
u^{\flat}:=\frac{1}{1+s_{d}^2},\quad\,  v^{\flat}:=\frac{s_{d}}{1+s_{d}^2},
\label{eqn:shockpoctrlimset}
\end{equation}
we have
\begin{equation}
\begin{aligned}
\lim_{M_{\infty}\rightarrow\infty} c_{d}^2&=c_{0}^2,\\
\lim_{M_{\infty}\rightarrow\infty}u^{\flat}&=\lim_{M_{\infty}\rightarrow\infty}u_{d}
=1-\frac{2c_{0}^2}{\gamma-1},\\
\lim_{M_{\infty}\rightarrow\infty} v^{\flat}&= \lim_{M_{\infty}\rightarrow\infty} v_{d}
=-\sqrt{\frac{2c_{0}^2}{\gamma-1-2c_{0}^2}}\,\Big(1-\frac{2c_{0}^2}{\gamma-1}\Big).
\end{aligned}
\label{eqn:shockpoctrllim}
\end{equation}

\noindent
\textbf{Proof.} For each $\rho_{0}$, it is direct to find $\rho_{d}$ such that (\ref{eqn:shockpolrho}) holds.
Moreover, $p_{d}\in (0,p^*)$ when $M_{\infty}$ is large enough.
By Lemma \ref{Lem:limshockpol}, we obtain the unique solution of (\ref{eqn:shockpolctrl}):
\begin{equation}
s_{d}=\frac{2(c_{0}^2-c_{\infty}^2)}{\gamma-1-2(c_{0}^2-c_{\infty}^2)},\quad\, u_{d}
=\frac{\rho_{d}+\rho_{\infty}s_{d}^2}{\rho_{d}(1+s_{d}^2)},
\quad\, v_{d}=\frac{(\rho_{d}-\rho_{\infty})s_{d}}{\rho_{d}(1+s_{d}^2)}.
\label{eqn:shockpolvelocityctrl}
\end{equation}
Then we have
\begin{align*}
s_{d}=-\sqrt{\frac{2c_{0}^2}{\gamma-1-2c_{0}^2}}
\Big(1-\frac{\gamma -1}{2c_{0}^2(\gamma-1-2c_{0}^2)}\, M_{\infty}^{-2}\Big)
+O(M_{\infty}^{-4})\qquad \text{as $M_{\infty}\rightarrow \infty$}.
\end{align*}
Substituting the above expansion into (\ref{eqn:shockpoctrlimset})--(\ref{eqn:shockpolvelocityctrl})
yields (\ref{eqn:shockpoctrllim}).\hfill{$\square$}

\medskip
Given $p_{0}\in(0,p^{*})$, we now solve the following problem for $s_{0}<\sigma<b_{0}$:
\begin{equation}\label{eqn:backsol}
\begin{cases}
\sigma^2\big(1-\frac{\tilde{u}^2}{\tilde{c}^2}\big) \tilde{u}_{\sigma}
-\frac{2\tilde{u}\tilde{v}\sigma^2}{c^2}\tilde{v}_{\sigma}
- \big(1-\frac{\tilde{v}^2}{\tilde{c}^2}\big) \tilde{v}_{\sigma}\sigma-\tilde{v}=0,\\
\tilde{u}_{\sigma}+\sigma \tilde{v}_{\sigma}=0,
\end{cases}
\end{equation}
with the boundary conditions:
\begin{equation}\label{eqn:backsol-2}
\begin{aligned}
&	\,\,\,\tilde{\rho}(\tilde{u}s_{0}-\tilde{v})=\rho_{\infty}s_{0},\quad
\tilde{u}+\tilde{v}s_{0}=1
\qquad\,\, &&\mbox{for $\sigma=s_{0}$},\\
&\,\,\,\tilde{\rho}=\rho_{0},\quad
\tilde{v}=b_{0}\tilde{u}
\qquad\,\, &&\mbox{for $\sigma=b_{0}$},
\end{aligned}
\end{equation}
and define $\big(\tilde{u}(\sigma;s_{0}),\tilde{v}(\sigma;s_{0})\big)=(1,0)$
for $\sigma<s_{0}$. Indeed, we have the following lemma.\\

\Lemma\label{Lem:asymsol}
Let $1<\gamma<3$ and $p_{0}\in(0,p^{*})$.
For $M_{\infty}>K_{1}$, problem (\ref{eqn:backsol}) has a unique solution
$((\tilde{u}(\sigma;s_{0}),\,\tilde{v}(\sigma;s_{0}))$ containing a supersonic conical shock-front
issuing from the vertex.
In addition,
\begin{equation}
\begin{aligned}
&\lim_{M_{\infty}\rightarrow\infty}(\sigma, \tilde{c}^2(\sigma;s_{0}))=(\tan\theta_{0}, c_{0}^2),\\		
&\lim_{M_{\infty}\rightarrow\infty}(\tilde{u}(\sigma;s_{0}),\tilde{v}(\sigma;s_{0})
=(\cos^2\theta_{0},\,\sin\theta_{0}\cos\theta_{0}),\\
\end{aligned}
\label{eqn:profilim}
\end{equation}
and
\begin{equation}
\lim_{M_{\infty}\rightarrow\infty}\frac{\tilde{u}(\sigma;s_{0})}{\tilde{c}(\sigma;s_{0})}
=\frac{\gamma-1-2c_{0}^2}{(\gamma-1)c_{0}}>1,\qquad
\cos(\theta_{0}\pm\theta_{m}^{0})>0,
\label{eqn:Machaangle}
\end{equation}
where $\theta_{0}=-\arctan\sqrt{\frac{2c_{0}^2}{\gamma-1-2c_{0}^2}}$ and
$\theta_{m}^{0}=\lim_{M_{\infty}\rightarrow\infty}\theta_{m}$
for $\sigma\in[s_{0},b_{0})$.

\smallskip
\noindent
\textbf{Proof.}
Given $p_{0}\in(0,p^{*})$, by the shooting method as in \cite{Courant1948},
problem (\ref{eqn:backsol}) has a unique solution $\big(\tilde{u}(\sigma;s_{0}),\tilde{v}(\sigma;s_{0})\big)$ with $\tilde{u}(s_{0};s_{0})>\tilde{c}(s_{0};s_{0})$, $\tilde{\rho}(b_{0};s_{0})=p_{0}^{\frac{1}{\gamma}}$, and $\tilde{v}(b_{0};s_{0})=b_{0}\tilde{u}(b_{0};s_{0})$.

We then focus on the asymptotic expansions (\ref{eqn:profilim}). Lemma \ref{Lem:proficomp} indicates that
\begin{align*}
\frac{1}{1+s_{0}^2}<\tilde{u}(\sigma;s_{0})\leq\tilde{u}(s_{0};s_{0}),\qquad	
\tilde{c}(s_{0};s_{0})\leq\tilde{c}(\sigma;s_{0})<c_{0}
<\sqrt{\frac{(\gamma-1)s_{0}^2}{2(1+s_{0}^2)}+\frac{1}{M_{\infty}^2}},
\end{align*}
for $\sigma\in[s_{0},b_{0})$.
Meanwhile, it follows from (\ref{eqn:shockpolrho}) that $\tilde{c}(s_{0};s_{0})>c_{d}$.
Then, due to Lemma \ref{Lem:monoshockpole}, we see that
$u^{\sharp}<\tilde{u}(s_{0};s_{0})<u_{d}$ and $s^{\sharp}<s_{0}<0$.
Therefore, we have
\begin{align*}
c_{d}<\tilde{c}(s_{0};s_{0})\leq\tilde{c}(\sigma;s_{0})<c,\qquad
u_{a}=\frac{1}{1+(s^{\sharp})^2}<\frac{1}{1+s_{0}^2}<&\tilde{u}(\sigma;s_{0})\leq\tilde{u}(s_{0};s_{0})<u_{d}.
\end{align*}
From Lemma \ref{Lem:limshockpol}--\ref{Lem:limshockpolctrl}, we have
\begin{align*}
\lim_{M_{\infty}\rightarrow\infty}\tilde{c}(\sigma;s_{0})=c_{0}^2,
\qquad \lim_{M_{\infty}\rightarrow\infty}\tilde{u}(\sigma;s_{0})=1-\frac{2c_{0}^2}{\gamma-1}.
\end{align*}	
Since $\tilde{v}(\sigma;s_{0})<0$, from the Bernoulli laws,
\begin{align*}
\frac{\tilde{u}^2+\tilde{v}^2}{2}+\frac{\tilde{c}^2}{\gamma-1}= \frac{1}{2}+\frac{c_{\infty}^2}{\gamma-1},\quad\,
\frac{(u^{\sharp})^2+(v^{\sharp})^2}{2}+\frac{c_{0}^2}{\gamma-1}=\frac{1}{2}+\frac{c_{\infty}^2}{\gamma-1},
\end{align*}
we conclude
\begin{align*}
\lim_{M_{\infty}\rightarrow \infty}\tilde{v}(\sigma;s_{0})=&-\sqrt{\frac{2c_{0}^2}{\gamma-1-2c_{0}^2}}\,
\Big(1-\frac{2c_{0}^2}{\gamma-1}\Big).
\end{align*}	
Again, by Lemma \ref{Lem:monoshockpole}, we know that
\begin{align*}
\frac{\tilde{v}(s_{0};s_{0})}{\tilde{u}(s_{0};s_{0})}>b_{0}>\sigma\geq s_{0}>s^{\sharp}.
\end{align*}
Combining all the expansions obtained above, we obtain (\ref{eqn:profilim}).

Furthermore, since
\begin{align*}
\lim_{M_{\infty}\rightarrow\infty}\big( \tilde{u}\sqrt{\tilde{u}^2+\tilde{v}^2-\tilde{c}^2}\big)^2
-( \tilde{v}\tilde{c})^2
=\lim_{M_{\infty}\rightarrow\infty}(\tilde{u}^2-\tilde{c}^2)(\tilde{u}^2+\tilde{v}^2)>0,
\end{align*}
we conclude
\begin{align*}
\cos(\theta_{0}\pm\theta^0_{ma})
=\lim_{M_{\infty}\rightarrow\infty}\cos(\theta\pm\theta_{m})
=\lim_{M_{\infty}\rightarrow\infty}\frac{ \tilde{u}\sqrt{\tilde{u}^2+\tilde{v}^2-\tilde{c}^2}\mp \tilde{v}\tilde{c}} {\tilde{u}^2+\tilde{v}^2}>0.
\end{align*}
This completes the proof.\hfill{$\square$}

\smallskip
Next, for $G(s)$, we have the following expansions, whose proof can be found in \cite{Wang2009}.

\smallskip
\Lemma\label{Lem:asymshockpol} For $G(s)=\big(\bar{u}(s),\bar{v}(s)\big)^\top$,
\begin{align*}
&\lim_{M_{\infty}\rightarrow\infty}(\bar{u}(s_{0}),\,\bar{v}(s_{0})) =(\cos^2\theta_{0},\,\cos\theta_{0}\sin\theta_{0}), \\
&\lim_{M_{\infty}\rightarrow\infty}(\bar{u}_s(s_{0}),\,\bar{u}_s(s_{0}))
=(-\sin2\theta_{0}\cos^2\theta_{0},\,\cos2\theta_{0}\cos^2\theta_{0}),
\end{align*}
where $\bar{u}_s(s)=\frac{{\rm d}\bar{u}(s)}{{\rm d} s}$ and $\bar{v}_s(s)=\frac{{\rm d}\bar{v}(s)}{{\rm d} s}$.

\medskip
Now, we introduce the elementary wave curves of system (\ref{eqn:simorig}).
We denote by $W(p_{0},p_{\infty})$ the curve formed by
$\tilde{U}(\sigma;s_{0})=\big(\tilde{u}(\sigma;s_{0}),\tilde{v}(\sigma;s_{0})\big)^\top$ for $s_{0}<\sigma<b_{0}$,
where $p_{0}$ is the corresponding pressure of the state at the endpoint.
As in \cite{Wang2009} (also {\it cf}. \cite{Zhang1999,Zhang2003}),
we parameterize the elementary $i$-wave curves for system (\ref{eqn:simorig}) in a neighborhood of $W(p_{0},p_{\infty})$:
$$
O_{r}\big(W(p_{0},p_{\infty})\big)=\bigcup\limits_{s_{0}<\sigma<b_{0}}\{U:\ |U-\tilde{U}(\sigma;s_{0})|<r\}
\qquad\mbox{for some $r>0$}
$$
by
\begin{align*}
\alpha_{i}\mapsto\Phi_{i}(\alpha_{i};U)
\end{align*}
with $\Phi_{i}\in \text{C}^2$ and
\begin{align*}
\left.\frac{\partial\Phi_{i}}{\partial \alpha_{i}}\right|_{\alpha_{i}=0}
=r_{i}(U) \qquad\mbox{for $U\in O_{r}(W(p_{0},p_{\infty}))$}, \,\, i=1,2.
\end{align*}

In the sequel, define
\begin{align*}
\Phi(\alpha_{1},\alpha_{2};U)=\Phi_{2}\big(\alpha_{2};\Phi_{1}(\alpha_{1};U)\big).
\end{align*}
Denote $\tilde{U}(\sigma;\sigma_{0},U_{l})$ the solution to the ODE system (\ref{eqn:selfsimilar}) with initial data
\begin{align*}
\tilde{U}|_{\sigma=\sigma_{0}}=U_{l}
\end{align*}
for $U_{l}\in O_{r}\big(W(p_{0},p_{\infty})\big)$. Then, as in \cite{Wang2009}, we have

\medskip
\Lemma\label{Lem:asymselfsi} For $p_{0}\in(0,p^{*})$,
\begin{align*}		
\lim_{M_{\infty}\rightarrow\infty}
\left.\frac{\dd\tilde{U}(\sigma;\sigma_{0})}{\dd\sigma}\right|_{\{\sigma=\sigma_{0},\,U_{l}\in W(p_{0},p_{\infty})\}}
=\big(\sin\theta_{0}\cos^3\theta_{0},-\cos^4\theta_{0}\big)^\top.
\end{align*}

With all the limits given above, we obtain the following lemma, which is essential in wave-interaction estimates.

\medskip
\Lemma \label{Lem:DET} For $U_{l}\in W(p_{0},p_{\infty})$,
\begin{align*}
&\lim_{M_{\infty}\rightarrow\infty}\det \big(r_{1}(U_{l}),r_{2}(U_{l})\big)=\dfrac{4\cos^2(\theta_{0}+\theta_{m}^{0})\cos^2(\theta_{0}-\theta_{m}^{0}) \cos^2\theta_{0}\cos^2\theta^{0}_{m}\sin(2\theta_{m}^{0})}{(\gamma+1)^2},\\
&\lim_{M_{\infty}\rightarrow\infty}\det \big(G'(s_{0}),r_{1}\big(G(s_{0})\big)\big)
=-\dfrac{2 \cos^2(\theta_{0}-\theta_{m}^{0})\cos\theta^{0}_{m}\cos^3\theta_{0}\sin(\theta_{0}+\theta_{m}^{0})}{\gamma+1},\\	
&\lim_{M_{\infty}\rightarrow\infty}\det \big(r_{2}\big(G(s_{0})\big),G'(s_{0})\big)
=\dfrac{2 \cos^2(\theta_{0}+\theta_{m}^{0})\cos\theta^{0}_{m}\cos^3\theta_{0}\sin(\theta_{0}-\theta_{m}^{0})}{\gamma+1},\\
&\lim_{M_{\infty}\rightarrow\infty}\det \Big(\frac{\dd\tilde{U}(\sigma;\sigma_{0},G(s_{0}))}{\dd\sigma},G'(s_{0})\Big)
=- \cos^5\theta_{0}\sin\theta_{0},\\
&\lim_{M_{\infty}\rightarrow\infty}\det \Big(r_{2}(U_{l}),\frac{\dd\tilde{U}(\sigma;\sigma_{0},U_{l})}{\dd\sigma}\Big)
=\frac{2}{\gamma+1} \cos^4\theta_{0}\cos^2(\theta_{0}+\theta_m^{0})\cos\theta_m^{0}\sin\theta_m^{0},\\
&\lim_{M_{\infty}\rightarrow\infty}\det \Big(r_{1}(U_{l}),\frac{\dd\tilde{U}(\sigma;\sigma_{0},U_{l})}{\dd\sigma}\Big)
=-\frac{2}{\gamma+1} \cos^4\theta_{0}\cos^2(\theta_{0}-\theta_m^{0})\cos\theta_m^{0}\sin\theta_m^{0}.
\end{align*}	

Furthermore, we have the following propositions,
which will be used in the construction of building blocks of our approximate solutions.

\medskip
\Proposition\label{Prop:RiemannP}
For $M_{\infty}$ sufficiently large, there exists $\varepsilon_{1}>0$ such that,
for any $U_{r}$ and $U_{l}$ lie in $O_{\varepsilon_{1}}\big(W(p_{0},p_{\infty})\big)$,
the Riemann problem (\ref{eqn:simorig}) with initial data
\begin{equation}
U|_{x=\bar{x}}=\left\{\begin{aligned}
&U_{r}&\,\,\, \text{for }y>\bar{y},\\
&U_{l}&\,\,\, \text{for }y<\bar{y},
\end{aligned}
\right.
\label{eqn:RPIni}
\end{equation}
admits a unique admissible solution consisting of at most two elementary waves $\alpha_{1}$
for the $1$-characteristic field and $\alpha_{2}$ for the $2$-characteristic field.
Moreover, states $U_{l}$ and $U_{r}$ are connected by
\begin{align*}
U_{r}=\Phi(\alpha_{1},\alpha_{2};U_{l}).
\end{align*}

\smallskip	
\Proposition\label{Prop:BdryP}
For $M_{\infty}$ sufficiently large, there exists $\varepsilon_{2}>0$ such that,
for any $U_{l}\in O_{\varepsilon_{2}}\big( W(p_{0},p_{\infty})\big)$ and
$p_{1}$, $p_{2}\in O_{\varepsilon_{2}}(p_{0})$, there is $\delta_{1}$ solving the equation:
\begin{equation}		
\dfrac{1}{2}|\Phi(\delta_{1},0;U_{l})|^2+\dfrac{\gamma }{\gamma-1}p_{2}^{\frac{\gamma-1}{\gamma}}
=\dfrac{1}{2}|U_{l}|^2+\dfrac{\gamma }{\gamma-1}p_{1}^{\frac{\gamma-1}{\gamma}}.
\label{eqn:BdryP}
\end{equation}

\noindent
\textbf{Proof.}
From (\ref{eqn:BdryP}), we have
\begin{align*}
\lim_{M_{\infty}\rightarrow\infty}\frac{1}{2}\,
\left.\frac{\partial |\Phi(\delta_{1},0;U_{l})|^2}{\partial\delta_{1}}\right|_{\delta_{1}=0}
=U_{l}\cdot r_{1}(U_{l})\neq0.
\end{align*}
By the implicit function theorem, there exists $\delta_{1}$ such that (\ref{eqn:BdryP})
holds, provided $\varepsilon_{2}$ sufficiently small.\hfill{$\square$}\\

\Proposition\label{Prop:RPSshock}
For $M_{\infty}$ sufficiently large, there exists $\varepsilon_{3}>0$ such that,
for any $U_{l}=U_{\infty}$ and
$U_{r}\in O_{\varepsilon_{3}}( W(p_{0},p_{\infty}))\cap O_{\varepsilon_{3}}(G(s_{0}))$,
the Riemann problem (\ref{eqn:simorig}) with initial data (\ref{eqn:RPIni})
admits a unique admissible solution that contains a strong $1$-shock $s_{1}$
and a $2$-weak wave $\beta_{2}$ of the 2-characteristic field.
Moreover, states $U_{l}$ and $U_{r}$ are connected by
\begin{equation}		
U_{r}=\Phi_{2}(\beta_{2};G(s_{1};U_{l})).
\label{eqn:SshockP}
\end{equation}

\noindent
\textbf{Proof.}
It follows from (\ref{eqn:SshockP}) and Lemma \ref{Lem:DET} that
\begin{align*}
\lim_{M_{\infty}\rightarrow\infty}
\det\Big(\left.\frac{\partial \Phi_{2}(\beta_{2};G(s_{1};U_{l}))}{\partial(s_{1},\beta_{2})}\right|_{\{s_{1}=s_{0}, \beta_{2}=0\}}\Big)
=-\lim_{M_{\infty}\rightarrow\infty}\det \big(r_{2}(G(s_{0}),G'(s_{0}))\big)\neq0.
\end{align*}
The existence of the solution of this Riemann problem is ensured by the implicit function theorem
for $\varepsilon_{3}$ sufficiently small.
\hfill{$\square$}

\medskip
To end this section, we introduce the following interaction estimate given
by Glimm \cite{Glimm1965} for weak waves (see also \cite{Smoller1983,Zhang2003,Wang2009}).

\smallskip
\Lemma\label{Lem:weakesti}
Let $U_{l}\in  W(p_{0},p_{\infty})$, $\alpha$, $\beta$, and $\delta$ satisfy
\begin{align*}
\Phi(\delta_{1},\delta_{2};U_{l})=\Phi(\beta_{1},\beta_{2};\Phi(\alpha_{1},\alpha_{2};U_{l})).
\end{align*}
Then
\begin{align*}
\delta=\alpha+\beta+O(1)Q^0(\alpha,\beta),
\end{align*}
where $Q^0(\alpha,\beta)=\sum\{|\alpha_{i}||\beta_{i}|:\text{$\alpha_{i}$ and $\beta_{j}$ approach}\}$,
and $O(1)$ depends continuously on $M_{\infty}<\infty$.

\section{Approximate Solutions}\label{Section-Glimm-Appro}	
In this section, we construct approximate solutions for system (\ref{eqn:Conseve1})
with (\ref{eqn:Boundary})--(\ref{eqn:Cauchyda}) by a modified Glimm scheme. Compared to the modified Glimm scheme
developed in \cite{Wang2009,Liu1999,Kuang2021}, in our construction, the grid points are fixed at the very beginning,
which are independent of the approximate solution and the random choice.\par

Given $\epsilon>0$ and $\Delta x>0$, there exist piece-wise constant functions $p^{b}_{\Delta x}$ such that
\begin{align*}
\text{T.V.}~p^{b}_{\Delta x}(\cdot)\leq\text{T.V.}~p^{b}(\cdot),
\qquad \|p^{b}_{\Delta x}-p^{b}\|_{\textbf{L}^\infty}\leq \epsilon,
\end{align*}
where
\begin{align*}
p^{b}_{\Delta x}(x)=\begin{cases}
p^{b}_{\Delta x,0}=p_{0} \,\,\,&\text{ for $x\in[0,x_{0})$},\\
p^{b}_{\Delta x,h+1} \,\,\,&\text{ for $x\in[x_{h},x_{h+1})$ and $h\in\mathbb{N}$},
\end{cases}
\end{align*}
with $p^{b}_{\Delta x,h+1}$ being constants on the corresponding intervals and $x_{h}=x_{0}+h\Delta x$ for
$h\in\mathbb{N}$.
Then, from Lemma \ref{Lem:asymsol}, for $p^{b}_{\Delta x,0}=p_{0}$, there exists
$\big(\tilde{u}(\sigma;s_{0}),\tilde{v}(\sigma;s_{0})\big)$ such that $\tilde{p}(b_{0};s_{0})=p_{0}$
and $\tilde{v}(b_{0};s_{0})=b_{0}\tilde{u}(b_{0};s_{0})$.

\medskip
We now define the difference scheme.
Choose $\vartheta=(\vartheta_{0},\vartheta_{1},\vartheta_{2},\dots,\vartheta_{h},\dots)$ randomly in $[0,1)$.
For $0< x< x_{0}$, let
\begin{align*}
b_{\Delta x,\vartheta}(x)=b_{0}x,\qquad\,\, \chi_{\Delta x,\vartheta}(x)=s_{0}x.
\end{align*}
We denote $\Gamma_{\Delta x,\vartheta,0}=\big\{(x,b_{\Delta x,\vartheta}(x))\,:\, 0\le x< x_{0}\big\}$,
$S_{\Delta x,\vartheta,0}=\big\{(x,\chi_{\Delta x,\vartheta}(x))\,:\, 0\le x< x_{0}\big\}$,
and $\Omega_{\Delta x,\vartheta,0}=\big\{(x,y)\,:\, y<b_{\Delta x,\vartheta}(x), 0\le x< x_{0}\big\}$.
In region $\Omega_{\Delta x,\vartheta,0}$, we then define
\begin{align*}
U_{\Delta x,\vartheta}(x,y)=
\begin{cases}
(u_{\Delta x,\vartheta}(x,y),v_{\Delta x,\vartheta}(x,y))^\top
\triangleq(\tilde{u}(\sigma;s_{0}),\tilde{v}(\sigma;s_{0}))^\top\,\,\, &\text{for $\frac{y}{x}=\sigma\in(s_{0},b_{0})$},\\
U_{\infty} \quad &\text{for $\frac{y}{x}=\sigma<s_{0}$},
\end{cases}
\end{align*}
and, on boundary $\Gamma_{\Delta x,\vartheta,0}$, we set
\begin{align*}
U_{\Delta x,\vartheta}(x, b_{\Delta x,\vartheta}(x))=U_{\Delta x,\vartheta}^b(x)=(u_{\Delta x,\vartheta}^b(x),v_{\Delta x,\vartheta}^b(x))^\top
\triangleq (\tilde{u}(b_{0};s_{0}),\tilde{v}(b_{0};s_{0}))^\top.
\end{align*}

On $x=x_{h}$ for $h\in \mathbb{N}$, the grid points are defined to
be the intersections of line $x=x_{h}$ with the self-similar rays
\begin{align*}
y=\left( b_{0}+n\Delta\sigma\right)x\qquad \mbox{for $n\in\mathbb{Z}$}.
\end{align*}
Here $\Delta\sigma>0$ is chosen so that $\Delta\sigma>\dfrac{4\Delta x}{x_{0}}\max_{i=1,2}\{|\lambda_i(G(s_{0}))|\}$,
and hence the numerical grids satisfies the usual Courant-Friedrichs-Lewy condition.
Then we define the approximate solution $U_{\Delta x,\vartheta}$ to be a piece-wise smooth solution to the self-similar
system (\ref{eqn:eqselfsimilar}), the approximate solution $U_{\Delta x,\vartheta}^b$ on the boundary,
the approximate boundary $\Gamma_{\Delta x,\vartheta}=\big\{(x,y)\,:\,y=b_{\Delta x,\vartheta}(x)\big\}$,
and the numerical grids inductively in $h$, $h=0,1,2,\cdots$.
\par

Suppose that the approximate solution has been defined on $x<x_{h}$. The grid points on $x=x_{h}$ are denoted by $y_n(h)$ for $n\in\mathbb{Z}$. Set
\begin{align*}
r_{h,n}=y_{n}(h)+\vartheta_{h}\big(y_{n+1}(h)-y_{n}(h)\big)\qquad\,\, \mbox{for $n\in\mathbb{Z}$}.
\end{align*}
Then the approximate solution $U_{\Delta x,\vartheta}(x_{h},y)$ for $y\in (y_{n}(h), y_{n+1}(h))$ is defined to
be the solution $U_{self,\Delta x,\vartheta}(\sigma(x,y))$ of (\ref{eqn:selfsimilar}) with the self-similar variable
$\sigma(x,y)=\frac{y}{x_{h}}$ and with the initial data:
\begin{align*}
\sigma=\frac{r_{h,n}}{x_{h}}:\quad U_{self,\Delta x,\vartheta}
=U_{\Delta x,\vartheta}(x_{h},r_{h,n})\triangleq U_{\Delta x,\vartheta}(x_{h}-,r_{h,n}+) \qquad\mbox{for $n\in\mathbb{Z}$}.	
\end{align*} \par

For the discontinuities at the grid points $\big(x_{h},y_{n}(h)\big)$ for $n\in\mathbb{Z}$,
we solve the Riemann problems for (\ref{eqn:simorig}) with the Riemann data:
\begin{equation}
U|_{x=x_{h}}=\left\{\begin{aligned}
&U_{\Delta x,\vartheta}(x_{h},y_{n}(h)-)\qquad \mbox{for $y<y_{n}(h)$}, \\
&U_{\Delta x,\vartheta}(x_{h},y_{n}(h)+)\qquad \mbox{for $y>y_{n}(h)$},
\end{aligned}
\right.
\label{eqn:RieIniDa}
\end{equation}
and the solution consisting of rarefaction waves and shock waves has form $U_{Rie}(\eta)$
with $\eta=\dfrac{y-y_n(h)}{x-x_{h}}$.
Setting $\sigma_{h,n+\frac{1}{2}}\triangleq\dfrac{1}{2x_{h}}\big(y_{n+1}(h+1)+y_{n}(h+1)\big)$ for $n\in\mathbb{Z}$,
then, in the region:
\begin{align*}
\Omega_{h+1,n}=\big\{(x,y)\, :\, x_{h}< x<x_{h+1}, \sigma_{h,n+\frac{1}{2}}>\sigma>\sigma_{h,n-\frac{1}{2}}\big\},
\end{align*}
along the ray
\begin{align*}
\big\{(x,y)\,:\, \frac{y-y_n(h)}{x-x_{h}}=\eta,\,x_{h}<x<x_{h+1}\big\},
\end{align*}
the approximate solution $U_{\Delta x,\vartheta}(x,y)$ is defined to be the solution:
$U_{self,\Delta x,\vartheta}\big(\sigma(x,y)\big)$ of (\ref{eqn:selfsimilar}) with the self-similar
variable $\sigma(x,y)=\frac{y}{x}$ and with the initial data:
\begin{align*}
\sigma=\frac{y_{n}}{x_{h}}:\quad U_{self,\Delta x,\vartheta}=U_{Rie}(\eta).
\end{align*}

The approximate boundary $\Gamma_{\Delta x,\vartheta}=\big\{(x,y)\,:\,y=b_{\Delta x,\vartheta}(x)\big\}$ is traced continuously;
see \cite{Liu1999,Wang2009,Kuang2021}.
For $x\in(0,x_{0})$, let $b_{\Delta x,\vartheta}(x)=b_{0}x$.
Suppose that the approximate solution is constructed for $x<x_{h}$ and
that $y_{n_{b,h}}<b_{\Delta x,\vartheta}(x_{h}-)<y_{n_{b,h}+1}$.
We call interval $y_{n_{b,h}-1}<y<y_{n_{b,h}+1}$ the boundary region at $x=x_{h}$.
In this boundary region, we first solve the self-similar problem (\ref{eqn:selfsimilar}) with the initial data:
\begin{align*}
\sigma=\frac{r_{h,n_{b}-1}}{x_{h}}:\quad U_{self}=U_{\Delta x,\vartheta}(x_{h}-,r_{h,n_{b}-1}+),
\end{align*}
and with the self-similar variable $\sigma(x_{h},y)=\frac{y}{x_{h}}$.
We denote the solution by $U_{self}\big(\sigma(x_{h},y)\big)$.
Given $p^{b}_{\Delta x,h+1}$, by Proposition \ref{Prop:BdryP}, there is $\beta_{1}$ such that
\begin{align*}
\frac{1}{2}\big|\Phi(\beta_{1},0;U_{self}(\sigma\big(x_{h},b_{\Delta x,\vartheta}(x_{h}))))\big|^2
+\frac{\gamma}{\gamma-1}(p^{b}_{\Delta x,h+1})^{\frac{\gamma-1}{\gamma}}
=\frac{1}{2}+\frac{\gamma}{\gamma-1}p_{\infty}^{\frac{\gamma-1}{\gamma}}.
\end{align*}
Then we define
\begin{align*}
U_{\Delta x,\vartheta}^b(x_{h})\triangleq\Phi(\beta_{1},0;U_{self}(\sigma(x_{h},b_{\Delta x,\vartheta}(x_{h})))),
\end{align*}
and
\begin{equation}
b_{\Delta x,\vartheta}(x)
=b_{\Delta x,\vartheta}(x_{h}-)+\frac{v^b_{\Delta x,\vartheta}(x_{h})}{u^b_{\Delta x,\vartheta}(x_{h})}\,(x-x_{h})
\qquad\mbox{for $x\in[x_{h},x_{h+1})$}.
\label{eqn:Defofbdry}
\end{equation}

\smallskip
Next, solve again the self-similar problem (\ref{eqn:selfsimilar}) with initial
data $U_{-}(\sigma(x_{h},b_{\Delta x,\vartheta}(x_{h})))=U_{\Delta x,\vartheta}^b(x_{h})$
and with the self-similar variable $\sigma(x_{h},y)=\frac{y}{x_{h}}$.
Denote the solution by $U_{-}(\sigma(x_{h},y))$. We define the approximate solution in the boundary region as
\begin{align*}
U_{\Delta x,\vartheta}(x_{h},y)=U_{-}(\sigma(x_{h},y))\qquad\mbox{for $x_{h}\leq x<x_{h+1}$}.
\end{align*}
The discontinuities at $(x_{h},y_{n_{b,h}-1})$ are resolved by the same methods as before.	

The leading strong conical shock $S_{\Delta x,\vartheta}=\big\{(x,y)\,:\,y=\chi_{\Delta x,\vartheta}(x)\big\}$
next to the uniform upstream flow is also traced continuously; see \cite{Liu1999,Wang2009,Kuang2021}.
For $x\in(0,x_{0})$, let $\chi_{\Delta x,\vartheta}(x)=s_{0}x$.
Suppose that the approximate solution is constructed for $x<x_{h}$ and
that $y_{n_{\chi,h}-1}<\chi_{\Delta x,\vartheta}(x_{h}-)<y_{n_{\chi,h}}$.
We call interval $y_{n_{\chi,h}-1}<y<y_{n_{\chi,h}+1}$ the front region at $x=x_{h}$.
In this front region, we first solve the self-similar problem (\ref{eqn:selfsimilar}) with the initial data
\begin{align*}
\sigma=\frac{r_{h,n_{\chi}}}{x_{h}}:\quad U_{self}=U_{\Delta x,\vartheta}(x_{h}-,r_{h,n_{\chi}}+),
\end{align*}
and with the self-similar variable $\sigma(x_{h},y)=\frac{y}{x_{h}}$.
Denote the solution by $U_{self}(\sigma(x_{h},y))$.
Then we solve the Riemann problem (\ref{eqn:simorig}) with the initial data
\begin{equation}
U(x_{h},y)=\left\{\begin{aligned}
&U_{\infty},\,\,\, &&y<\chi_{\Delta x,\vartheta}(x_{h}-), \\
&U_{self}(\sigma(x_{h},\chi_{\Delta x,\vartheta}(x_{h}))),
\,\,\,&&\chi_{\Delta x,\vartheta}(x_{h}-)<y<y_{n_{\chi,h}+1}.
\end{aligned}
\right.
\label{eqn:RieIniDaSS}
\end{equation}
The solution $U(x,y)$ contains a weak 2-wave $\beta_{2}$ and a relatively strong 1-shock wave
$s_{\Delta x,\vartheta}(h+1)$ such that
\begin{align*}
U_{self}(\sigma(x_{h},\chi_{\Delta x,\vartheta}(x_{h})))=\Phi(0,\beta_{2};G(s_{\Delta x,\vartheta}(h+1);U_{\infty})).
\end{align*}
Then, let
\begin{equation}
\chi_{\Delta x,\vartheta}(x)=\chi_{\Delta x,\vartheta}(x_{h}-)+s_{\Delta x,\vartheta}(h+1)(x-x_{h})
\qquad\mbox{for $x\in[x_{h},x_{h+1})$}.
\label{eqn:DefofSshock}
\end{equation}

Next, solve again the self-similar problem (\ref{eqn:selfsimilar}) with initial
data $U_{+}(\sigma(x_{h},\chi_{\Delta x,\vartheta}(x_{h})))=G(s_{\Delta x,\vartheta}(h+1);U_{\infty})$
and with the self-similar variable $\sigma(x_{h},y)=\frac{y}{x_{h}}$.
Denote the solution by $U_{+}(\sigma(x_{h},y))$. We define the approximate solution in the front region as
\begin{align*}
U_{\Delta x,\vartheta}(x_h,y)=\left\{\begin{aligned}
&U_{\infty},\,\,\, &&y<\chi_{\Delta x,\vartheta}(x_{h}), \\
&U_{+}(\sigma(x_{h},y)),\,\,\,&&\chi_{\Delta x,\vartheta}(x_{h})<y<y_{n_{\chi,h}+1}.
\end{aligned}
\right.
\end{align*}
The discontinuities at $(x_{h},y_{n_{\chi,h}})$ are resolved by the same methods as before.

\section{Riemann-Type Problems and Interaction Estimates}\label{Section-RIemann-Estimate}	

Let $\Omega_{\Delta x,\vartheta,h}=\{(x,y)\,:\,y<b_{\Delta x,\vartheta}, x_{h-1}\leq x< x_{h}\}$
and $h\in\mathbb{N}_{+}$.
In order to define the approximate solutions in
$\Omega_{\Delta x,\vartheta}\triangleq\bigcup_{k=0}^\infty\Omega_{\Delta x,\vartheta,k}$,
the approximate boundary $\Gamma_{\Delta x,\vartheta}\triangleq\bigcup_{k=0}^\infty \Gamma_{\Delta x,\vartheta,k}$,
and the approximate leading shock
$S_{\Delta x,\vartheta}\triangleq\bigcup_{k=0}^\infty
S_{\Delta x,\vartheta,k}$,
we need a
uniform bound of them to
ensure that all the Riemann problems and the differential equations (\ref{eqn:selfsimilar})
are solvable.
To achieve this,  the following formulas are used:
\begin{enumerate}
\item[\rm (i)]  If $f\in\text{C}^1(\mathbb{R})$, then
\begin{equation}
f(t)-f(0)=t\int_{0}^1f_t(\mu t)\dd\mu\qquad\,\, \mbox{for $t\in\mathbb{R}$}.
\label{eqn:Formula1}
\end{equation}	

\item[\rm (ii)]  If $f\in\text{C}^2(\mathbb{R})$, then
\begin{equation}
f(s,t)-f(s,0)-f(0,t)+f(0,0)=st\int_{0}^1\int_{0}^1f_{st}(\mu s,\lambda t)\dd\mu \dd\lambda\qquad\,\,
\mbox{for $(s,t)\in\mathbb{R}^2$}.
\label{eqn:Formula2}
\end{equation}	
\end{enumerate}

From now on, we use Greek letters $\alpha$, $\beta$, $\gamma$, and $\delta$ to denote the elementary
waves in the approximate solution, and $\alpha_{i}$, $\beta_{i}$, $\gamma_{i}$, and $\delta_{i}$
stand for the corresponding $i$-th components for $i=1,2$.
As in \cite{Chen2006,Smoller1983,Wang2009,Zhang2003}, a curve $I$ is called a mesh curve provided
that $I$ is a space-like curve and consists of the line segments joining the random points one by one in turn.
$I$ divides region $\Omega_{\Delta x, \vartheta}$ into two parts: $I^-$ and $I^+$,
where $I^-$ denotes the part containing line $x=x_{0}$.
For any two mesh curves $I$ and $J$, we use $J>I$ to represent that every mesh point of curve $J$
is either on $I$ or contained in $I^+$.
We say $J$ is an immediate successor to $I$ if $J>I$ and every mesh point of $J$ except one is on $I$ in general but three
when these points are near the approximate boundary
or the approximate shock.

Assume now that $U_{\Delta x,\vartheta}$ has been
defined in $\bigcup_{k=0}^h\Omega_{\Delta x,\vartheta,k}$
and the following conditions are satisfied:
\begin{itemize}[leftmargin=5em]
\item[$H_{1}(h)$:] $\{S_{\Delta x,\vartheta,k}\}_{k=0}^{h}$ forms an approximate strong shock
$S_{\Delta x,\vartheta}|_{0\leq x<x_{h}}$, and $\{\Gamma_{\Delta x,\vartheta,k}\}_{k=0}^{h}$
forms an approximate boundary $\Gamma_{\Delta x,\vartheta}|_{0\leq x<x_{h}}$,
both of which emanate from the origin.

\item[$H_{2}(h)$:]  In each $\Omega_{\Delta x,\vartheta,k}$ for $0\leq k\leq h$,
the strong 1-shock $S_{\Delta x, \vartheta, k}$ divides $\Omega_{\Delta x,\vartheta,k}$
into two parts: $\Omega_{\Delta x,\vartheta,k}^-$ and $\Omega_{\Delta x,\vartheta,k}^+$,
where $\Omega_{\Delta x,\vartheta,k}^+$ is the part between
$S_{\Delta x, \vartheta, k}$ and $\Gamma_{\Delta x, \vartheta, k}$;

\item[$H_{3}(h)$:] $U_{\Delta x,\vartheta}|_{\Omega_{\Delta x,\vartheta,k}^-}=U_{\infty}$,
$\,\,U_{\Delta x,\vartheta}|_{\Omega_{\Delta x,\vartheta,k}^+}
\in O_{\varepsilon_{0}}(G(s_{0}))\cap O_{\varepsilon_{0}}( W(p_{0},p_{\infty}))$, and
$$
U_{\Delta x,\vartheta}(x,b_{\Delta x, \vartheta}(x)-)=U_{\Delta x,\vartheta, k}^{b}
\in  O_{\varepsilon_{0}}(G(s_{0}))\cap O_{\varepsilon_{0}}( W(p_{0},p_{\infty}))
$$
for $x_{k}\leq x<x_{k+1}$,
$0\leq k\leq h$, and $0<\varepsilon_{0}<\min\{\varepsilon_{j},j=1,2,3\}$,
where $\varepsilon_{j}$ are introduced in Propositions \ref{Prop:RiemannP}--\ref{Prop:RPSshock}
for $j=1,2,3$.
\end{itemize}

We prove that $U_{\Delta x,\vartheta}$ can be defined
in $\Omega_{\Delta x,\vartheta,h+1}$ satisfying $H_{1}(h+1)$ -- $H_3(h+1)$.
As in \cite{Glimm1965} (see also \cite{Chen2006,Smoller1983,Kuang2021}),
we consider a pair of the mesh curves $(I,J)$ lying in
$\{x_{h-1}<x<x_{h+1}\}\cap\Omega_{\Delta x, \vartheta}$ with $J$ being an immediate successor of $I$.

\medskip
Now, let $\Lambda$ be the region between $I$ and $J$, and let
\begin{align*}
U_{\Delta x,\vartheta}\in  O_{\varepsilon_{0}}(G(s_{0}))\cap O_{\varepsilon_{0}}(W(p_{0},p_{\infty})).
\end{align*}

\medskip
\Case\label{case:1} {\it $\Lambda$ is between $\Gamma_{\Delta x,\vartheta}$ and $S_{\Delta x, \vartheta}$}.
In this case, we consider the interactions between weak waves.
From the construction of the approximate solutions, the waves entering $\Lambda$ issuing
from $(x_{h-1},y_{n-1}(h-1))$ and from $(x_{h-1},y_{n}(h-1))$ are denoted by
$\alpha=(\alpha_{1},\alpha_{2})$ and $\beta=(\beta_{1},\beta_{2})$, respectively. We denote
\begin{align*}
\sigma_{0}&=\frac{r_{h-1,n}}{x_{h-1}},\
&\bar{\sigma}_{0}&=\frac{r_{h-1,n-1}}{x_{h-1}},\
&\hat{\sigma}_{0}&=\frac{r_{h-1,n-2}}{x_{h-1}},\\
\sigma_{1}&=\frac{y_{n}(h-1)}{x_{h-1}}=\frac{y_{n}(h)}{x_{h}},\
&\sigma_{2}&=\frac{y_{n-1}(h-1)}{x_{h-1}}=\frac{y_{n-1}(h)}{x_{h}},
\end{align*}
and
\begin{align*}
U_{1}=U_{\Delta x,\vartheta}(x_{h-1}-,r_{h-1,n}+),\ U_{2}=U_{\Delta x,\vartheta}(x_{h-1}-,r_{h-1,n-1}+),
\ U_3=U_{\Delta x,\vartheta}(x_{h-1}-,r_{h-1,n-2}+).
\end{align*}

\begin{figure}
\begin{center}
\includegraphics{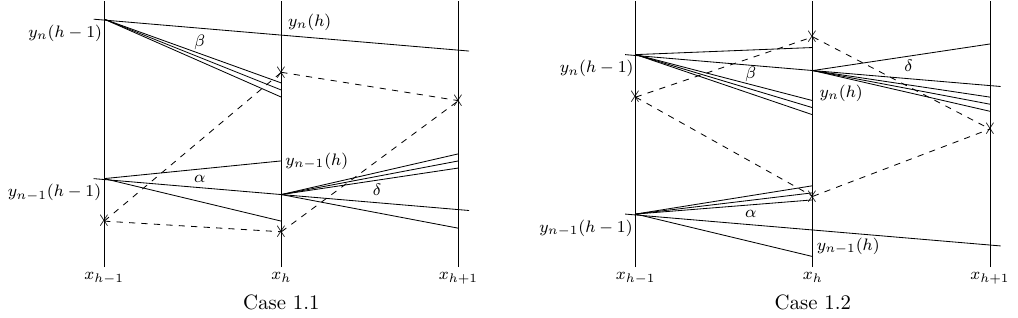}
\end{center}
\caption{Interaction between weak waves}
\label{fig:weak}
\end{figure}

\textbf{Case 1.1.} Let $\delta=(\delta_{1},\delta_{2})$ be the waves
issuing from $(x_{h},y_{n-1}(h))$; see Fig. \ref{fig:weak}.
Then we need to solve the following equations of
$\delta=(\delta_{1},\delta_{2})$:
\begin{equation}
\tilde{U}(\sigma_{1};\sigma_{2},\Phi(\delta_{1},\delta_{2};U_{l}))
=\Phi(\beta_{1},0;\tilde{U}(\sigma_{1};\sigma_{2},\Phi(\alpha_{1},\alpha_{2};U_{l}))),
\label{eqn:weakrela}
\end{equation}
where $U_{l}=\tilde{U}(\sigma_{2};\hat{\sigma}_{0},U_3)$.

\smallskip
\Lemma\label{Lem:Interweak} Equation (\ref{eqn:weakrela}) has a unique solution $\delta=(\delta_{1},\delta_{2})$ such that
\begin{align*}
\delta_{1}=\alpha_{1}+\beta_{1}+O(1)Q(\Lambda),\qquad \delta_{2}=\alpha_{2}+O(1)Q(\Lambda),
\end{align*}
where
\begin{align*}
Q(\Lambda)=Q^0(\Lambda)+Q^1(\Lambda)
\end{align*}
with
\begin{align*}
Q^0(\Lambda)=\sum\big\{|\alpha_{j}||\beta_k|\,:\, \alpha_{j}\text{ and }\beta_k\text{ approach}\big\},\qquad
Q^1(\Lambda)=|\beta_{1}||\Delta\sigma|,
\end{align*}
and $\Delta\sigma=\sigma_{1}-\sigma_{2}$,
where $O(1)$ depends continuously on $M_{\infty}$ but independent of
$(\alpha, \beta, \Delta\sigma)$.

\smallskip
\noindent\textbf{Proof.} Lemma \ref{Lem:DET} yields
\begin{align*}
&\lim_{M_{\infty}\rightarrow\infty}
\det\bigg(\left.\frac{\partial \Phi(\delta_{1},\delta_{2};U_{l})}{\partial(\delta_{1},\delta_{2})}\right|_{\{\delta_{1}=\delta_{2}=0,\, U_{l}\in W(p_{0,U_{\infty}})\}}\bigg)\\[1mm]
&=\dfrac{4\cos^2(\theta_{0}+\theta_{m}^0)\cos^2(\theta_{0}-\theta_{m}^0)\cos^2\theta_{0}\cos^2\theta^0_{m}\sin(2\theta_{m}^0)}{(\gamma+1)^2}.
\end{align*}
Then, by the implicit function theorem,
system (\ref{eqn:weakrela}) has a unique C$^2$--solution:
\begin{align*}
\delta=\delta(\alpha,\beta,\Delta\sigma;U_{l})
\end{align*}
in a neighborhood of $(\alpha,\beta,\Delta\sigma,U_{l})=(0,0,0,G(s_0))$.
Due to (\ref{eqn:Formula2}), we have
\begin{align*}
\delta_{i}(\alpha,\beta,\Delta\sigma;U_{l})
&=\delta_{i}(\alpha,0,\Delta\sigma;U_{l})+\delta_{i}(\alpha,\beta,0;U_{l})
-\delta_{i}(\alpha,0,0;U_{l})+O(1)|\beta||\Delta\sigma|\\
&=\alpha_{i}+\beta_{i}+O(1)Q^0(\Lambda)+O(1)|\beta||\Delta\sigma|
\qquad \mbox{for $i=1,2$},
\end{align*}
where $\beta_{2}=0$. Then the proof is complete.
\hfill{$\square$}

\smallskip
\textbf{Case 1.2.} Let $\delta=(\delta_{1},\delta_{2})$ be the waves
issuing from $(x_{h},y_{n}(h))$; see Fig. \ref{fig:weak}.
Then we need to solve the following equations of
$\delta=(\delta_{1},\delta_{2})$:
\begin{equation}
\Phi(\delta_{1},\delta_{2};\tilde{U}(\sigma_{1};\sigma_{2},U_{l}))
=\Phi(\beta_{1},\beta_{2};\tilde{U}(\sigma_{1},\sigma_{2};\Phi(0,\alpha_{2};U_{l}))),
\label{eqn:weakrela2}
\end{equation}
where $U_{l}$ satisfies $\tilde{U}(\bar{\sigma}_{0};\sigma_{1},\Phi(0,\alpha_{2};U_{l}))=U_{2}$. Similarly, we have the following lemma.

\smallskip
\Lemma\label{Lem:Interweak2} Equation (\ref{eqn:weakrela2}) has a unique solution $\delta=(\delta_{1},\delta_{2})$ such that
\begin{align*}
\delta_{1}=\beta_{1}+O(1)Q(\Lambda),\qquad \delta_{2}=\alpha_{2}+\beta_{2}+O(1)Q(\Lambda),
\end{align*}
where
\begin{align*}
Q(\Lambda)=Q^0(\Lambda)+Q^1(\Lambda)
\end{align*}
with
\begin{align*}
Q^0(\Lambda)=\sum\{|\alpha_{j}||\beta_k|\,:\,\alpha_{j}\text{ and }\beta_k\text{ approach}\},\qquad
Q^1(\Lambda)=|\alpha||\Delta\sigma|,
\end{align*}
and $\Delta\sigma=\sigma_{1}-\sigma_{2}$. Here $O(1)$ depends continuously on $M_{\infty}$ but independent of $(\alpha, \beta, \Delta\sigma)$.

\begin{figure}[t]
\begin{center}
\includegraphics{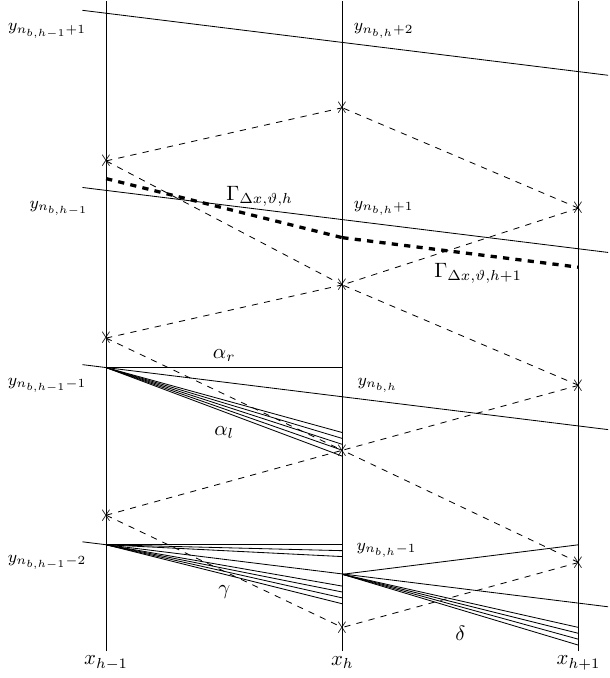}
\end{center}
\caption{Reflection at the boundary}
\label{fig:Reflbdry}
\end{figure}

\medskip
\Case\label{case:2} $\Lambda_b$ covers the part of $\Gamma_{\Delta x,\vartheta}$ but none of $S_{\Delta x, \vartheta}$.
We take three diamonds at the same time, as shown in Fig. \ref{fig:Reflbdry}.
Let $\Delta_{h,n_{b,h}-1}$, $\Delta_{h,n_{b,h}}$, and $\Delta_{h,n_{b,h}+1}$ denote the diamonds centering
in $(x_{h},y_{n_{b,h}-1})$, $(x_{h},y_{n_{b,h}})$, and $(x_{h},y_{n_{b,h}+1})$, respectively,
and denote $\Lambda_b=\Delta_{h,n_{b,h}-1}\cup\Delta_{h,n_{b,h}}\cup\Delta_{h,n_{b,h}+1}$.
Let $\alpha$ and $\gamma$ be the weak waves issuing
from $(x_{h-1},y_{n_{b,h-1}-1})$ and $(x_{h-1},y_{n_{b,h-1}-2})$ respectively, and entering $\Lambda_b$.
We divide $\alpha=(\alpha_{1},\alpha_{2})$ into parts $\alpha_{l}=(\alpha_{l,1},0)$ and $\alpha_{r}=(\alpha_{r,1},\alpha_{r,2})$,
where $\alpha_{l}$ and $\alpha_{r}$ entering $\Delta_{h,n_{b,h}-1}$ and $\Delta_{h,n_{b,h}}$, respectively.
Moreover, let $\gamma=(\gamma_{1},\gamma_{2})$ issuing
from $(x_{h-1},y_{n_{b,h}-2})$, and let $\delta$ be the outgoing wave issuing from $(x_{h},y_{n_{b,h}-1})$.
For simplicity of notation, we denote
\begin{align*}
&\sigma_{\alpha}=\sigma(x_{h-1},y_{n_{b,h-1}-1}),&\,\,\,\,
&\sigma_{b}(h-1)=\sigma(x_{h-1},b_{\Delta x,\vartheta}(x_{h-1})),&\,\,\, &\sigma_{b}(h)=\sigma(x_{h},b_{\Delta x,\vartheta}(x_{h})),&\\
&\Delta\sigma_{\alpha}=\sigma_{b}(h-1)-\sigma_{\alpha},&\,\,\,\,
&\Delta\bar{\sigma}_{\alpha}=\sigma_{b}(h)-\sigma_{\alpha},&\,\,\,\,
&\Delta\sigma_{b_{h}}=\sigma_{b}(h)-\sigma_{b}(h-1),&\\
&\sigma_{\gamma}=\sigma(x_{h-1},y_{n_{b,h-1}-2}),&\,\,\,\,
&\Delta \sigma_{\gamma}= \sigma_{\alpha}-\sigma_{\gamma},&\,\,\,\, &\sigma_{0}=\sigma(x_{h-1},r_{h-1,n_b-2}),&
\end{align*}
and $U_{1}=U_{\Delta x,\vartheta}(x_{h-1}-,r_{h-1,n_b-2}+)$. Let $U_{l}=\Phi(\alpha_{l,1},0;\tilde{U}(\sigma_{\alpha};\sigma_{0},U_{1}))$.

\par

\smallskip
To gain the estimates of $\delta$, we need to deal with the equation:
\begin{align}	
&\dfrac{1}{2}|\Phi(\beta_{1},0;\tilde{U}(\sigma_{b}(h);\sigma_{\alpha},U_{l}))|^2
+\dfrac{\gamma }{\gamma-1}(p^{b}_{\Delta x,h+1})^{\frac{\gamma-1}{\gamma}}\nonumber\\
&=
\dfrac{1}{2}|\tilde{U}(\sigma_{b}(h-1);\sigma_{\alpha},\Phi(\alpha_{r,1},\alpha_{r,2};U_{l}))|^2
+\dfrac{\gamma }{\gamma-1}(p^{b}_{\Delta x,h})^{\frac{\gamma-1}{\gamma}},
\label{eqn:reflbry1}
\end{align}
and then we obtain the following lemma.

\smallskip
\Lemma\label{Lem:Refbound} Equation (\ref{eqn:reflbry1}) has a unique solution
$\beta_{1}=\beta_{1}(\alpha_{r,1},\alpha_{r,2},\Delta \sigma_\alpha,\Delta\bar{\sigma}_\alpha,\omega_{h+1};U_{l})\in\text{C}^2$ in a neighborhood of
$(\alpha_{r,1},\alpha_{r,2},\Delta \sigma_\alpha,\Delta\bar{\sigma}_\alpha,\omega_{h+1},U_{l})=(0,0,0,0,0,G(s_{0}))$
with
$\omega_{h+1}=p^{b}_{\Delta x,h+1}-p^{b}_{\Delta x,h}$ such that
\begin{equation}
\begin{aligned}
\delta_{1}&=\alpha_{r,1}+\alpha_{l,1}+\gamma_{1}
+K_{r,1}\alpha_{r,2}
+K_{\sigma,1}\Delta\sigma_{b_{h}}+K_{b,1}\omega_{h+1}+O(1)Q(\Lambda_b)\\
\delta_{2}&=\gamma_{2}+K_{r,2}\alpha_{r,2}
+K_{\sigma,2}\Delta\sigma_{b_{h}}+K_{b,2}\omega_{h+1}+O(1)Q(\Lambda_b),
\end{aligned}
\label{eqn:Refbdryesti}
\end{equation}
with
\begin{equation}	
Q(\Lambda_b)=Q^0((\alpha_{1},0),\gamma)+|\alpha_{1}||\Delta\sigma_\gamma|+|\alpha_{r,1}||\Delta\sigma_\alpha|,
\label{eqn:Refbdryesti1}
\end{equation}
where $O(1)$ depends continuously on $M_{\infty}$.
Moreover, when $\alpha_{r,1}=\alpha_{r,2}=\Delta\sigma_\alpha=\Delta\bar{\sigma}_\alpha=\omega_{h+1}=0$, $p^{b}_{\Delta x,h+1}=p_{0}$, and $U_{l}=G(s_{0})$,
\begin{equation}
\lim_{M_{\infty}\rightarrow\infty}K_{r,1}=-\dfrac{\cos^2(\theta_{0}+\theta_{m}^0)}{\cos^2(\theta_{0}-\theta_{m}^0)},\
\lim_{M_{\infty}\rightarrow\infty}|K_{b,i}|<\infty,\
\lim_{M_{\infty}\rightarrow\infty}K_{r,2}=0,\
\lim_{M_{\infty}\rightarrow\infty}K_{\sigma,i}=0,
\label{eqn:Refbdryesti2}
\end{equation}
for $i=1,2$.

\smallskip
\noindent
\textbf{Proof.} A direct computation leads to
\begin{align*}	
\left.\dfrac{1}{2}
\frac{\partial (|\Phi(\beta_{1},0;\tilde{U}(\sigma_{b}(h);\sigma_{\alpha},U_{l}))|^2)}{\partial\beta_{1}}\right|
_{\{\delta_{1}=\Delta\bar{\sigma}_\alpha=0, U_{l}=G(s_{0})\}}=r_{1}(G(s_{0}))\cdot G(s_{0}).
\end{align*}
Lemma \ref{Lem:asymsol}, together with the implicit function theorem,
implies that there is a unique $\text{C}^2-$solution
\begin{equation}
\beta_{1}=\beta_{1}(\alpha_{r,1},\alpha_{r,2},\Delta\sigma_\alpha,\Delta\bar{\sigma}_\alpha,\omega_{h+1};U_{l})
\label{eqn:rebdrsol}
\end{equation}
in a neighborhood of
$(\alpha_{r,1},\alpha_{r,2},\Delta\sigma_\alpha,\Delta\bar{\sigma}_\alpha,\omega_{h+1},U_{l})=(0,0,0,0,0,G(s_{0}))$.\par

Using (\ref{eqn:Formula1})--(\ref{eqn:Formula2}), we have
\begin{align*}
\beta_{1}&=\beta_{1}(\alpha_{r,1},\alpha_{r,2},\Delta\sigma_\alpha,\Delta\sigma_\alpha,\omega_{h+1};U_{l})
+\bar{K}_{\sigma,1}(\Delta\bar{\sigma}_\alpha-\Delta\sigma_\alpha)\\
&=\beta_{1}(\alpha_{r,1},0,\Delta\sigma_\alpha,\Delta\sigma_\alpha,0;U_{l})
+\bar{K}_{\sigma,1}\Delta\sigma_{b_{h}}+\bar{K}_{r,1}\alpha_{r,2}+\bar{K}_{b,1}\omega_{h+1}\\
&=\alpha_{r,1}
+\bar{K}_{\sigma,1}\Delta\sigma_{b_{h}}+\bar{K}_{r,1}\alpha_{r,2}+\bar{K}_{b,1}\omega_{h+1}+O(1)|\alpha_{r,1}||\Delta\sigma_\alpha|.
\end{align*}
Taking derivative with respect to $\Delta\sigma_{b_{h}}$ in (\ref{eqn:reflbry1})
at $(\alpha_{r,1},\alpha_{r,2},\Delta\sigma_\alpha,\Delta\bar{\sigma}_\alpha,\omega_{h+1},U_{l})=(0,0,0,0,0,G(s_{0}))$,
we obtain
\begin{align*}	
G(s_{0})\cdot r_{1}(G(s_{0}))\frac{\partial\beta_{1}}{\partial\Delta\sigma_{b_{h}}}
+G(s_{0})\cdot\frac{\partial\tilde{U}(\sigma_\alpha+\Delta\bar{\sigma}_\alpha+\Delta\sigma_{b_{h}};\sigma_\alpha,G(s_{0}))}{\partial\Delta\sigma_{b_{h}}}=0,
\end{align*}
which yields
\begin{align*}	
\lim_{M_{\infty}\rightarrow\infty}
\left.\frac{\partial\beta_{1}}{\partial\Delta\sigma_{b_{h}}}\right|
_{\{\alpha_{r,1}=\alpha_{r,2}=\Delta\sigma_\alpha=\Delta\bar{\sigma}_\alpha=\omega_{h+1}=0,\,p^{b}_{\Delta x,h+1}=p_{0},\,U_{l}=G(s_{0})\}}=0.
\end{align*}
Similarly, we have
\begin{align*}	
&\lim_{M_{\infty}\rightarrow\infty}\left.\frac{\partial\beta_{1}}{\partial\omega_{h+1}}\right|
_{\{\alpha_{r,1}=\alpha_{r,2}=\Delta\sigma_\alpha=\Delta\bar{\sigma}_\alpha=\omega_{h+1}=0,\,p^{b}_{\Delta x,h+1}=p_{0},\,U_{l}=G(s_{0})\}}
=\lim_{M_{\infty}\rightarrow\infty}\dfrac{-p_{0}^{-\frac{1}{\gamma}}}{r_{1}(G(s_{0}))\cdot G(s_{0})}>-\infty,\\
&\lim_{M_{\infty}\rightarrow\infty}\left.\frac{\partial\beta_{1}}{\partial\alpha_{r,2}}\right|
_{\{\alpha_{r,1}=\alpha_{r,2}=\Delta\sigma_\alpha=\Delta\bar{\sigma}_\alpha=\omega_{h+1}=0,\,p^{b}_{\Delta x,h+1}=p_{0},\,U_{l}=G(s_{0})\}}
=-\dfrac{\cos^2(\theta_{0}+\theta_{m}^0)}{\cos^2(\theta_{0}-\theta_{m}^0)}.
\end{align*}
By the construction of the approximate solution, we have
\begin{align*}
\tilde{U}(\sigma_{b}(h);\sigma_\gamma,\Phi(\delta_{1},\delta_{2};U_m))
=\Phi(\beta_{1},0;\tilde{U}(\sigma_{b}(h);\sigma_{\alpha},\Phi(\alpha_{l,1},0;\tilde{U}(\sigma_\alpha;\sigma_\gamma,\Phi(\gamma_{1},\gamma_{2};U_m)))))
\end{align*}
with $U_m=U_{\Delta x,\vartheta}(x_{h},y_{n_{b,h-1}-2}-)$. Then, a similar argument as in Case \ref{case:1}
gives (\ref{eqn:Refbdryesti})--(\ref{eqn:Refbdryesti2}).
This completes the proof.\hfill{$\square$}

\medskip
\Lemma\label{Lem:DistofCent} In Case \ref{case:2},
for $b'_{h}:=b_{\Delta x,\vartheta}'(x_{h}-)$ for $h\in\mathbb{N}_{+}$,
\begin{align*}
b'_{h+1}-b'_{h}=K_{c,2}\alpha_{r,2}+K_{c,\sigma}\Delta\sigma_{b_{h}}+O(1)\omega_{h+1}+O(1)|\alpha_{r,1}||\Delta\sigma_{\alpha}|
\end{align*}
with $O(1)$ depending continuously on $p_{0}$ such that
\begin{align*}
&\lim_{M_{\infty}\rightarrow\infty}K_{c,\sigma}|_{\{\alpha_{r,2}=\alpha_{r,1}=\Delta\sigma_\alpha=\Delta\bar{\sigma}_\alpha=\omega_{h+1}=0,
p^{b}_{\Delta x,h+1}=p_{0},U_{l}=G(s_{0})\}}=-1,\\
&\lim_{M_{\infty}\rightarrow\infty}K_{c,2}|_{\{\alpha_{r,2}=\alpha_{r,1}=\Delta\sigma_\alpha=\Delta\bar{\sigma}_\alpha=\omega_{h+1}=0,
p^{b}_{\Delta x,h+1}=p_{0},U_{l}=G(s_{0})\}}
=-\frac{4}{\gamma+1}\,\frac{\cos^2\theta^0_m\cos^2(\theta^0_m+\theta_{0})}{\cos^2\theta_{0}}.
\end{align*}

\smallskip
\noindent
\textbf{Proof.} From (\ref{eqn:Defofbdry}), we have
\begin{align*}
b'_{h+1}-b'_{h}
&=\frac{\Phi^{(2)}(\beta_{1}(\alpha_{r,1},\alpha_{r,2},\Delta\sigma_\alpha,\Delta\bar{\sigma}_\alpha,\omega_{h+1};U_{l}),0;
\tilde{U}(\sigma_{b_{h}};\sigma_{\alpha},U_{l}))}
{\Phi^{(1)}(\beta_{1}(\alpha_{r,1},\alpha_{r,2},\Delta\sigma_\alpha,\Delta\bar{\sigma}_\alpha,\omega_{h+1};U_{l}),0;\tilde{U}(\sigma_{b_{h}};\sigma_{\alpha},U_{l}))}\\
&\quad-\frac{\tilde{U}^{(2)}(\sigma_{b_{h-1}};\sigma_{\alpha},\Phi(\alpha_{r,1},\alpha_{r,2};U_{l}))}
{\tilde{U}^{(1)}(\sigma_{b_{h-1}};\sigma_{\alpha},\Phi(\alpha_{r,1},\alpha_{r,2};U_{l}))}.
\end{align*}
By (\ref{eqn:Formula1})--(\ref{eqn:Formula2}), we obtain
\begin{align*}
b'_{h+1}-b'_{h}
&=\frac{\Phi^{(2)}(\beta_{1}(\alpha_{r,1},0,\Delta\sigma_\alpha,\Delta\bar{\sigma}_\alpha,0;U_{l}),0;\tilde{U}(\sigma_{b_{h}};\sigma_{\alpha},U_{l}))}
{\Phi^{(1)}(\beta_{1}(\alpha_{r,1},0,\Delta\sigma_\alpha,\Delta\bar{\sigma}_\alpha,0;U_{l}),0;\tilde{U}(\sigma_{b_{h}};\sigma_{\alpha},U_{l}))}
-\frac{\tilde{U}^{(2)}(\sigma_{b_{h-1}};\sigma_{\alpha},\Phi(\alpha_{r,1},0;U_{l}))}
{\tilde{U}^{(1)}(\sigma_{b_{h-1}};\sigma_{\alpha},\Phi(\alpha_{r,1},0;U_{l}))}\\
&\quad +K_{c,2}\alpha_{r,2}+O(1)\omega_{h+1}\\
&=\frac{\Phi^{(2)}(\beta_{1}(\alpha_{r,1},0,\Delta\sigma_\alpha,\Delta\sigma_\alpha,0;U_{l}),0;\tilde{U}(\sigma_{b_{h-1}};\sigma_{\alpha},U_{l}))}
{\Phi^{(1)}(\beta_{1}(\alpha_{r,1},0,\Delta\sigma_\alpha,\Delta\sigma_\alpha,0;U_{l}),0;\tilde{U}(\sigma_{b_{h-1}};\sigma_{\alpha},U_{l}))}
-\frac{\tilde{U}^{(2)}(\sigma_{b_{h-1}};\sigma_{\alpha},\Phi(\alpha_{r,1},0;U_{l}))}
{\tilde{U}^{(1)}(\sigma_{b_{h-1}};\sigma_{\alpha},\Phi(\alpha_{r,1},0;U_{l}))}\\
&\quad +K_{c,\sigma}\Delta\sigma_{b_{h}}+K_{c,2}\alpha_{r,2}+O(1)\omega_{h+1}\\
=&K_{c,\sigma}\Delta\sigma_{b_{h}}+K_{c,2}\alpha_{r,2}+O(1)\omega_{h+1}+O(1)|\alpha_{r,1}||\Delta\sigma_{\alpha}|.
\end{align*}
By similar calculation in Lemma \ref{Lem:Refbound}, we have
\begin{align*}
&\lim_{M_{\infty}\rightarrow\infty}K_{c,\sigma}|_{\{\alpha_{r,2}=\alpha_{r,1}=\Delta\sigma_\alpha=\Delta\bar{\sigma}_\alpha=\omega_{h+1}=0,
\,p^{b}_{\Delta x,h+1}=p_{0},\,U_{l}=G(s_{0})\}}=-1,\\
&\lim_{M_{\infty}\rightarrow\infty}K_{c,2}|_{\{\alpha_{r,2}=\alpha_{r,1}=\Delta\sigma_\alpha=\Delta\bar{\sigma}_\alpha=\omega_{h+1}=0,
\,p^{b}_{\Delta x,h+1}=p_{0},\,U_{l}=G(s_{0})\}}
=-\frac{4}{\gamma+1}\,\frac{\cos^2\theta^0_m\cos^2(\theta^0_m+\theta_{0})}{\cos^2\theta_{0}}.
\end{align*}
Then the proof is complete.\hfill{$\square$}

\medskip
\Lemma\label{Lem:RoughEstbdry} For $\Delta x$ sufficiently small,
\begin{align*}
|b'_{h}-\sigma_{b}(h-1)| \geq6|\Delta\sigma_{b_h}|,
\end{align*}
where $\sigma_{b}(h)=\frac{b_{\Delta x,\vartheta}(x_{h})}{x_{h}}$.

\smallskip
\noindent\textbf{Proof.} Using the notation as in Case \ref{case:2}, we have
\begin{align*}
b'_{h}=\dfrac{b_{\Delta x,\vartheta}(x_{h})-b_{\Delta x,\vartheta}(x_{h-1})}{\Delta x}.
\end{align*}
Then a direct computation leads to
\begin{align*}
|b'_h-\sigma_{b}(h-1)|
&=\left|\dfrac{b_{\Delta x,\vartheta}(x_{h})-b_{\Delta x,\vartheta}(x_{h-1})}{\Delta x}-\sigma_{b}(h-1)\right|\\
&=\left|\dfrac{\sigma_{b}(h)x_{h}-\sigma_{b}(h-1)x_{h-1}}{\Delta x}-\sigma_{b}(h-1)\right|\\
&=\left|\dfrac{x_{h}}{\Delta x}\right||\sigma_{b}(h)-\sigma_{b}(h-1)|\\[1mm]
&\geq 6\left|\sigma_{b}(h)-\sigma_{b}(h-1)\right|
\end{align*}
for $\Delta x$ small enough.\hfill{$\square$}

\medskip	
Denote $\theta_{b}(h)=|\sigma_{b}(h-1)-b'_{h}|$ that measures the angle
between boundary $\Gamma_{\Delta x,\vartheta,h}$ and the ray issuing from
the origin and passing through $(x_{h-1},b_{\Delta x,\vartheta}(x_{h-1}))$. Then we have the following estimate for $\theta_{b}(h)$.

\smallskip
\Lemma\label{Lem:Estchabdry} For $M_{\infty}$ sufficiently large and $\Delta x$ sufficiently small,
\begin{align*}
\theta_{b}(h)-\theta_{b}(h+1)\geq |\Delta \sigma|-|K_{c,2}||\alpha_{r,2}|-C|\omega_{h+1}|-C|\alpha_{r,1}||\Delta\sigma_{\alpha}|,
\end{align*}
where $h\in\mathbb{N}_{+}$, and constant $C>0$ is independent of $M_{\infty}$ and $\Delta x$.	

\smallskip
\noindent\textbf{Proof.} We consider the following two different cases.

\smallskip
1. $\sigma_{b}(h-1)<b'_{h}$ so that $\sigma_{b}(h)>\sigma_{b}(h-1)$.
\begin{itemize}
\item If $b'_{h+1}>\sigma_{b}(h)$, then it follows from Lemma \ref{Lem:DistofCent}
that
\begin{align*}
\theta_{b}(h)-\theta_{b}(h+1)
&= b'_{h}-\sigma_{b}(h-1)-\big(b'_{h+1}-\sigma_{b}(h)\big)\\
&= (1-K_{c,\sigma})\Delta\sigma_{b_{h}}-K_{c,2}\alpha_{r,2}+O(1)\omega_{h+1}+O(1)|\alpha_{r,1}||\Delta\sigma_{\alpha}|\\
&\geq |\Delta\sigma_{b_{h}}|-|K_{c,2}||\alpha_{r,2}|-C|\omega_{h+1}|-C|\alpha_{r,1}||\Delta\sigma_{\alpha}|.
\end{align*}	

\item If $b'_{h+1}<\sigma_{b}(h)$, then, from
Lemma \ref{Lem:DistofCent}--\ref{Lem:RoughEstbdry}, we have
\begin{align*}
\theta_{b}(h)-\theta_{b}(h+1)
&= b'_{h}-\sigma_{b}(h-1)-\big(\sigma_{b}(h)-b'_{h+1}\big)\\
&= 2\big(b'_{h}-\sigma_{b}(h-1)\big)+b'_{h+1}-\sigma_{b}(h)-\big(b'_{h}-\sigma_{b}(h-1)\big)\\
&\geq (11+K_{c,\sigma})|\Delta\sigma_{b_{h}}|+K_{c,2}\alpha_{r,2}+O(1)\omega_{h+1}+O(1)|\alpha_{r,1}||\Delta\sigma_{\alpha}|\\
&\geq |\Delta\sigma_{b_{h}}|-|K_{c,2}||\alpha_{r,2}|-C|\omega_{h+1}|-C|\alpha_{r,1}||\Delta\sigma_{\alpha}|.
\end{align*}	
\end{itemize}

2. $\sigma_{b}(h-1)>b'_{h}$ so that $\sigma_{b}(h)<\sigma_{b}(h-1)$.

\begin{itemize}
\item If $b'_{h+1}>\sigma_{b}(h)$, then it follows
from Lemma \ref{Lem:DistofCent}--\ref{Lem:RoughEstbdry} that
\begin{align*}
\theta_{b}(h)-\theta_{b}(h+1)
&=\sigma_{b}(h-1)-b'_{h}-\big(b'_{h+1}-\sigma_{b}(h)\big)\\
&= 2(\sigma_{b}(h-1)-b'_{h})+\sigma_{b}(h)-b'_{h+1}-(\sigma_{b}(h-1)-b'_{h})\\
&\geq (11+K_{c,\sigma})|\Delta\sigma_{b_{h}}|-K_{c,2}\alpha_{r,2}-O(1)\omega_{h+1}-O(1)|\alpha_{r,1}||\Delta\sigma_{\alpha}|\\
&\geq |\Delta\sigma_{b_{h}}|-|K_{c,2}||\alpha_{r,2}|-C|\omega_{h+1}|-C|\alpha_{r,1}||\Delta\sigma_{\alpha}|.
\end{align*}	

\item If $b'_{h+1}<\sigma_{b}(h)$, then, from Lemma \ref{Lem:DistofCent}, we have
\begin{align*}
\theta_{b}(h)-\theta_{b}(h+1)
&=\sigma_{b}(h-1)-b'_{h}-\big(\sigma_{b}(h)-b'_{h+1}\big)\\
&=(-1+K_{c,\sigma})\Delta\sigma_{b_{h}}+K_{c,2}\alpha_{r,2}+O(1)\omega_{h+1}+O(1)|\alpha_{r,1}||\Delta\sigma_{\alpha}|\\
&\geq |\Delta\sigma_{b_{h}}|-|K_{c,2}||\alpha_{r,2}|-C|\omega_{h+1}|
-C|\alpha_{r,1}||\Delta\sigma_{\alpha}|.
\end{align*}	
\end{itemize}
Note that we have used the fact that
$$
\lim_{M_{\infty}\rightarrow\infty}K_{c,\sigma}|_{\{\alpha_{r,2}=\beta_{1}=\Delta\sigma_\alpha=\Delta\bar{\sigma}_\alpha=\omega_{h+1}=0,
p^{b}_{\Delta x,h+1}=p_{0},U_{l}=G(s_{0})\}}=-1
$$
in above estimates.
This completes the proof.\hfill{$\square$}

\begin{figure}[t]
\begin{center}
\includegraphics{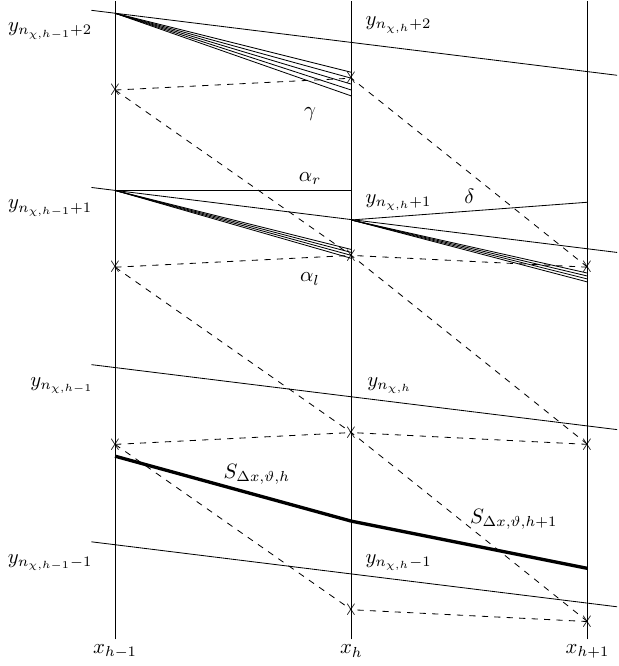}
\end{center}
\caption{Near the strong shock wave}
\label{fig:Strongshock}
\end{figure}

\medskip
\Case\label{case:3} $\Lambda_{s}$ covers the part of $S_{\Delta x,\vartheta}$ but none of $\Gamma_{\Delta x, \vartheta}$.
We take three diamonds at the same time, as shown in Fig. \ref{fig:Strongshock}.
Let $\Delta_{h,n_{\chi,h}-1}$, $\Delta_{h,n_{\chi,h}}$, and
$\Delta_{h,n_{\chi,h}+1}$ be the diamonds centering in $(x_{h},y_{n_{\chi,h}-1})$, $(x_{h},y_{n_{\chi,h}})$, and $(x_{h},y_{n_{\chi,h}+1})$, respectively.
Denote $\Lambda_{s}=\Delta_{h,n_{\chi,h}-1}\cup\Delta_{h,n_{\chi,h}}\cup\Delta_{h,n_{\chi,h}+1}$.
Let $\alpha$ and $\gamma$ be the weak waves issuing from
$(x_{h-1},y_{n_{\chi,h-1}+1})$ and $(x_{h-1},y_{n_{\chi,h-1}+2})$ respectively and entering $\Lambda_{s}$. We divide $\alpha$ into
parts $\alpha_{l}=(\alpha_{l,1},0)$ and $\alpha_{r}=(\alpha_{r,1},\alpha_{r,2})$ where $\alpha_{l}$ and $\alpha_{r}$
enter $\Delta_{h,n_{\chi,h}}$ and $\Delta_{h,n_{\chi,h}+1}$, respectively.
Moreover, let $\gamma=(\gamma_{1},0)$, and
let $\delta$ be the outgoing wave issuing from $(x_{h},y_{n_{\chi,h}+1})$.
\par
Then, for simplicity of notation, we denote
\begin{align*}
&\sigma_{\alpha}=\sigma(x_{h-1},y_{n_{\chi,h-1}+1}),&\,\,\,
&\sigma_{\chi}(h-1)=\sigma(x_{h-1},\chi_{\Delta x,\vartheta}(x_{h-1})),&\,\,\, &\sigma_{\chi}(h)=\sigma(x_{h},\chi_{\Delta x,\vartheta}(x_{h})),&\\
&\Delta\sigma_{\alpha}=\sigma_{\alpha}-\sigma_{\chi}(h-1),&\,\,\,
&\Delta\bar{\sigma}_{\alpha}=\sigma_{\alpha}-\sigma_{\chi}(h),&\,\,\,
&\Delta\sigma_{\chi_{h}}=\sigma_{\chi}(h)-\sigma_{\chi}(h-1),&\\
&\sigma_{\gamma}=\sigma(x_{h-1},y_{n_{\chi,h-1}+2}),&\,\,\,
&\Delta \sigma_{\gamma}= \sigma_{\gamma}-\sigma_{\alpha}.&
\end{align*}

To gain the estimates of $(s_{h+1},\delta)$, we need to deal with the equation:
\begin{equation}
\tilde{U}(\sigma_\alpha;\sigma_{\chi}(h),\Phi(0,\beta_{2};G(s_{h+1};U_{\infty})))
=\Phi(\alpha_{l,1},0;\tilde{U}(\sigma_\alpha;\sigma_{\chi}(h-1),G(s_{h};U_{\infty}))),
\label{eqn:Interstrong}
\end{equation}
to obtain the following lemma.

\medskip
\Lemma\label{Lem:Intstrong} Equation (\ref{eqn:Interstrong}) has a unique solution $(s_{h+1},\beta_{2})$ in a neighborhood of
\begin{align*}
(\alpha_{l,1},\alpha_{r},\gamma,\Delta\sigma_{\alpha},\Delta \sigma_{\chi_{h}},s_{h})=(0,0,0,0,0,s_{0}),
\end{align*}
such that
\begin{equation}
\begin{aligned}
&\delta_{1}=\alpha_{r,1}+\gamma_{1}+\mu_{w,1}\Delta\sigma_{\chi_{h}}+K_{w,1}\alpha_{l,1}+O(1)Q(\Lambda_{s}),\\
&\delta_{2}=\alpha_{r,2}+\mu_{w,2}\Delta\sigma_{\chi_{h}}+K_{w,2}\alpha_{l,1}+O(1)Q(\Lambda_{s}),\\
&s_{h+1}= s_{h}+K_{s}\alpha_{l,1}+\mu_{s}\Delta \sigma_{\chi_{h}},
\end{aligned}
\label{eqn:Refshockesti}
\end{equation}
with
\begin{equation}	
Q(\Lambda_{s})=|\gamma_{1}||\Delta\sigma_{\gamma}|+Q^0(\alpha_{r},\gamma),
\label{eqn:Refshockesti1}
\end{equation}
where $O(1)$ depends continuously on $M_{\infty}$.
In addition, for $\alpha_{l}=0$, $\Delta\sigma_{\alpha}=\Delta\sigma_{\chi_{h}}=0$, and $s_{h}=s_{0}$,
denoting the derivative of $G$ by $G_{s}$, then
\begin{equation}
\begin{aligned}
&\lim_{M_{\infty}\rightarrow \infty}K_{w,1}=0,\quad K_{w,2}=\dfrac{\det(r_1(G(s_0)),\,G_{s}(s_0))}{\det(r_2(G(s_0)),\,G_{s}(s_0))},\quad K_s=\dfrac{\det(r_2(G(s_0)),r_1(G(s_0))}{\det(r_2(G(s_0)),G_{s}(s_0))},\\
&\lim_{M_{\infty}\rightarrow\infty}\mu_{w,1}=0,\quad
\mu_{w,2}=\dfrac{\det(\partial\tilde{U}/\partial(\Delta\sigma_{\chi_h}),\,G_{s}(s_0))}{\det(r_2(G(s_0)),\,G_{s}(s_0))},\quad
\lim_{M_{\infty}\rightarrow\infty}\mu_s\in(-1,0).
\end{aligned}	
\label{eqn:Refshockesti2}
\end{equation}

\noindent
\textbf{Proof.} From Lemma \ref{Lem:DET} and the implicit function theorem, (\ref{eqn:Interstrong})
has a unique $\text{C}^2$--solution $(s_{h+1},\beta_{2})$ such that
\begin{align*}
&s_{h+1}=s_{h+1}(\alpha_{l,1},\Delta\sigma_{\alpha},\Delta\bar{\sigma}_{\alpha},\Delta\sigma_{\chi_{h}},s_{h}),\\
&\beta_{2}=\beta_{2}(\alpha_{l,1},\Delta\sigma_{\alpha},\Delta\bar{\sigma}_{\alpha},\Delta\sigma_{\chi_{h}},s_{h}).
\end{align*}
A direct computation leads to
\begin{align*}
\beta_{2}
=\beta_{2}(\alpha_{l,1},\Delta\sigma_{\alpha},\Delta\bar{\sigma}_{\alpha},\Delta\sigma_{\chi_{h}},s_{h})
&=\mu_{w,2}\Delta\sigma_{\chi_{h}}+K_{w,2}\alpha_{l,1}
+\beta_{2}(0,\Delta\sigma_{\alpha},\Delta\sigma_{\alpha},0,s_{h})\\
&=\mu_{w,2}\Delta\sigma_{\chi_{h}}+K_{w,2}\alpha_{l,1}.
\end{align*}
Similarly, we have
\begin{align*}	
s_{h+1}=s_{h+1}(\alpha_{l,1},\Delta\sigma_{\alpha},\Delta\bar{\sigma}_{\alpha},\Delta\sigma_{\chi_{h}},s_{h})
=\mu_s\Delta\sigma_{\chi_{h}}+K_{s}\alpha_{l,1}+s_{h}.
\end{align*}

Next, we compute the coefficients: $K_s$, $K_{w,2}$, $\mu_{w,2}$, and $\mu_s$. Differentiating equation (\ref{eqn:Interstrong})
with respect to $\alpha_{l,1}$ and $\Delta\sigma_{\chi_{h}}$, and then letting $\alpha_{l,1}=\Delta\sigma_{\alpha}=\Delta\sigma_{\chi_{h}}=0$ and $s_{h}=s_{0}$, we can obtain
\begin{align*}
&r_{2}(G(s_{0}))K_{w,2}+G_s(s_{0})K_s=r_{1}(G(s_{0})),\\
&r_{2}(G(s_{0}))\mu_{w,2}+G_s(s_{0})\mu_s
=\dfrac{\partial\tilde{U}}{\partial (\Delta\sigma_{\chi_{h}})}(\sigma_{\chi}(h);\sigma_{\chi}(h),G(s_{0})).
\end{align*}
Then Cramer's rule gives the result.
Moreover, since $\theta_{0}<0<\theta^0_{ma}$ and $\theta_{0}\pm\theta^0_{ma}\in(-\frac{\pi}{2},\frac{\pi}{2})$,
\begin{align*}
\lim_{M_{\infty}\rightarrow\infty}\mu_s=\frac{\cos\theta_{0}\sin\theta_{m}^0}{\sin(\theta_{0}-\theta_{m}^0)}\in(-1,0).
\end{align*}

By the construction of the approximate solution, we have
\begin{align*}
&\tilde{U}(\sigma_{\gamma};\sigma_{\alpha},\Phi(\delta_{1},\delta_{2};\tilde{U}(\sigma_\alpha;\sigma_{\chi}(h),U_m)))\\
&=\Phi(\gamma_{1},0; \tilde{U}(\sigma_{\gamma};\sigma_{\alpha},
\Phi(\alpha_{r,1},\alpha_{r,2};
\tilde{U}(\sigma_\alpha;\sigma_{\chi}(h),
\Phi(0,\beta_{2}(\alpha_{l,1},\Delta\sigma_{\alpha},\Delta\bar{\sigma}_{\alpha},\Delta\sigma_{\chi_{h}},s_{h});U_m))))),
\end{align*}
with
\begin{align*}
U_m=G(s_{h+1};U_{\infty}).
\end{align*}
Then, by similar arguments as in Case \ref{case:1}, we obtain
\begin{align*}
\delta_{1}&=\alpha_{r,1}+\gamma_{1}+\mu_{w,1}\Delta\sigma_{\chi_{h}}+K_{w,1}\alpha_{l,1}+O(1)|\gamma_{1}||\Delta\sigma_{\gamma}|+O(1)Q^0(\alpha_{r},\gamma),\\
\delta_{2}&=\alpha_{r,2}+\mu_{w,2}\Delta\sigma_{\chi_{h}}+K_{w,2}\alpha_{l,1}+O(1)|\gamma_{1}||\Delta\sigma_{\gamma}|+O(1)Q^0(\alpha_{r},\gamma).
\end{align*}
This completes the proof.\hfill{$\square$}

\medskip
\Lemma\label{Lem:RoughEst} For $\Delta x$ sufficiently small,
\begin{align*}
|s_h-\sigma_{\chi}(h-1)|\geq6|\Delta\sigma_{\chi_{h}}|.
\end{align*}

\noindent\textbf{Proof.}
Using the notation as in Case \ref{case:3}, we have
\begin{align*}
\sigma_{\chi}(h)=\frac{\chi_{\Delta x,\vartheta}(x_{h})}{x_{h}},\qquad s_{h}=\dfrac{\chi_{\Delta x,\vartheta}(x_{h})-\chi_{\Delta x,\vartheta}(x_{h-1})}{\Delta x}.
\end{align*}
Then a direct computation leads to
\begin{align*}
|s_h-\sigma_{\chi}(h-1)|
&=\left|\dfrac{\chi_{\Delta x,\vartheta}(x_{h})-\chi_{\Delta x,\vartheta}(x_{h-1})}{\Delta x}-\sigma_{\chi}(h-1)\right|\\
&=\left|\dfrac{\sigma_{\chi}(h)x_{h}-\sigma_{\chi}(h-1)x_{h-1}}{\Delta x}-\sigma_{\chi}(h-1)\right|\\
&=\left|\dfrac{x_{h}}{\Delta x}\right||\sigma_{\chi}(h)-\sigma_{\chi}(h-1)|\\
&\geq 6\left|\sigma_{\chi}(h)-\sigma_{\chi}(h-1)\right|,
\end{align*}
for $\Delta x$ small enough.\hfill{$\square$}	

\medskip	
Denote $\theta_{\chi}(h)=|\sigma_{\chi}(h-1)-s_{h}|$ that measures
the angle between the leading shock $S_{\Delta x,\vartheta,h}$ and the ray issuing from the origin and passing through
$(x_{h-1},\chi_{\Delta x,\vartheta}(x_{h-1}))$.
Then we have the following estimate for $\theta_{\chi}(h)$.

\smallskip
\Lemma\label{Lem:Estchashock} For $M_{\infty}$ sufficiently large and $\Delta x$ sufficiently small,
\begin{align*}
\theta_{\chi}(h)-\theta_{\chi}(h+1)\geq |\Delta \sigma_{\chi_{h}}|-|K_s||\alpha_{l,1}|,
\end{align*}
with $h\geq0$.

\smallskip
\noindent\textbf{Proof.} We consider the following two different cases:

\smallskip
1. $\sigma_{\chi}(h-1)<s_{h}$ so that $\sigma_{\chi}(h)>\sigma_{\chi}(h-1)$.
\begin{itemize}
\item If $s_{h+1}>\sigma_{\chi}(h)$, then it follows from
Lemma \ref{Lem:Intstrong} that
\begin{align*}
\theta_{\chi}(h)-\theta_{\chi}(h+1)
&= s_{h}-\sigma_{\chi}(h-1)-\big(s_{h+1}-\sigma_{\chi}(h)\big)
=(1-\mu_s)\Delta\sigma_{\chi_{h}}-K_s\alpha_{l,1}\\
&\geq |\Delta\sigma_{\chi_{h}}|-|K_s||\alpha_{l,1}|.
\end{align*}	

\item If $s_{h+1}<\sigma_{\chi}(h)$, then,
from Lemmas \ref{Lem:Intstrong}--\ref{Lem:RoughEst}, we have
\begin{align*}
\theta_{\chi}(h)-\theta_{\chi}(h+1)
&=s_{h}-\sigma_{\chi}(h-1)-\big(\sigma_{\chi}(h)-s_{h+1}\big)\\
&= 2(s_{h}-\sigma_{\chi}(h-1))+s_{h+1}-\sigma_{\chi}(h)-\big(s_{h}-\sigma_{\chi}(h-1)\big)\\
&\geq (11+\mu_s)|\Delta\sigma_{\chi_{h}}|+K_s\alpha_{l,1}\\
&\geq |\Delta\sigma_{\chi_{h}}|-|K_s||\alpha_{l,1}|.
\end{align*}	
\end{itemize}

2. $\sigma_{\chi}(h-1)>s_{h}$ so that $\sigma_{\chi}(h)<\sigma_{\chi}(h-1)$.
\begin{itemize}
\item If $s_{h+1}>\sigma_{\chi}(h)$, then it follows from
Lemmas \ref{Lem:Intstrong}--\ref{Lem:RoughEst} that
\begin{align*}
\theta_{\chi}(h)-\theta_{\chi}(h+1)
&=\sigma_{\chi}(h-1)-s_{h}-\big(s_{h+1}-\sigma_{\chi}(h)\big)\\
&=2(\sigma_{\chi}(h-1)-s_{h})+\sigma_{\chi}(h)-s_{h+1}-(\sigma_{\chi}(h-1)-s_{h})\\
&\geq (11+\mu_s)|\Delta\sigma_{\chi_{h}}|-K_s\alpha_{l,1}\\
&\geq |\Delta\sigma_{\chi_{h}}|-|K_s||\alpha_{l,1}|.
\end{align*}	

\item If $s_{h+1}<\sigma_{\chi}(h)$,
then, from Lemma \ref{Lem:Intstrong}, we have
\begin{align*}
\theta_{\chi}(h)-\theta_{\chi}(h+1)
&=\sigma_{\chi}(h-1)-s_{h}-\big(\sigma_{\chi}(h)-s_{h+1}\big)\\
&=(-1+\mu_s)\Delta\sigma_{\chi_{h}}+K_{s}\alpha_{l,1}\\
&\geq |\Delta\sigma_{\chi_{h}}|-|K_s||\alpha_{l,1}|.
\end{align*}	
\end{itemize}
Note that we have used the fact that $\mu_s\in(-1,0)$
as $M_{\infty}\rightarrow \infty$ in above estimates.
This completes the proof.\hfill{$\square$}

\section{Glimm-Type Functional and Compactness of the Approximate Solutions}\label{Section-Decrease}

For each $I\subset\cup_{k=1}^{h+1}\Omega_{\Delta x,\vartheta,k}$,
there exists $k_{I}$ with $1\leq k_{I}\leq h+1$ such that
$I\cap \Gamma_{\Delta x,\vartheta, k_{I}}\neq \emptyset$.
Next, as in \cite{Zhang1999,Liu1999}, we assign each mesh curve
$I\subset\bigcup_{k=1}^{h+1}\Omega_{\Delta x,\vartheta,k}$ with a Glimm-type functional $F_s(I)$; see also \cite{Wang2009,Kuang2021}:

\medskip
\Definition (Weighted total variation). Define
\begin{align*}
&L_{0}^{(i)}(I)=\sum\{|\alpha_{i}|:\alpha_{i} \text{ is the weak $i$-wave crossing $I$}\} \qquad \mbox{for $i=1,\,2$},\\
&L_{1}(I)=\sum\{|\omega_{k}|: k>k_{I}\},\\
&L_s(I)=\theta_\chi(I) \quad \text{ for $\theta_{\chi}(I)=\theta_\chi(h)$ in Lemma \ref{Lem:Estchashock} when $S_{\Delta x,\vartheta}$ crossing $I$},\\
&L_b(I)=\theta_b(I) \quad\,\text{ for $\theta_{b}(I)=\theta_b(h)$ in Lemma \ref{Lem:Estchabdry} when $\Gamma_{\Delta x,\vartheta}$ crossing $I$}.
\end{align*}
Then the weighted total variation is defined as
\begin{align*}
L(J)=L_{0}^{(1)}(I)+K_{2}L_{0}^{(2)}(I)+K_{1}L_{1}(I)+K_3L_s(I)+K_4L_b(I),
\end{align*}
where $K_{l}$ are positive constants for $l=1,2,3,4$.

\medskip	
Let
\begin{equation}
\sigma^*=b_{0}+C_{1}\sum_{h=1}^{\infty}|\omega_{h}|,\qquad \sigma_*=s_{0}-\varpi,
\label{eqn:sigmabound}
\end{equation}
where $s_{0}$ is the velocity of the leading shock of the background solution, $\varpi$ and $C_{1}$
are constants to be determined; see also \cite{Chen2004,Wang2009,Kuang2021}.
Note that $\varpi$ and $\sum_{h\geq1}|\omega_{h}|$ are chosen so small
that the largeness of $M_{\infty}$ implies the smallness of $b_{0}-s_{0}$,
which leads to the smallness of $\sigma^*-\sigma_*$.
We now define the total interaction potential.\par

\medskip
\Definition (Total interaction potential). Define
\begin{align*}
&Q_{0}(I)
=\sum \{|\alpha||\beta|\,:\,\text{$\alpha$ and $\beta$ are weak waves crossing $I$ and approach}\},\\	
&Q_{1}(I)
=\sum\{|\alpha||\sigma_\alpha-\sigma_*|\,:\,\text{$\alpha$ is a weak $1-$wave crossing $I$}\},\\
&Q_{2}(I)=\sum\{|\alpha||\sigma^*-\sigma_\alpha|\,:\,\text{$\alpha$ is a weak $2-$wave crossing $I$}\},
\end{align*}
where $\sigma_\alpha$ is the $\sigma$-coordinate of the grid point where $\alpha$ issues. Then the total interaction potential is defined as
\begin{align*}
Q(I)=Q_{0}(I)+2Q_{1}(I)+2Q_{2}(I).
\end{align*}

\smallskip
Now, we are able to define the Glimm-type functional.

\medskip
\Definition (Glimm-type functional).  Let
\begin{align*}
F(I)=L(I)+KQ(I),
\end{align*}
where $K$ is a big real number to be chosen later.

\medskip
Let
\begin{equation}
E_{\Delta x,\vartheta}(\Lambda)=\left\{\begin{aligned}
&Q(\Lambda)  &&\text{(defined in Case \ref{case:1}),}&\\
&\xi\big(|\alpha_{r,2}|+|\omega_{h+1}|+|\Delta\sigma_{b_{h}}|+Q(\Lambda_b)\big)
&&\text{(defined in Case \ref{case:2}),}&\\
&\xi\big(|\alpha_{l,1}|+|\Delta\sigma_{\chi_{h}}|+Q(\Lambda_s)\big)
&&\text{(defined in Case \ref{case:3}),}&
\end{aligned}
\right.
\label{eqn:DiffGLimmFunc}
\end{equation}
with $\xi>0$ sufficiently small and to be chosen later.

\medskip	
In order to make the Glimm-type functional monotonically decreasing, we have to choose the weights carefully in the functional, based on the underlying
features of the wave interactions governed by the system. Indeed, we have the following lemma ({\it cf.}
\cite{Wang2009}).

\medskip
\Lemma\label{Lem:coeffi}
Let $K_{r,1}$, $K_{w,2}$, $K_s$ and $\mu_{w,2}$ be given
by Lemmas \ref{Lem:Refbound} and \ref{Lem:Intstrong}. Then
\begin{align*}
\lim_{M_{\infty}\rightarrow\infty}
\big(|K_{r,1}||K_{w,2}|+|K_{r,1}||K_s||\mu_{w,2}|\big)<1.
\end{align*}

\smallskip
\noindent\textbf{Proof.}
Lemmas \ref{Lem:asymshockpol}--\ref{Lem:asymselfsi} give
\begin{align*}
&\lim_{M_{\infty}\rightarrow\infty}|K_{r,1}||K_{w,2}|=\left|\dfrac{\sin(\theta_{0}+\theta_{m}^0)}{\sin(\theta_{0}-\theta_{m}^0)} \right|,\\[2mm]
&\lim_{M_{\infty}\rightarrow\infty}|K_{r,1}||K_s||\mu_{w,2}|\\
&\,\,\,=\dfrac{\cos^2(\theta_{0}+\theta_{m}^0)}{\cos^2(\theta_{0}-\theta_{m}^0)}\lim_{M_{\infty}\rightarrow\infty}
 \left|\dfrac{\det \big(r_{1}(U),r_{2}(U)\big)}{\det \big(r_{2}(G(s_{0}),G'(s_{0}))\big)} \right|
 \left| \dfrac{\det((\partial\tilde{U})/(\partial \Delta\sigma_{\chi_{h}}),G_s(s_{0};U_{\infty}))}{\det(r_{2}(G(s_{0};U_{\infty})),G_s(s_{0};U_{\infty}))}\right|
\\
&\,\,\,=\dfrac{\sin2\theta_{m}^0\cos\theta_{0}|\sin\theta_{0}|}{\sin^2(\theta_{0}-\theta_{m}^0)}.
\end{align*}
Note that $\theta_{0}\in(-\frac{\pi}{2},0)$
and $\theta_{0}\pm\theta_{m}^0\in(-\frac{\pi}{2},\frac{\pi}{2})$, and that $\theta_{m}^0\in(0,\frac{\pi}{2})$.
Then, when $\theta_{0}+\theta_{m}^0<0$,
\begin{align*}
\lim_{M_{\infty}\rightarrow\infty}
\big(|K_{r,1}||K_{w,2}|+|K_{r,1}||K_s||\mu_{w,2}|\big)
<\dfrac{2\sin\theta_{m}^0\cos\theta_{0}-\sin(\theta_{0}+\theta_{m}^0)}{\sin(\theta_{m}^0-\theta_{0})}=1;
\end{align*}
when $\theta_{0}+\theta_{m}^0>0$,
\begin{align*}
\lim_{M_{\infty}\rightarrow\infty}\big(|K_{r,1}||K_{w,2}|+|K_{r,1}||K_s||\mu_{w,2}|\big)
<\dfrac{2\cos\theta_{m}^0|\sin\theta_{0}|+\sin(\theta_{0}+\theta_{m}^0)}{\sin(\theta_{m}^0-\theta_{0})}=1.
\end{align*}
This implies the expected result.\hfill{$\square$}

\medskip
At this stage, we are able to choose the coefficients in the Glimm-type functional ({\it cf.} \cite{Wang2009}).

\smallskip
\Lemma\label{Lem:GlimWeigh} There exist positive constants $K_{2}$ and $K_3$ such that
\begin{align*}
\lim_{M_{\infty}\rightarrow\infty}\big(K_{2}|K_{w,2}|+K_3|K_s|\big)<1,\quad \lim_{M_{\infty}\rightarrow\infty}\big(K_{2}|\mu_{w,2}|-K_3\big)<0,
\quad
\lim_{M_{\infty}\rightarrow\infty}\big(K_{2}-|K_{r,1}|\big)>0.
\end{align*}

\smallskip
\noindent
\textbf{Proof.}
Let $K_{r,1}^*=\lim_{M_{\infty}\rightarrow\infty}|K_{r,1}|$, $K_{w,2}^*=\lim_{M_{\infty}\rightarrow\infty}|K_{w,2}|$, $K_s^*=\lim_{M_{\infty}\rightarrow\infty}|K_s|$,
and $\mu_{w,2}^*=\lim_{M_{\infty}\rightarrow\infty}|\mu_{w,2}|$.
Then, by Lemma \ref{Lem:coeffi},
\begin{align*}
K_{r,1}^*\big(K_{w,2}^*+K_s^*\mu_{w,2}^*\big)<1.
\end{align*}
Hence, we choose $K_{2}$ such that
\begin{align*}
K_{2}>K_{r,1}^*,\qquad K_{2}\big(K_{w,2}^*+K_s^*\mu_{w,2}^*\big)<1,
\end{align*}
which implies
\begin{align*}
K_{2}K_s^*\mu_{w,2}^*<1-K_{2}K_{w,2}^*.
\end{align*}
Then we take $K_3$ such that
\begin{align*}
K_3>K_{2}\mu_{w,2}^*,\qquad K_3K_s^*<1-K_{2}K_{w,2}^*,
\end{align*}
and the proof is complete.\hfill{$\square$}

\medskip
With the coefficients chosen properly, we can derive a decay property for the Glimm-type functional.

\medskip
\Proposition\label{Prop:Estimate}
Let $M_{\infty}$ be sufficiently large, and let $\sigma^{*}-\sigma_{*}$ and $\sum_{h\geq1}|\omega_{h}|$ be sufficiently small.
Let $I$ and $J$ be a pair of space-like mesh curves with $J$ being an immediate successor of $I$.
The region bounded by the difference between $I$ and $J$ is denoted as $\Lambda$.
Then there exist positive constants $\epsilon_\infty$, $K$,
and $K_{l}$ for $l=1,2,3,4$, such that, if $F(I)<\epsilon_\infty$, then
\begin{align*}
F(I)\le F(J)-\frac{1}{4}E_{\Delta x, \vartheta}(\Lambda),
\end{align*}
where $E_{\Delta x, \vartheta}(\Lambda)$ is given by (\ref{eqn:DiffGLimmFunc}).

\medskip
\noindent
\textbf{Proof.}
When $M_{\infty}$ is large enough, according to Lemma \ref{Lem:GlimWeigh}, there are constants $K_{2}$ and $K_3$ so that
\begin{align*}
K_{2}|K_{w,2}|+K_3|K_s|<1-\xi_{0},\quad K_{2}|\mu_{w,2}|-K_3<-\xi_{0},\quad K_{2}-|K_{r}|-K_4|K_{c,2}|>\xi_{0},
\end{align*}
for some $\xi_{0}>0$.\par

Now, as in \cite{Liu1999}, we prove the result inductively; see also \cite{Wang2009,Kuang2021}.
We consider three special cases as
in \S \ref{Section-RIemann-Estimate},
depending on the location of  $\Lambda$.
From now on, we use $C$ to denote a universal constant depending only on the system, which may be different at each occurrence.\par

\smallskip
\noindent\textbf{Case 1.} $\Lambda$ lies between $\Gamma_{\Delta x,\vartheta}$
and $S_{\Delta x, \vartheta}$.
We consider the case as in Lemma \ref{Lem:Interweak}.
Notice that
\begin{align*}
&(L_{0}^{(1)}+K_{2}L_{0}^{(2)})(J)-(L_{0}^{(1)}+K_{2}L_{0}^{(2)})(I)\leq CQ(\Lambda),\\
&L_b(J)-L_b(I)=0,\\
&(K_{1}L_{1}+K_3L_s)(J)-(K_{1}L_{1}+K_3L_s)(I)=0.
\end{align*}
Then we obtain
\begin{align*}
L(J)-L(I)\leq CQ(\Lambda).
\end{align*}
For the terms contained in $Q$, we have
\begin{align*}
Q_{0}(J)-Q_{0}(I)\leq CL(I)Q(\Lambda)-Q^0(\Lambda).
\end{align*}

\smallskip
For \textbf{Case 1.1:}
\begin{align*}
(Q_{1}+Q_{2})(J)-(Q_{1}+Q_{2})(I)
&= |\delta_{1}|(\sigma_{2}-\sigma_*)-|\alpha_{1}|(\sigma_{2}-\sigma_*)-|\beta_{1}|(\sigma_{1}-\sigma_*)\\
&\quad +|\delta_{2}|(\sigma^*-\sigma_{2})-|\alpha_{2}|(\sigma^*-\sigma_{2})\\
&\leq C(\sigma^*-\sigma_*)Q(\Lambda)-|\Delta\sigma||\beta_{1}|.
\end{align*}

\medskip
For \textbf{Case 1.2:}
\begin{align*}
(Q_{1}+Q_{2})(J)-(Q_{1}+Q_{2})(I)&\leq|\delta_{1}|(\sigma_{1}-\sigma_*)-|\beta_{1}|(\sigma_{1}-\sigma_*)\\
&\quad+|\delta_{2}|(\sigma^*-\sigma_{1})-|\alpha_{2}|(\sigma^*-\sigma_{2})-|\beta_{2}|(\sigma^*-\sigma_{1})\\
&\leq C(\sigma^*-\sigma_*)Q(\Lambda)-|\Delta\sigma||\alpha_{2}|,
\end{align*}
which gives
\begin{align*}
Q(J)-Q(I)\leq&-\big(1-C(L(I)+\sigma^*-\sigma_*)\big)Q(\Lambda).
\end{align*}
When $L(I)$ and $\sigma^*-\sigma_*$ are small enough, and $K$ is sufficiently large, it follows that
\begin{align*}
F(J)-F(I)\leq-\left\{K\big(1-C(L(I)+\sigma^*-\sigma_*)\big)-C\right\}Q(\Lambda)\leq-\frac{1}{4}Q(\Lambda).
\end{align*}

\smallskip
\noindent\textbf{Case 2.} $\Lambda_b=\Delta_{h,n_{b,h}-1}\cup\Delta_{h,n_{b,h}}\cup\Delta_{h,n_{b,h}+1}$
covers a part of $\Gamma_{\Delta x,\vartheta}$ but none of $S_{\Delta x, \vartheta}$. Direct computation shows that
\begin{align*}
&L_{0}^{(1)}(J)-L_{0}^{(1)}(I)\leq |K_{r,1}||\alpha_{r,2}|
+|K_{\sigma,1}||\Delta\sigma_{b_{h}}|+|K_{b,1}||\omega_{h+1}|+CQ(\Lambda_b),\\
&L_{0}^{(2)}(J)-L_{0}^{(2)}(I)\leq -|\alpha_{r,2}|+|K_{r,2}||\alpha_{r,2}|
+|K_{\sigma,2}||\Delta\sigma_{b_{h}}|+|K_{b,2}||\omega_{h+1}|+CQ(\Lambda_b),\\
&L_{1}(J)-L_{1}(I)=-|\omega_{h+1}|,\\
&L_s(J)-L_s(I)=0,\\
&L_b(J)-L_b(I)=-|\Delta\sigma_{b_{h}}|+|K_{c,2}||\alpha_{r,2}|+C|\omega_{h+1}|+C|\alpha_{r,1}||\Delta\sigma_{\alpha}|.
\end{align*}
Combining the above estimates together, we obtain
\begin{align*}
L(J)-L(I)&\leq-\big(K_{2}-|K_{r,1}|-K_4|K_{c,2}|-K_{2}|K_{r,2}|\big)|\alpha_{r,2}|\\
&\quad -\big(K_{1}-|K_{b,1}|-K_{2}|K_{b,2}|-CK_4\big)|\omega_{h+1}|\\
&\quad -\big(K_4-|K_{\sigma,1}|-K_{2}|K_{\sigma,2}|\big)|\Delta \sigma_{b_{h}}|+CQ(\Lambda_b)+C|\alpha_{r,1}||\Delta\sigma_{\alpha}|.
\end{align*}
For the terms contained in $Q$,
noting that $|\Delta\sigma_{\alpha}|\leq|\Delta\sigma_{\gamma}|$,
we have
\begin{align*}
Q_{0}(J)-Q_{0}(I)&\leq -Q^0((\alpha_{1},0),\gamma)+CL(I)\big(|\alpha_{r,2}|
+|\Delta\sigma_{b_{h}}|+|\omega_{h+1}|+CQ(\Lambda_b)\big),\\
Q_{1}(J)-Q_{1}(I)&=|\delta_{1}|(\sigma_{\gamma}-\sigma_*)-(|\alpha_{l,1}|+|\alpha_{r,1}|)(\sigma_{\alpha}-\sigma_*)-|\gamma_{1}|(\sigma_{\gamma}-\sigma_*)\\
&\leq-(|\alpha_{l,1}|+|\alpha_{r,1}|)|\Delta\sigma_{\gamma}|+C(\sigma^*-\sigma_*)\big(|\alpha_{r,2}|
+|\Delta\sigma_{b_{h}}|+|\omega_{h+1}|+CQ(\Lambda_b)\big),\\
Q_{2}(J)-Q_{2}(I)&=|\delta_{2}|(\sigma^*-\sigma_{\gamma})-|\alpha_{r,2}|(\sigma^*-\sigma_{\alpha})-|\gamma_{2}|(\sigma^*-\sigma_{\gamma})\\
&\leq C(\sigma^*-\sigma_*)\big(|\alpha_{r,2}|
+|\Delta\sigma_{b_{h}}|+|\omega_{h+1}|+CQ(\Lambda_b)\big).
\end{align*}
Then we conclude
\begin{align*}
Q(J)-Q(I)&\leq
-Q^0((\alpha_{1},0),\gamma)+CL(I)\big(|\alpha_{r,2}|
+|\Delta\sigma_{b_{h}}|+|\omega_{h+1}|+CQ(\Lambda_b)\big)-|\alpha_{1}||\Delta\sigma_{\gamma}|\\
&\quad -|\alpha_{r,1}||\Delta\sigma_{\alpha}|+2C(\sigma^*-\sigma_*)\big(|\alpha_{r,2}|
+|\Delta\sigma_{b_{h}}|+|\omega_{h+1}|+CQ(\Lambda_b)\big)\\
&\leq -\big(1-C(L(I)+\sigma^*-\sigma_*)\big)Q(\Lambda_b)+C\big(L(I)+\sigma^*-\sigma_*\big)\big(|\alpha_{r,2}|
+|\Delta\sigma_{b_{h}}|+|\omega_{h+1}|\big).
\end{align*}
Finally, combining all the estimates above together, we obtain
\begin{align*}
F(J)-F(I)&\leq-\left\{K\big(1-C(L(I)+\sigma^*-\sigma_*)\big)-
C\right\}Q(\Lambda_b)\\
&\quad-\left\{K_{2}-|K_{r,1}|-K_4|K_{c,2}|-K_{2}|K_{r,2}|-KC\big(L(I)+\sigma^*-\sigma_*\big)\right\}|\alpha_{r,2}|\\
&\quad-\left\{K_{1}-|K_{b,1}|-K_{2}|K_{b,2}|-CK_4-KC\big(L(I)+\sigma^*-\sigma_*\big)\right\}|\omega_{h+1}|\\
&\quad-\left\{K_4-|K_{\sigma,1}|-K_{2}|K_{\sigma,2}|-KC\big(L(I)+\sigma^*-\sigma_*\big)\right\}|\Delta \sigma_{b_{h}}|.
\end{align*}
Taking suitably large $K_{1}$, then,
when $K$ is sufficiently large, and $L(I)$ and $\sigma^*-\sigma_*$ are sufficiently small, we conclude
\begin{align*}
F(J)-F(I)\leq -\frac{\xi}{4}\big(|\alpha_{r,2}|+|\omega_{h+1}|+|\Delta\sigma_{b_{h}}|+Q(\Lambda_b)\big),
\end{align*}
for some $\xi>0$ small enough.

\medskip
\noindent\textbf{Case 3.} $\Lambda_{s}=\Delta_{h,n_{\chi,h}-1}\cup\Delta_{h,n_{\chi,h}}\cup\Delta_{h,n_{\chi,h}+1}$ covers a part of $S_{\Delta x,\vartheta}$
but none of $\Gamma_{\Delta x, \vartheta}$.
A direct computation shows that
\begin{align*}
&L_{0}^{(1)}(J)-L_{0}^{(1)}(I)\leq -|\alpha_{l,1}|+|K_{w,1}||\alpha_{l,1}|+|\mu_{w,1}||\Delta\sigma_{\chi_{h}}|+CQ(\Lambda_{s}),\\
&L_{0}^{(2)}(J)-L_{0}^{(2)}(I)
\leq |K_{w,2}||\alpha_{l,1}|+|\mu_{w,2}||\Delta \sigma_{\chi_{h}}|+CQ(\Lambda_{s}),\\
&L_{1}(J)-L_{1}(I)=0,\\
&L_s(J)-L_s(I)\leq -|\Delta\sigma_{\chi_{h}}|+|K_s||\alpha_{l,1}|,\\
&L_b(J)-L_b(I)=0.
\end{align*}
Combine the above estimates together, we obtain
\begin{align*}
&L(J)-L(I)\\
&\leq -(1-|K_{w,1}|-K_{2}|K_{w,2}|-K_3|K_s|)|\alpha_{l,1}|-\big(K_3-|\mu_{w,1}|-K_{2}|\mu_{w,2}|\big)|\Delta\sigma_{\chi_{h}}|+CQ(\Lambda_{s})\\
&\leq -\big(1-K_{2}|K_{w,2}|-K_3|K_s|-|K_{w,1}|\big)|\alpha_{l,1}|
-\big(K_3-K_{2}|\mu_{w,2}|-|\mu_{w,1}|\big)|\Delta\sigma_{\chi_{h}}|+CQ(\Lambda_{s}).
\end{align*}

For the terms contained in $Q$, we have
\begin{align*}
Q_{0}(J)-Q_{0}(I)&\leq -Q^0(\alpha_{r},\gamma)
+CL(I)\big(|\alpha_{l,1}|+|\Delta\sigma_{\chi_{h}}|+Q(\Lambda_{s})\big),\\
Q_{1}(J)-Q_{1}(I)&= |\delta_{1}|(\sigma_{\alpha}-\sigma_*)-(|\alpha_{l,1}|+|\alpha_{r,1}|)(\sigma_{\alpha}-\sigma_*)-|\gamma_{1}|(\sigma_{\gamma}-\sigma_*)\\
&\leq -|\gamma_{1}||\Delta\sigma_{\gamma}|+C(\sigma^*-\sigma_*)\big(|\alpha_{l,1}|+|\Delta\sigma_{\chi_{h}}|+Q(\Lambda_{s})\big),\\
Q_{2}(J)-Q_{2}(I)&= |\delta_{2}|(\sigma^*-\sigma_{\gamma})-|\alpha_{r,2}|(\sigma^*-\sigma_{\alpha})\\
&\leq C(\sigma^*-\sigma_*)\big(|\alpha_{l,1}|+|\Delta\sigma_{\chi_{h}}|+Q(\Lambda_{s})\big).
\end{align*}
Then we deduce that
\begin{align*}
Q(J)-Q(I)\leq
-\big(1-C(L(I)+\sigma^*-\sigma_*)\big)Q(\Lambda_{s})+C\big(L(I)+\sigma^*-\sigma_*\big)\big(|\alpha_{l,1}|+|\Delta\sigma_{\chi_{h}}|\big).
\end{align*}
Finally, combining all the estimates above together, we obtain
\begin{align*}
F(J)-F(I)\leq&-\left\{K\big(1-C(L(I)+\sigma^*-\sigma_*)\big)
-C\right\}Q(\Lambda_{s})\\
&-\Big\{1-K_{2}|K_{w,2}|-K_3|K_s|-|K_{w,1}|-CK\big(L(I)+\sigma^*-\sigma_*\big)\Big\}|\alpha_{l,1}|\\
&-\Big\{K_3-K_{2}|\mu_{w,2}|-|\mu_{w,1}|-CK\big(L(I)+\sigma^*-\sigma_*\big)\Big\}|\Delta\sigma_{\chi_{h}}|.
\end{align*}
When $K$ is sufficiently large, and $L(I)$ and $\sigma^*-\sigma_*$ are sufficiently small, we conclude that
\begin{align*}
F(J)-F(I)\leq&-\frac{\xi}{4}(Q(\Lambda_{s})+|\alpha_{l,1}|+|\Delta\sigma_{\chi_{h}}|)
\end{align*}
for some $\xi>0$ small enough.
Combining the above three cases, we conclude our result.\hfill{$\square$}\\

Now, let $I_{h}$ be the mesh curve in the stripe:
$\{(x,y)\,:\, x_{h-1}\leq x\leq x_{h}\}$ for $h\in\mathbb{N}_{+}$;
that is, $I_{h}$ connects all the mesh points in
the strip.
Let $I$ and $J$ be any pair of mesh curves with $I_{h}<I<J<I_{h+1}$,
and let $J$ be an immediate successor of $I$.
That is, the mesh points on $J$ differ from those on $I$ by only one point
generally (except three points near the approximate boundary or near the approximate shock),
and the region bounded by the difference between $I$ and $J$ is denoted by $\Lambda$.
Proposition \ref{Prop:Estimate} suggests that the total variation of the approximate solution is uniformly bounded.

\medskip
Moreover, we have the following estimates for the approximate boundary and the approximate leading shock:

\smallskip
\Proposition\label{Prop:EstimateBdr}
There exists a constant $\bar{C}>0$, independent of $\Delta x$, $\vartheta$,
and $U_{\Delta x, \vartheta}$, such that
\begin{align*}
&\text{T.V.}\{s_{\Delta x,\vartheta}:[0,\infty)\}=\sum_{h=0}^{\infty}|s_{h+1}-s_{h}|\leq \bar{C}\sum_{h\geq1}|\omega_{h}|,\\
&\text{T.V.}\{b'_{\Delta x,\vartheta}:[0,\infty)\}=\sum_{h=0}^{\infty}|b'_{h+1}-b'_{h}|\leq \bar{C}\sum_{h\geq1}|\omega_{h}|.
\end{align*}

\noindent\textbf{Proof.}
Notice that
\begin{align*}
\text{T.V.}\{s_{\Delta x,\vartheta}:[0,\infty)\}
&=\sum_{h=0}^{\infty}|s_{h+1}-s_{h}|\leq O(1)\sum_{\Lambda_s}E_{\Delta x,\vartheta}(\Lambda_s)\\
&\leq O(1)\sum_{\Lambda}F(I)-F(J)\leq O(1)F(I_{1}).
\end{align*}
Similarly, we have
\begin{align*}
\text{T.V.}\{b'_{\Delta x,\vartheta}:[0,\infty)\}\leq O(1)F(I_{1}).
\end{align*}
Therefore, $\bar{C}$ in the statement can be determined.
\hfill{$\square$}

\smallskip
We choose $C_{1}=2\bar{C}$ and $\varpi=2\bar{C}\sum_{h\geq1}|\omega_{h}|$ in (\ref{eqn:sigmabound}).
The largeness of $M_{\infty}$ and the smallness of $\sum_{h\geq1}|\omega_{h}|$ imply the smallness of $\sigma^{*}-\sigma_{*}$.
Then, following \cite{Chen2006,Zhang2003}, we conclude

\smallskip
\Theorem \label{Thm:compactness}
Under assumptions (\ref{Asum:pres})--(\ref{Asum:super}),
if $M_{\infty}$ is sufficiently large and $\sum_{h\geq1}|\omega_{h}|$ is
sufficiently small,
then, for any $\vartheta\in\Pi_{h=0}^{\infty}[0,1)$ and $\Delta x>0$,
the modified Glimm scheme introduced above defines a sequence of global approximate solutions $U_{\Delta x,\vartheta}(x,y)$ such that
\begin{align*}
&\sup_{x>0}\text{T.V.}\{U_{\Delta x,\vartheta}(x,y):(-\infty,b_{\Delta x,\vartheta}(x))\}<\infty,\\
&\int_{-\infty}^0|U_{\Delta x,\vartheta}(x_{1},y+b_{\Delta x,\vartheta}(x_{1}))-U_{\Delta x,\vartheta}(x_{2},y+b_{\Delta x,\vartheta}(x_{2}))|\,\dd y< L_{1}|x_{1}-x_{2}|,
\end{align*}
for some $L_{1}>0$ independent of $U_{\Delta x,\vartheta}$, $\Delta x$, and $\vartheta$.

\section{Convergence of the Approximate Solutions}\label{Section-final}
In \S \ref{Section-Decrease}, the uniform bound of the total variation of the approximate solutions $U_{\Delta x,\vartheta}$ has been obtained.
Then, by Propositions \ref{Prop:Estimate}--\ref{Prop:EstimateBdr},
the existence of convergent subsequences of the approximate solutions
$\{U_{\Delta x,\vartheta}\}$ follows.
Now we are going to prove that there is a convergent subsequence of
the approximate solutions $\{U_{\Delta x,\vartheta}\}$ whose limit
is an entropy solution to our problem.

\smallskip
Take $\Delta x=2^{-m}$, $m=0,1,2,\cdots$.
For any randomly chosen sequences
$\vartheta=(\vartheta_{0},\vartheta_{1},\vartheta_{2},\dots,\vartheta_{h},\dots)$,
we obtain a set of approximate solutions, which are denoted by $\{(u_m,v_m)\}$. It suffices to prove that there is a subsequence (still denoted by) $\{(u_m,v_m)\}$ such that
\begin{equation}
\iint_{\Omega_{\Delta x,\vartheta}}\Big(\phi_x\rho_m u_m+\phi_y\rho_m v_m
-\frac{\rho_m v_m\phi}{y}\Big)\,\dd x \dd y
+\int_{-\infty}^{y_{0}(0)}\phi(x_{0},y)\rho(x_{0},y)u(x_{0},y)\,\dd y\rightarrow0
\label{eqn:Conver1}
\end{equation}
for any $\phi(x,y)\in \text{C}_{0}^1(\mathbb{R}^2;\mathbb{R})$,
and
\begin{equation}
\iint_{\Omega_{\Delta x,\vartheta}}\left(\phi_x v_m-\phi_y u_m\right)\dd x \dd y\rightarrow0
\label{eqn:Conver2}
\end{equation}
for any $\phi(x,y)\in \text{C}_{0}^1(\Omega;\mathbb{R})$, as $m\rightarrow\infty$.
 We now prove (\ref{eqn:Conver1}) only, since (\ref{eqn:Conver2}) can be deduced analogously.\par

For simplicity,
we drop the subscript of $(u_m,v_m)$, and rewrite (\ref{eqn:Conver1}) as
\begin{align*}
&\iint_{\Omega_{\Delta x,\vartheta}}
\Big(\phi_x\rho u+\phi_y\rho v-\frac{\rho v\phi}{y}\Big)\,\dd x \dd y
+\int_{-\infty}^{y_{0}(0)}\phi(x_{0},y)\rho(x_{0},y)u(x_{0},y)\,\dd y\\
&=\sum_{h=1}^{\infty}\iint_{\Omega_{\Delta x,\vartheta,h}}
\Big(\phi_x\rho u+\phi_y\rho v-\frac{\rho v\phi}{y}\Big)\,\dd x \dd y
+\int_{-\infty}^{y_{0}(0)}\phi(x_{0},y)\rho(x_{0},y)u(x_{0},y)\,
\dd y.
\end{align*}
By the shock waves and the upper/lower edges of rarefaction waves,
each $\Omega_{\Delta x,\vartheta,h}$ can be divided into smaller polygons:
$\Omega_{\Delta x,\vartheta,h,j}$, $j=0,-1,-2,\cdots$, alternatively,
where $\Omega_{\Delta x,\vartheta,h,0}$ is the uppermost area below the approximate boundary $\Gamma_{\Delta x,\vartheta,h}$.
Then we have
\begin{equation}
\begin{aligned}
&\iint_{\Omega_{\Delta x,\vartheta}}
\Big(\phi_x\rho u+\phi_y\rho v-\frac{\rho v\phi}{y}\Big)\,\dd x\dd y
+\int_{-\infty}^{y_{0}(0)}\phi(x_{0},y)\rho(x_{0},y)u(x_{0},y)\,\dd y\\
&=\sum_{h=1}^{\infty}\sum_{j=0}^{-\infty}\iint_{\Omega_{\Delta x,\vartheta,h,j}}
\Big(\phi_x\rho u+\phi_y\rho v-\frac{\rho v\phi}{y}\Big)\,\dd x \dd y
+\int_{-\infty}^{y_{0}(0)}\phi(x_{0},y)\rho(x_{0},y)u(x_{0},y)\,\dd y\\
&=-\sum_{h,j}\iint_{\Omega_{\Delta x,\vartheta,h,j}}
\phi\Big((\rho u)_x+(\rho v)_y+\frac{\rho v}{y}\Big)\,\dd x \dd y
+\sum_{h,j}\iint_{\Omega_{\Delta x,\vartheta,h,j}}
\big(	(\phi\rho u)_x+(\phi\rho v)_y\big)\,\dd x \dd y\\
&\quad +\int_{-\infty}^{y_{0}(0)}\phi(x_{0},y)\rho(x_{0},y)u(x_{0},y)\,\dd y\\[1mm]
&=:\text{\uppercase\expandafter{\romannumeral1}}+\text{\uppercase\expandafter{\romannumeral2}}
+\text{\uppercase\expandafter{\romannumeral3}}.
\end{aligned}
\label{eqn:Conver1Prf}
\end{equation}

We first have

\smallskip
\Proposition\label{Prop:Conver1} \,\,
\uppercase
\expandafter{\romannumeral1}$\rightarrow0$ as $\Delta x\rightarrow 0$.

\smallskip
\noindent\textbf{Proof.}
To deal with the first term \uppercase\expandafter{\romannumeral1} in (\ref{eqn:Conver1Prf}), we use the transform:
\begin{equation}
\sigma=\frac{y}{x},\qquad \eta=\frac{y-y_n(h)}{x-x_{h}},
\label{eqn:Conver1trs}
\end{equation}
where $(x_{h},y_n(h))$ are the center of the Riemann problem,
and $n$ depends on $j$. Then we obtain
\begin{equation}
\begin{aligned}
\text{\uppercase\expandafter{\romannumeral1}}=&-\sum_{h,j}\iint_{\Omega_{\Delta x,\vartheta,h,j}}\frac{(\sigma-\eta)\phi\rho}{\sigma(y_n(h)-\eta x_{h})}
\Big(-\sigma^2\big(1-\frac{u^2}{c^2}\big) u_{\sigma}+\frac{2uv\sigma^2}{c^2}v_{\sigma}
+\big(1-\frac{v^2}{c^2}\big) v_{\sigma}\sigma+v\Big)\,\dd\eta \dd\sigma\\
&-\sum_{h,j}\iint_{\Omega_{\Delta x,\vartheta,h,j}}\frac{(\eta-\sigma)\phi}{\sigma  x_{h}-y_n(h)}
\big(-\eta(\rho u)_{\eta}+(\rho v)_{\eta}\big)\,\dd\eta \dd\sigma.
\end{aligned}
\label{eqn:Conver1ter1}
\end{equation}
From the construction of the approximate solutions, the first term of (\ref{eqn:Conver1ter1}) vanishes. For the second term, we have
\begin{align*}
-\eta(\rho u)_{\eta}+(\rho v)_{\eta}=O(1)\Delta\sigma,
\end{align*}
where $\Delta\sigma$ is the change of the $\sigma-$coordinate in domain $\Omega_{\Delta x,\vartheta,h,j}$.
Denote the rarefaction waves in $\Omega_{\Delta x,\vartheta,h}$ alternatively by $\alpha_{R,h,i}$. Then we have
\begin{align*}
\text{\uppercase\expandafter{\romannumeral1}}=O(1)\sum_{h,j}\Delta\eta(\Delta\sigma)^2,
\end{align*}
with $\Delta\eta=O(1)\alpha_{R,h,i}$.
According to Proposition \ref{Prop:Estimate},
the total strength $\sum_{i}|\alpha_{R,h,i}|$ of rarefaction waves in $\Omega_{\Delta x,\vartheta,h}$ is bounded, so that
\begin{equation}
\text{\uppercase\expandafter{\romannumeral1}}=O(1)\,\text{diam}(\text{supp}\phi)\,\Delta x,
\label{eqn:Conver1con1}
\end{equation}
which gives desired result.
\hfill{$\square$}

\medskip
Next, applying Green's formula in each $\Omega_{\Delta x,\vartheta,h,j}$,
we obtain
\begin{equation}
\begin{aligned}
\text{\uppercase\expandafter{\romannumeral2}}+\text{\uppercase\expandafter{\romannumeral3}}
&=\sum_{h=1}^{\infty}\int_{-\infty}^{b_{\Delta x,\vartheta}(x_{h})}\phi(x_{h},y)\big(\rho(x_{h}-,y)u(x_{h}-,y)-\rho(x_{h}+,y)u(x_{h}+,y)\big)\,\dd y\\
&\quad +\sum_{h=0}^{\infty}\int_{x_{h}}^{x_{h+1}}\phi(x,b(x))
\rho(x,b(x)) \big(v(x,b(x))-u(x,b(x))b'(x)\big)\dd x\\
&\quad +\sum_{h,i}\int_{W_{h,i}}
\big(s_{h,i}(\rho^+u^+-\rho^-u^-)-(\rho^+v^+-\rho^-v^-)\big)\phi\,\dd x\\
&=: \text{\uppercase\expandafter{\romannumeral4}}+\text{\uppercase\expandafter{\romannumeral5}}+\text{\uppercase\expandafter{\romannumeral6}},
\end{aligned}
\label{eqn:Conver1ter2}
\end{equation}
where
$W_{h,i}=\{(x,y)\,:\,y=w_{h,i}(x)=s_{h,i}(x-x_{h})+y_n(h)\,\text{ for some }n\}$ are shock waves or upper/lower edges of rarefaction waves
lying in $\Omega_{\Delta x,\vartheta,h}$,
and $\rho^\pm=\rho(x,w_{i,h}(x)\pm)$, $u^\pm=u(x,w_{i,h}(x)\pm)$, and $v^\pm=v(x,w_{i,h}(x)\pm)$.\par

\medskip
We now show

\smallskip
\Proposition\label{Prop:Conver2}
There exists
a subsequence of $\{(u_m,v_m)\}$ such that
\uppercase\expandafter{\romannumeral4}$\rightarrow0$ as $m\rightarrow \infty$.

\smallskip
\noindent\textbf{Proof.}
The first term on the right hand of (\ref{eqn:Conver1ter2}) can be rewritten as
\begin{align*}
\text{\uppercase\expandafter{\romannumeral4}}=\sum_{h\geq1}V_{h}
\end{align*}
with
\begin{align*}
V_{h}&=\sum_{n=n_{\chi,h}+1}^{n_{b,h}-1}\int_{y_{n-1}(h)}^{y_{n}(h)}\phi(x_{h},y)\big(\rho(x_{h}-,y)u(x_{h}-,y)-\rho(x_{h}+,y)u(x_{h}+,y)\big)\,\dd y\\
&\quad+\int_{\chi_{\Delta x,\vartheta}(x_{h})}^{y_{n_{\chi,h}}(h)+1}\phi(x_{h},y)\big(\rho(x_{h}-,y)u(x_{h}-,y)-\rho(x_{h}+,y)u(x_{h}+,y)\big)\,\dd y\\
&\quad+\int_{y_{n_{b,h}}(h)-1}^{b_{\Delta x,\vartheta}(x_{h})}\phi(x_{h},y)\big(\rho(x_{h}-,y)u(x_{h}-,y)-\rho(x_{h}+,y)u(x_{h}+,y)\big)\,\dd y.
\end{align*}
To show $\text{\uppercase\expandafter{\romannumeral4}}\rightarrow0$ for some subsequence $\{(u_m,v_m)\}$, we now introduce
\begin{align*}
\widetilde{V}=\sum_{h\geq1}\widetilde{V}_{h}
\end{align*}
with
\begin{align*}
\widetilde{V}_{h}
&=\sum_{n=n_{\chi,h}+1}^{n_{b,h}-1}\int_{y_{n-1}(h)}^{y_{n}(h)}\phi(x_{h},y_n(h))\big(\rho(x_{h}-,y)u(x_{h}-,y)-\rho(x_{h}+,y)u(x_{h}+,y)\big)\,\dd y\\
&\quad +\int_{y_{n_{b,h}}(h)}^{b_{\Delta x,\vartheta}(x_{h})}\phi(x_{h},y_{n_{b,h}+1}(h))\big(\rho(x_{h}-,y)u(x_{h}-,y)-\rho(x_{h}+,y)u(x_{h}+,y)\big)\,\dd y\\
&\quad+\int_{y_{n_{b,h}-1}(h)}^{y_{n_{b,h}}(h)}\phi(x_{h},y_{n_{b,h}}(h))\big(\rho(x_{h}-,y)u(x_{h}-,y)-\rho(x_{h}+,y)u(x_{h}+,y)\big)\dd y\\	
&\quad+\int_{y_{n_{\chi,h}}(h)}^{y_{n_{\chi,h}+1}(h)}\phi(x_{h},y_{n_{\chi,h}+1}(h))\big(\rho(x_{h}-,y)u(x_{h}-,y)-\rho(x_{h}+,y)u(x_{h}+,y)\big)\,\dd y\\	
&\quad+\int_{\chi_{\Delta x,\vartheta}(x_{h})}^{y_{n_{\chi,h}}(h)}\phi(x_{h},y_{n_{\chi,h}}(h))\big(\rho(x_{h}-,y)u(x_{h}-,y)-\rho(x_{h}+,y)u(x_{h}+,y)\big)\,\dd y\\
&=:\sum_{n=n_{\chi,h}+1}^{n_{b,h}-1}\int_{y_{n-1}(h)}^{y_{n}(h)}\phi(x_{h},y_n(h))\big(\rho(x_{h}-,y)u(x_{h}-,y)-\rho(x_{h}+,y)u(x_{h}+,y)\big)\,\dd y\\[1mm]
&\quad+\check{V}_{h}^{(0)}+\check{V}_{h}^{(1)}+\hat{V}_{h}^{(1)}+\hat{V}_{h}^{(0)}.
\end{align*}
From the construction of approximate solutions, we have
\begin{align*}
\check{V}_{h}^{(0)}&=O(1)\Delta x\big(|\alpha_{r,2}|+|\omega_{h+1}|+|\Delta\sigma_{b_{h}}|+Q(\Lambda_b)\big)&&\text{(see Case \ref{case:2}),}&\\
\hat{V}_{h}^{(0)}&=O(1)\Delta x\big(|\alpha_{l,1}|+|\Delta\sigma_{\chi_{h}}|+Q(\Lambda_{s})\big)&&\text{(see Case \ref{case:3}).}&
\end{align*}
From Proposition \ref{Prop:Estimate}, we obtain
\begin{equation}
\sum_{h=1}^{\infty}\big(\check{V}_{h}^{(0)}+\hat{V}_{h}^{(0)}\big)=O(1)\Delta xF(I_{1}).
\label{eqn:ConverTail}
\end{equation}
Then we write $\check{V}_{h}^{(1)}$ and $\hat{V}_{h}^{(1)}$ as
\begin{align*}
\check{V}_{h}^{(1)}
&=\int_{y_{n_{b,h}-1}(h)}^{y_{n_{b,h}}(h)}
\phi(x_{h},y_{n_{b,h}}(h))
\big(\rho(x_{h}-,y)u(x_{h}-,y)-\check{\rho}(x_{h}+,y)\check{u}(x_{h}+,y)\big)
\,\dd y\\
&\quad +\int_{y_{n_{b,h}-1}(h)}^{y_{n_{b,h}}(h)}
\phi(x_{h},y_{n_{b,h}}(h))\big(\check{\rho}(x_{h}+,y)\check{u}(x_{h}+,y)-\rho(x_{h}+,y)u(x_{h}+,y)\big)\,\dd y
\end{align*}
and
\begin{align*}
\hat{V}_{h}^{(1)}
&=\int_{y_{n_{\chi,h}}(h)}^{y_{n_{\chi,h}+1}(h)}
\phi(x_{h},y_{n_{\chi,h}}(h))
\big(\rho(x_{h}-,y)u(x_{h}-,y)-\hat{\rho}(x_{h}+,y)\hat{u}(x_{h}+,y)\big)\,\dd y\\
&\quad +\int_{y_{n_{\chi,h}}(h)}^{y_{n_{\chi,h}+1}(h)}
\phi(x_{h},y_{n_{\chi,h}}(h))\big(\hat{\rho}(x_{h}+,y)\hat{u}(x_{h}+,y)-\rho(x_{h}+,y)u(x_{h}+,y)\big)\,\dd y,
\end{align*}
where
\begin{align*}
&\check{U}(x_{h}+,y)= \tilde{U}(\dfrac{y}{x_{h}};\dfrac{r_{h,n_{b,h}-1}}{x_{h}},U(x_{h}+,r_{h,n_{b,h}-1})),\\
&\hat{U}(x_{h}+,y)=\tilde{U}(\dfrac{y}{x_{h}};\dfrac{r_{h,n_{\chi,h}-1}}{x_{h}},U(x_{h}+,r_{h,n_{\chi,h}-1})),
\end{align*}
and $\check{\rho}$ and $\hat{\rho}$ are determined via Bernoulli's equation.
By the construction of the approximate solutions near the boundary and near the leading shock, we have
\begin{align*}
&\int_{y_{n_{b,h}-1}(h)}^{y_{n_{b,h}}(h)}\phi\big(x_{h},y_{n_{b,h}}(h)\big)\big(\check{\rho}(x_{h}+,y)\check{u}(x_{h}+,y)-\rho(x_{h}+,y)u(x_{h}+,y)\big)\dd y\\
&\,\,\, =O(1)\Delta x\big(|\alpha_{r,2}|+|\omega_{h+1}|+|\Delta\sigma_{b_{h}}|+Q(\Lambda_b)\big)&&\text{(see Case \ref{case:2}),}&\\
&\int_{y_{n_{\chi,h}}(h)}^{y_{n_{\chi,h}+1}(h)}\phi\big(x_{h},y_{n_{\chi,h}}(h)\big)\big(\hat{\rho}(x_{h}+,y)\hat{u}(x_{h}+,y)-\rho(x_{h}+,y)u(x_{h}+,y)\big)\dd y\\
&\,\,\, =O(1)\Delta x\big(|\alpha_{l,1}|+|\Delta\sigma_{\chi_{h}}|+Q(\Lambda_{s})\big) &&\text{(see Case \ref{case:3}).}&
\end{align*}
Similarly, by Proposition \ref{Prop:Estimate}, we conclude
\begin{equation}
\begin{aligned}
&\sum_{h=1}^{\infty}\int_{y_{n_{b,h}-1}(h)}^{y_{n_{b,h}}(h)}\phi(x_{h},y_{n_{b,h}}(h))\big(\check{\rho}(x_{h}+,y)\check{u}(x_{h}+,y)
-\rho(x_{h}+,y)u(x_{h}+,y)\big)\,\dd y
=O(1)\Delta xF(I_{1}),\\
&\sum_{h=1}^{\infty}\int_{y_{n_{\chi,h}}(h)}^{y_{n_{\chi,h}+1}(h)}
\phi(x_{h},y_{n_{\chi,h}}(h))\big(\hat{\rho}(x_{h}+,y)\hat{u}(x_{h}+,y)-\rho(x_{h}+,y)u(x_{h}+,y)\big)\,\dd y
=O(1)\Delta xF(I_{1}).
\end{aligned}
\label{eqn:ConverTailre}
\end{equation}
Setting
\begin{align*}
\bar{V}_{h}
&= \sum_{n=n_{\chi,h}+1}^{n_{b,h}-1}\int_{y_{n-1}(h)}^{y_{n}(h)}
\phi(x_{h},y_n(h))\big(\rho(x_{h}-,y)u(x_{h}-,y)
-\rho(x_{h}+,y)u(x_{h}+,y)\big)\,\dd y\\
&\quad +\int_{y_{n_{b,h}-1}(h)}^{y_{n_{b,h}}(h)}
\phi(x_{h},y_{n_{b,h}}(h))
\big(\rho(x_{h}-,y)u(x_{h}-,y)
-\check{\rho}(x_{h}+,y)\check{u}(x_{h}+,y)\big)\,\dd y\\
&\quad +\int_{y_{n_{\chi,h}}(h)}^{y_{n_{\chi,h}+1}(h)}
\phi(x_{h},y_{n_{\chi,h}}(h))\big(\rho(x_{h}-,y)u(x_{h}-,y)
-\hat{\rho}(x_{h}+,y)\hat{u}(x_{h}+,y)\big)\,\dd y.
\end{align*}

As in \cite{Glimm1965} (see also \cite{Chen2004}),
let
$$
H=\prod_{h=0}^{\infty}[0,1)
=\{\vartheta=(\vartheta_{0},\vartheta_{1},\vartheta_{2},\dots,
\vartheta_{h},\dots)\,:\,\vartheta_{h}\in[0,1),\ h=0,1,2,\cdots\}.
$$
Denoting $\bar{y}=y_{n-1}(h)+\vartheta_{h}\big(y_{n}(h)-y_{n-1}(h)\big)$,
we obtain from (\ref{eqn:Formula2}) that
\begin{align*}
&\rho(x_{h}-,y)u(x_{h}-,y)-\rho(x_{h}+,y)u(x_{h}+,y)\\
&=\rho(x_{h}-,y)u(x_{h}-,y)
-\rho(x_{h}-,\bar{y})u(x_{h}-,\bar{y})+\rho(x_{h}+,\bar{y})u(x_{h}+,\bar{y})
-\rho(x_{h}+,y)u(x_{h}+,y)\\
&=O(1)|\alpha|+O(1)|\alpha||\Delta\sigma_\alpha|+O(1)|\Delta\sigma_{\alpha}|\\
&=O(1)(|\alpha|+|\Delta\sigma_\alpha|),
\end{align*}
where $\alpha$ is an elementary wave in $\Omega_{\Delta x,\vartheta,h,j}$ and $\Delta\sigma_{\alpha}$
is the change of the $\sigma-$coordinate in the elementary wave $\alpha$.
Denote the elementary waves in $\Omega_{\Delta x,\vartheta,h}$ by $\alpha_{h,i}$. Then we
\begin{equation}
\bar{V}_{h}=O(1)\Big(\sum_{i\leq0}|\alpha_{h,i}|+\sigma^*-\sigma_*\Big)\Delta x,
\label{eqn:Conver2est1}
\end{equation}
which implies
\begin{align*}
\sum_{h\geq1}\int_{H}\bar{V}_{h}^2\dd\vartheta=O(1)\text{diam}(\text{supp}\phi)\Big(\sum_{i\leq0}|\alpha_{h,i}|+\sigma^*-\sigma_*\Big)^2\Delta x.
\end{align*}

\smallskip
Next, we need the following lemma.

\smallskip
\Lemma\label{Lem:ConvRan} The approximate solutions
$\{U_{\Delta x,\vartheta}(x,y)\}$ satisfy
\begin{equation}
\int_{0}^1\int_{y_{n-1}(h)}^{y_n(h)}
\big(U_{\Delta x,\vartheta}(x_{h}-,y)-U_{\Delta x,\vartheta}(x_{h}+,y)\big)\,
\dd y \dd\vartheta_{h}=O(1)(\Delta x)^3+O(1)(|\alpha|+|\beta|)(\Delta x)^2.
\label{eqn:Conver1Lem}
\end{equation}

\noindent\textbf{Proof.} We now give a proof when $\alpha$ and $\beta$ are both shock waves, since the remaining cases can be obtained similarly.

Suppose that $\alpha$ and $\beta$ issue from $(x_{h-1},y_{n-1}(h-1))$
and $(x_{h-1},y_{n}(h-1))$,
and end at $(x_{h},r_{1})$ and $(x_{h},r_{2})$, respectively. Set $a_{1}=\frac{r_{1}-y_{n-1}(h)}{y_{n}(h)-y_{n-1}(h)}$
and $a_{2}=\frac{r_{2}-y_{n-1}(h)}{y_{n}(h)-y_{n-1}(h)}$.
From the construction of approximate solutions, we have
\begin{align*}
&\int_{0}^1\int_{y_{n-1}(h)}^{y_n(h)}
\big(U_{\Delta x,\vartheta}(x_{h}-,y)-U_{\Delta x,\vartheta}(x_{h}+,y)\big)\,\dd y \dd\vartheta_{h}\\
&=\int_{0}^1\int_{y_{n-1}(h)}^{y_n(h)}U_{\Delta x,\vartheta}(x_{h}-,y)\,\dd y \dd\vartheta_{h}-
\int_{0}^1\int_{y_{n-1}(h)}^{y_n(h)}U_{\Delta x,\vartheta}(x_{h}+,y)\,\dd y \dd\vartheta_{h}\\
&=\int_{y_{n-1}(h)}^{r_{1}}\tilde{U}(\frac{y}{x_{h}};\sigma_{2},U_{l})\,\dd y
+\int_{r_{1}}^{r_{2}}\tilde{U}(\frac{y}{x_{h}};\sigma_{2},\Phi(0,\alpha_{2};U_{l}))\,\dd y\\
&\quad +\int_{r_{2}}^{y_n(h)}\tilde{U}(\frac{y}{x_{h}};\sigma_{1},\Phi(\beta_{1},0;\tilde{U}(\sigma_{1};\sigma_{2},\Phi(0,\alpha_{2};U_{l}))))\,\dd y
-a_{1}\int_{y_{n-1}(h)}^{y_n(h)}\tilde{U}(\frac{y}{x_{h}};\sigma_{2},U_{l})\,\dd y\\
&\quad -(a_{2}-a_{1})\int_{y_{n-1}(h)}^{y_n(h)}\tilde{U}(\frac{y}{x_{h}};\sigma_{2},\Phi(0,\alpha_{2};U_{l}))\,\dd y\\
&\quad -(1-a_{2})\int_{y_{n-1}(h)}^{y_n(h)}\tilde{U}(\frac{y}{x_{h}};\sigma_{1},\Phi(\beta_{1},0;\tilde{U}(\sigma_{1};\sigma_{2},\Phi(0,\alpha_{2};U_{l}))))\,\dd y.
\end{align*}
Since $\tilde{U}(\frac{y}{x_{h}};\sigma_{2},\Phi(0,\alpha_{2};U_{l}))
=\tilde{U}(\frac{y}{x_{h}};\tilde{U}(\sigma_{1};\sigma_{2},\Phi(0,\alpha_{2};U_{l})))$, we obtain
\begin{align*}
&\int_{0}^1\int_{y_{n-1}(h)}^{y_n(h)}(U_{\Delta x,\vartheta}(x_{h}-,y)-U_{\Delta x,\vartheta}(x_{h}+,y))\,\dd y \dd\vartheta_{h}\\
&=\int_{y_{n-1}(h)}^{r_{1}}(1-a_{1})\big(\tilde{U}(y/x_{h};\sigma_{2},U_{l})-\tilde{U}(\frac{y}{x_{h}};\sigma_{2},\Phi(0,\alpha_{2};U_{l}))\big)\,\dd y\\
&\quad -\int_{y_{n-1}(h)}^{r_{1}}(1-a_{2})\big(\tilde{U}(\frac{y}{x_{h}};\sigma_{1},\Phi(\beta_{1},0;\tilde{U}(\sigma_{1};\sigma_{2},\Phi(0,\alpha_{2};U_{l}))))
-\tilde{U}(\frac{y}{x_{h}};\tilde{U}(\sigma_{1};\sigma_{2},\Phi(0,\alpha_{2};U_{l})))\big)\,\dd y\\
&\quad +\int_{r_{1}}^{r_{2}}a_{1}\big(\tilde{U}(\frac{y}{x_{h}};\sigma_{2},\Phi(0,\alpha_{2};U_{l}))-\tilde{U}(\frac{y}{x_{h}};\sigma_{2},U_{l})\big)\,\dd y\\
&\quad -\int_{r_{1}}^{r_{2}}(1-a_{2})\big(\tilde{U}(\frac{y}{x_{h}};\sigma_{1},\Phi(\beta_{1},0;\tilde{U}(\sigma_{1};\sigma_{2},\Phi(0,\alpha_{2};U_{l}))))
-\tilde{U}(\frac{y}{x_{h}};\tilde{U}(\sigma_{1};\sigma_{2},\Phi(0,\alpha_{2};U_{l})))\big)\,\dd y\\
&\quad +\int_{r_{2}}^{y_n(h)}a_{2}\big(\tilde{U}(\frac{y}{x_{h}};\sigma_{1},\Phi(\beta_{1},0;\tilde{U}(\sigma_{1};\sigma_{2},\Phi(0,\alpha_{2};U_{l}))))
-\tilde{U}(\frac{y}{x_{h}};\tilde{U}(\sigma_{1};\sigma_{2},\Phi(0,\alpha_{2};U_{l})))\big)\, \dd y\\
&\quad -\int_{r_{2}}^{y_n(h)}a_{1}\big(\tilde{U}(\frac{y}{x_{h}};\sigma_{2},U_{l})-\tilde{U}(\frac{y}{x_{h}};\sigma_{2},\Phi(0,\alpha_{2};U_{l}))\big)\,\dd y.
\end{align*}
Then, by Taylor's expansion, we have
\begin{align*}
&\tilde{U}(\frac{y}{x_{h}};\sigma_{2},U_{l})
-\tilde{U}(\frac{y}{x_{h}};\sigma_{2},\Phi(0,\alpha_{2};U_{l}))\\
&\,\,\,=U_{l}-\Phi(0,\alpha_{2};U_{l})+A_{1}(y-y_{n-1}(h))+O(1)(y-y_{n-1}(h))^2,
\\[2mm]
&\tilde{U}(\frac{y}{x_{h}};\sigma_{1},\Phi(\beta_{1},0;\tilde{U}(\sigma_{1};\sigma_{2},\Phi(0,\alpha_{2};U_{l}))))
-\tilde{U}(\frac{y}{x_{h}};\tilde{U}(\sigma_{1};\sigma_{2},\Phi(0,\alpha_{2};U_{l})))\\
&\,\,\,=\Phi(\beta_{1},0;\tilde{U}(\sigma_{1};\sigma_{2},\Phi(0,\alpha_{2};U_{l})))
-\tilde{U}(\sigma_{1};\sigma_{2},\Phi(0,\alpha_{2};U_{l}))+A_{2}(y-y_{n}(h))+O(1)(y-y_{n}(h))^2,
\end{align*}
with
\begin{align*}
A_{1}=&\left. \partial_y\Big(\tilde{U}(\frac{y}{x_{h}};\sigma_{2},U_{l})-\tilde{U}(\frac{y}{x_{h}};\sigma_{2},\Phi(0,\alpha_{2};U_{l}))\Big)
\right|_{y=y_{n-1}(h)},\\
A_{2}=&\left. \partial_y\Big(\tilde{U}(\frac{y}{x_{h}};\sigma_{1},\Phi(\beta_{1},0;\tilde{U}(\sigma_{1};\sigma_{2},\Phi(0,\alpha_{2};U_{l}))))
-\tilde{U}(\frac{y}{x_{h}};\tilde{U}(\sigma_{1};\sigma_{2},\Phi(0,\alpha_{2};U_{l})))\Big)
\right|_{y=y_{n}(h)}.
\end{align*}
A direct computation leads to
\begin{align*}
&\int_{0}^1\int_{y_{n-1}(h)}^{y_n(h)}
\big(U_{\Delta x,\vartheta}(x_{h}-,y)-U_{\Delta x,\vartheta}(x_{h}+,y)\big)\,\dd y \dd\vartheta_{h}\\
&= O(1)(y_{n}(h)-y_{n-1}(h))^3+\dfrac{1}{2}A_{1}(r_{1}-y_{n-1}(h))(r_{1}-y_{n}(h))-\dfrac{1}{2}A_{2}(r_{2}-y_{n-1}(h))(r_{2}-y_{n}(h)).
\end{align*}
Noting that $A_{1}=O(1)|\alpha|$ and $A_{2}=O(1)|\beta|$, together with the Courant-Friedrichs-Lewy condition, we conclude  (\ref{eqn:Conver1Lem}).\hfill{$\square$}

\medskip	
Substituting $U_{\Delta x,\vartheta}$ in Lemma \ref{Lem:ConvRan} by $\rho u$ and carrying out the same process lead to
\begin{align*}
\int_{0}^1\bar{V}_{h}\dd\vartheta_{h}
=O(1)\Big(\text{diam(supp$\phi$)}+\sum_{i\leq0}|\alpha_{h,i}|\Big)(\Delta x)^2=O(1)(\Delta x)^2.
\end{align*}
As in (\ref{eqn:Conver2est1}), we obtain
\begin{align*}
\bar{V}_k=O(1)\Big(\sum_{i\leq0}|\alpha_{k,i}|+\sigma^*-\sigma_*\Big)\Delta x.
\end{align*}
Then
\begin{align*}
\sum_{h>k}\int_{H}\bar{V}_{h}	\bar{V}_k\,\dd\vartheta
\leq&\sum_{h>k}\left|\int_{0}^1	\bar{V}_{h}d\vartheta_{h}\right|\int_{0}^1|	\bar{V}_k|\,\dd \hat{\vartheta}_{h}
=O(1)\big(\text{diam(supp$\phi$)}\big)^2\Delta x,
\end{align*}
where $\dd\hat{\vartheta}_{h}
=\dd\vartheta_{0}\cdots \dd\vartheta_{h-1}\dd\vartheta_{h+1}\cdots$.\par

Since
\begin{align*}
\|\bar{V}\|_{L^2(H)}
=\sum_{h\geq1}\int_{H}\bar{V}_{h}^2\dd\vartheta
+2\sum_{h>k}\int_{H}\bar{V}_{h}\bar{V}_k\,\dd\vartheta,
\end{align*}
we conclude that
\begin{align*}
\|\bar{V}\|_{L^2(H)}\rightarrow0\qquad \text{as }\Delta x\rightarrow0,
\end{align*}
which, combining with (\ref{eqn:ConverTail})--(\ref{eqn:ConverTailre}), gives a subsequence (still denoted by) $\{(u_m,v_m)\}$ such that $\widetilde{V}\rightarrow0$ almost everywhere.
Meanwhile, we have
\begin{align*}
V_{h}-\widetilde{V}_{h}
&= \sum_{n=n_{\chi,h}+1}^{n_{b,h}-1}\int_{y_{n-1}(h)}^{y_n(h)}\big(\phi(x_{h},y_n(h))-\phi(x_{h},y)\big)\big(\rho(x_{h}-,y)u(x_{h}-,y)-\rho(x_{h}+,y)u(x_{h}+,y)\big)\,\dd y\\
&\quad +\int_{y_{n_{b,h}-1}(h)}^{y_{n_{b,h}}(h)}
\big(\phi(x_{h},y_{n_{b,h}}(h))-\phi(x_{h},y)\big)
\big(\rho(x_{h}-,y)u(x_{h}-,y)-\rho(x_{h}+,y)u(x_{h}+,y)\big)\,\dd y\\	
&\quad+\int_{y_{n_{b,h}}(h)}^{b_{\Delta x,\vartheta}(x_{h})}
\big(\phi(x_{h},y_{n_{b,h}+1}(h))-\phi(x_{h},y)\big)\big(\rho(x_{h}-,y)u(x_{h}-,y)-\rho(x_{h}+,y)u(x_{h}+,y)\big)\,\dd y\\
&\quad +\int_{\chi_{\Delta x,\vartheta}(x_{h})}^{y_{n_{\chi,h}}(h)}
\big(\phi(x_{h},y_{n_{\chi,h}}(h))-\phi(x_{h},y)\big)\big(\rho(x_{h}-,y)u(x_{h}-,y)-\rho(x_{h}+,y)u(x_{h}+,y)\big)\,\dd y\\
&\quad +\int_{y_{n_{\chi,h}}(h)}^{y_{n_{\chi,h}+1}(h)}
\big(\phi(x_{h},y_{n_{\chi,h}+1}(h))-\phi(x_{h},y)\big)\big(\rho(x_{h}-,y)u(x_{h}-,y)-\rho(x_{h}+,y)u(x_{h}+,y)\big)\,\dd y\\	
&=O(1)\Delta x\sum_{n=n_{\chi,h}+1}^{n_{b,h}-1}\int_{y_{n-1}(h)}^{y_n(h)}\big(\rho(x_{h}-,y)u(x_{h}-,y)-\rho(x_{h}+,y)u(x_{h}+,y)\big)\,\dd y
+O(1)(\Delta x)^2\\
&=O(1)\Big(\sum_{i\leq0}|\alpha_{h,i}|+\sigma^*-\sigma_*+1\Big)(\Delta x)^2,
\end{align*}
which leads to
\begin{align*}
\text{\uppercase\expandafter{\romannumeral4}}-\widetilde{V}=&\sum_{h\geq1}V_{h}-\widetilde{V}_{h}
=O(1)\,\text{diam}(\text{supp}\phi)\Delta x.
\end{align*}
Thus, $\text{\uppercase\expandafter{\romannumeral4}}\rightarrow0$
as $m\rightarrow\infty$ for some subsequence $\{(u_m,v_m)\}$.
\hfill{$\square$}

\medskip
\Proposition\label{Prop:Conver3}
\,\,\uppercase\expandafter{\romannumeral5}, \uppercase\expandafter{\romannumeral6}$\rightarrow0$ \, as $\Delta x\rightarrow 0$.

\smallskip
\noindent\textbf{Proof.} Since
\begin{align*}
b'(x)=\frac{v\big(x_{h}+,b(x_{h})-\big)}{u\big(x_{h}+,b(x_{h})-\big)}\qquad \text{for }x\in(x_{h},x_{h+1}),
\end{align*}
it follows from the construction of our approximate solution that
\begin{align*}
v\big(x,b(x)\big)-u\big(x,b(x)\big)b'(x)=O(1)\Delta x.
\end{align*}
Therefore, we have
\begin{align*}
\text{\uppercase\expandafter{\romannumeral5}}
=O(1)\,\text{diam}(\text{supp}\phi)\,\Delta x\rightarrow0
\qquad\text{as }\Delta x\rightarrow0.
\end{align*}

As for \uppercase\expandafter{\romannumeral6}, we divide this term into three parts.
The first part is the integral along the leading shock, where $W_{h,i}=S_{\Delta x,\vartheta,h}$.
For this part, by similar arguments in treating \uppercase\expandafter{\romannumeral5}, we have
\begin{align*}
\sum_{h}\int_{S_{\Delta x,\vartheta,h}}[s_{h}(\rho^+u^+-\rho^-u^-)-(\rho^+v^+-\rho^-v^-)]\phi\,\dd x
=O(1)\Delta x.
\end{align*}
The second part is the integral along the upper or lower edges of rarefaction waves and therefore vanishes automatically.
The third part is the integral along the weak shock waves, that is, $W_{h,i}\neq S_{\Delta x,\vartheta,h}$.
In this case, by (\ref{eqn:Formula1}), we have
\begin{align*}
&\rho^+u^+-\rho^-u^-=(\rho^+u^+-\rho^-u^-)|_{x=x_{h}+}+O(1)(\rho^+u^+-\rho^-u^-)|_{x=x_{h}+}\Delta x,\\
&\rho^+v^+-\rho^-v^-=(\rho^+v^+-\rho^-v^-)|_{x=x_{h}+}+O(1)(\rho^+v^+-\rho^-v^-)|_{x=x_{h}+}\Delta x.
\end{align*}
Thus, in view of the Rankine-Hugoniot conditions, we obtain
\begin{align*}
\sum_{i}s_{h,i}(\rho^+u^+-\rho^-u^-)-(\rho^+v^+-\rho^-v^-)=O(1)\sum_{i}|\alpha_{S,h,i}|\Delta x,
\end{align*}
where $\alpha_{S,h,i}$ are the weak shock waves in
$\Omega_{\Delta x,\vartheta,h}$. Combining all the three parts together,
we have
\begin{align*}
\text{\uppercase\expandafter{\romannumeral4}}
=O(1)\,\text{diam}(\text{supp}\phi)\,\sum_{i}|\alpha_{S,h,i}|\Delta x+O(1)\Delta x.
\end{align*}
By Proposition \ref{Prop:Estimate}, $\sum_{i}|\alpha_{S,h,i}|$ is uniformly bounded with respect to $h$.
Therefore, $\text{\uppercase\expandafter{\romannumeral4}}\rightarrow0$
as $\Delta x\rightarrow0$.\hfill{$\square$}\\

With all the arguments stated above, a standard procedure as in \cite{Chen2004,Liu1999} gives the following main theorem.

\smallskip
\Theorem\label{Thm:main}
Suppose that (\ref{Asum:pres})--(\ref{Asum:super}), $1<\gamma<3$, and $0<p_{0}<p^*$.
Then, when $M_{\infty}$ is sufficiently large, there are $\epsilon_0>0$
and  a null set $\mathcal{N}$ such that,
if T.V.\,$p^{b}<\epsilon_0$, for each $\vartheta\in\prod_{h=0}^{\infty}[0,1)\backslash\mathcal{N}$, there exist both a subsequence
$\{\Delta_{i}\}_{i=0}^{\infty}\subset\{\Delta x\}$ of the mesh size with $\Delta_{i}\rightarrow0$ as $i\rightarrow\infty$
and a triple of functions $b_{\vartheta}(x)$ with $b_{\vartheta}(0)=0$, $\chi_{\vartheta}(x)$ with $\chi_{\vartheta}(0)=0$,
and $U_{\vartheta}(x,y)\in O_{\varepsilon_{0}}\big(G(s_{0})\cap W(p_{0},p_{\infty})\big)$ such that
\newcounter{Man}
\begin{list}{(\roman{Man})}{\usecounter{Man}}
\item $b_{\Delta_{i},\vartheta}(x)$ converges to $b_{\vartheta}(x)$ uniformly in any bounded $x$-interval;
\item $\chi_{\Delta_{i},\vartheta}(x)$ converges to $\chi_{\vartheta}(x)$ uniformly in any bounded $x$-interval;
\item $b'_{\Delta_{i},\vartheta}(x)$ converges to $b'_{\vartheta}(x)\in BV([0,\infty))$ a.e.
 and $\dfrac{\dd b_{\vartheta}(x)}{\dd x}=b'(x)$, a.e. with $|b'_{\vartheta}-b_0|<C\epsilon_0$;
\item $s_{\Delta_{i},\vartheta}(x)$ converges to $s_{\vartheta}(x)\in BV([0,\infty))$ a.e.
and $\dfrac{\dd\chi_{\vartheta}(x)}{\dd x}=\sigma_{\vartheta}(x)$, a.e. with $|s_{\vartheta}-s_0|<C\epsilon_0$;
\item $U_{\Delta_{i},\vartheta}(x,\cdot)$ converges to $U_{\vartheta}\in\textbf{L}^1_{\text{loc}}(-\infty,b_{\vartheta}(x))$ for every $x>0$,
and $U_{\vartheta}$ is a global entropy solution of the inverse problem (\ref{eqn:orig})--(\ref{eqn:Bernoulli}) and
satisfies (\ref{eqn:Boundary})--(\ref{eqn:Cauchyda}).
\end{list}

\section{Asymptotic Behavior of Global Entropy Solutions}\label{Section-asymp}

To establish the asymptotic behavior of global entropy solutions,
we need further estimates of the approximate solutions.

\smallskip
\Lemma\label{Lem:aysmest}
There exists a constant $M_{1}$, independent of $U_{\Delta x,\vartheta}$, $\Delta x$, and $\vartheta$, such that
\begin{equation}
\sum_{\Lambda}E_{\Delta x,\vartheta}(\Lambda)<M_{1}
\label{eqn:aysmest1}
\end{equation}
for $E_{\Delta x,\vartheta}(\Lambda)$ given as in (\ref{eqn:DiffGLimmFunc}).

\smallskip
\noindent\textbf{Proof.}
By Proposition \ref{Prop:Estimate}, for any interaction
region $\Lambda\subset\{(h-1)\Delta x\leq (h+1)\Delta x\}$ for $h\geq1$,
we have
\begin{align*}
\sum_{\Lambda}E_{\Delta x,\vartheta}(\Lambda)\leq 4\sum_{\Lambda}\big(F(I)-F(J)\big)\leq 4F(I_{1}).
\end{align*}
Thus, choosing $M_{1}=4F(I_{1})+1$, the proof is complete.\hfill{$\square$}

\medskip
For any $t>0$, let $\mathcal{L}_{j,\vartheta}(t-)$, $j=1,2$, be the total variation of $j-$weak waves
in $U_{\vartheta}$ crossing line $x=t$, and let $\mathcal{L}_{j,\Delta x, \vartheta}(t-)$, $j=1,2$, be the total variation of $j-$weak waves
in $U_{\Delta x,\vartheta}$ crossing line $x=t$. Then we have

\smallskip
\Lemma\label{Lem:aysmest1} As $x\rightarrow\infty$,
\begin{align*}
\sum_{j=1}^{2}\mathcal{L}_{j,\vartheta}(x-)\rightarrow0.
\end{align*}

\noindent\textbf{Proof.}
Let $U_{\Delta_{i},\vartheta}$ be a sequence of the approximate solutions introduced in Theorem \ref{Thm:main},
and let the corresponding term $E_{\Delta x,\vartheta}(\Lambda)$ be defined in (\ref{eqn:DiffGLimmFunc}).
As in \cite{Glimm1970}, denoted by $\dd E_{\Delta x,\vartheta}$
the measure of assigning quantities $E_{\Delta x,\vartheta}(\Lambda)$ to the center of $\Lambda$.
Then, by Lemma \ref{Lem:aysmest}, we can choose a subsequence
(still denoted as) $\dd E_{\Delta_{i},\vartheta}$ such that
\begin{align*}
\dd E_{\Delta_{i},\vartheta}\rightarrow dE_{\vartheta}
\qquad \text{as } \Delta_{i}\rightarrow0
\end{align*}
with $E_{\vartheta}(\Lambda)<\infty$.\par

Therefore, for $\epsilon_{1}>0$ sufficiently small, we can choose $x_{\epsilon_{1}}$ (independent of $U_{\Delta_{i},\vartheta}$),
$\Delta_{i}$, and $\vartheta$ such that
\begin{align*}
\sum_{h>[x_{\epsilon_{1}}/\Delta x]} E_{\Delta_{i},\vartheta}(\Lambda_{h,n})<\epsilon_{1}.
\end{align*}
Let $X_{\epsilon_{1}}^{1}=(x_{\epsilon_{1}},\chi_{\Delta_{i},\vartheta}(x_{\epsilon_{1}}))$
and $X_{\epsilon_{1}}^{2}=(x_{\epsilon_{1}},b_{\Delta_{i},\vartheta}(x_{\epsilon_{1}}))$ be the two points lying in the approximate
leading shock $y=\chi_{\Delta_{i},\vartheta}(x)$ and the approximate boundary $y=b_{\Delta_{i},\vartheta}(x)$, respectively.
Let $\chi_{\Delta_{i},\vartheta}^{j}$ be the approximate $j$--generalized characteristic issuing from $X_{\epsilon_{1}}^{j}$ for $j=1,2$, respectively.
According to the construction of the approximate solutions,
there exist constants $\hat{M}_{j}>0$, $j=1,2$, independent of $U_{\Delta_{i},\vartheta}$, $\Delta_{i}$, and $\vartheta$,  such that
\begin{align*}
|\chi_{\Delta_{i},\vartheta}^{j}(x_{1})-\chi_{\Delta_{i},\vartheta}^{j}(x_{2})|\leq \hat{M}_{j}\big(|x_{1}-x_{2}|+\Delta_{i}\big)
\qquad \text{for } x_{1},x_{2}>x_{\epsilon_{1}}.
\end{align*}
Then we choose a subsequence (still denoted by) $\Delta_{i}$ such that
\begin{align*}
\chi_{\Delta_{i},\vartheta}^{j}\rightarrow\chi_{\vartheta}^{j}\qquad \text{as } \Delta_{i}\rightarrow0
\end{align*}
for some $\chi_{\vartheta}^{j}\in$Lip with $(\chi_{\vartheta}^{j})'$ bounded.\par

Let two characteristics $\chi_{\vartheta}^{1}$ and $\chi_{\vartheta}^{2}$ intersect with the cone boundary $\Gamma_{\vartheta}$
and the leading shock $S_{\vartheta}$ at points
$(t_{\epsilon_{1}}^1,\chi_{\vartheta}^{1}(t_{\epsilon_{1}}^1))$ and $(t_{\epsilon_{1}}^2,\chi_{\vartheta}^{2}(t_{\epsilon_{1}}^2))$
 for some $t_{\epsilon_{1}}^1$ and $t_{\epsilon_{1}}^2$, respectively.
Then, as in \cite{Glimm1970}, we apply the approximate conservation law
to the domain below $\chi_{\Delta_{i},\vartheta}^{1}$
and above $\chi_{\Delta_{i},\vartheta}^{1}$ and use Lemma \ref{Lem:aysmest} to obtain
\begin{align*}
\mathcal{L}_{j,\Delta_{i},\vartheta}(x-)\leq C\sum_{h>[x_{\epsilon_{1}}/\Delta x]} E_{\Delta_{i},\vartheta}(\Lambda_{h,n})<C\epsilon_{1}
\end{align*}
for $j=1,2$, $x>t_{\epsilon_{1}}^1+t_{\epsilon_{1}}^2$. This completes the proof.\hfill{$\square$}\\

\Theorem\label{Thm:asymp} For $p^{b}_{\infty}:=\lim_{x\rightarrow\infty}p^{b}(x)$, $s_{\infty}:=\lim_{x\rightarrow\infty}s_{\vartheta}(x)$,
and $b'_{\infty}=\lim_{x\rightarrow \infty}(b_{\vartheta})'_{+}(x)$,
\begin{align*}
&\lim_{x\rightarrow\infty}\sup\left\{\big|U_{\vartheta}(x,y)-\tilde{U}(\sigma;s_{\infty},G(s_{\infty}))\big|\,:\,\chi_{\vartheta}(x)<y<b_{\vartheta}(x)\right\}=0,\\
&\dfrac{1}{2}\big|\tilde{U}(b_{\infty}';s_{\infty},G(s_{\infty}))\big|^2+\dfrac{\gamma (p^{b}_{\infty})^{\frac{\gamma-1}{\gamma}}}{\gamma-1}=\dfrac{1}{2}+\dfrac{\gamma p_{\infty}^{\frac{\gamma-1}{\gamma}}}{\gamma-1},\\
&\tilde{U}(b_{\infty}';s_{\infty},G(s_{\infty}))\cdot(-b_{\infty}',1)=0.
\end{align*}

\noindent\textbf{Proof.} For every $x\in[x_{k-1},x_{k})$, we have
\begin{align*}
&\big|U_{\vartheta}(x,y)-\tilde{U}(\sigma;s_{\Delta_{i},\vartheta},G(s_{\Delta_{i},\vartheta}))\big|
+\big|\tilde{U}(b_{\Delta_{i},\vartheta}';s_{\Delta_{i},\vartheta},G(s_{\Delta_{i},\vartheta}))\cdot(-b_{\Delta_{i},\vartheta}',1)\big|\\
&\quad +\Big|\dfrac{1}{2}\big|\tilde{U}(b_{\Delta_{i},\vartheta}';s_{\Delta_{i},\vartheta},G(s_{\Delta_{i},\vartheta}))\big|^2
+\dfrac{\gamma (p^{b}_{\Delta x,k})^{\frac{\gamma-1}{\gamma}}}{\gamma-1}-\dfrac{1}{2}-\dfrac{\gamma p_{\infty}^{\frac{\gamma-1}{\gamma}}}{\gamma-1}\Big|\\
&\leq C \Big(\sum_{j=1}^{2}\mathcal{L}_{j,\Delta_{i},\vartheta}(x-)+|\Delta_{i}|\Big).
\end{align*}
By Theorem \ref{Thm:main}, letting $i\rightarrow\infty$, we obtain
\begin{align*}
&\sup_{\chi_{\vartheta}(x)<y<b_{\vartheta}(x)}
\big|U_{\vartheta}(x,\cdot)-\tilde{U}\big(\sigma;s_{\vartheta},G(s_{\vartheta})\big)\big|
+\big|\tilde{U}((b_{\vartheta})_{+}';s_{\vartheta},G(s_{\vartheta}))\cdot(-b_{\vartheta}',1)\big|\\
&\quad +\Big|\dfrac{1}{2}\big|\tilde{U}((b_{\vartheta})_{+}';s_{\vartheta},G(s_{\vartheta}))\big|^2+\dfrac{\gamma (p^{b})^{\frac{\gamma-1}{\gamma}}}{\gamma-1}-\dfrac{1}{2}-\dfrac{\gamma p_{\infty}^{\frac{\gamma-1}{\gamma}}}{\gamma-1}\Big|\\
&\leq C \sum_{j=1}^{2}\mathcal{L}_{j,\vartheta}(x-).
\end{align*}
Then, using Lemma \ref{Lem:aysmest1} and noting that
$\tilde{U}(\sigma;s,G(s))$ is a continuous function with respect
to $\sigma$ and $s$, we conclude our result.\hfill{$\square$}\\

\bigskip
\noindent
{\bf Acknowledgments}.
This work was initiated when Yun Pu studied at the University of Oxford as a recognized DPhil student
through the Joint Training Ph.D. Program between the University of Oxford and Fudan University -- He would like to express his sincere thanks
to both the home and the host universities for providing him with such a great opportunity.
The research of Gui-Qiang G. Chen is partially supported
by the UK Engineering and Physical Sciences Research Council Awards
EP/L015811/1, EP/V008854/1, and EP/V051121/1.
The research of Yun Pu was partially supported by the Joint Training Ph.D. Program of China Scholarship Council, No. 202006100104.
The research of Yongqian Zhang is supported in part by the NSFC Project 12271507.

\medskip
\noindent
{\bf Conflict of Interest}. The authors declare that they have no conflict of interest. The authors also
declare that this manuscript has not been previously published, and will not be submitted elsewhere
before your decision.

\medskip
\noindent
{\bf Data Availability}: Data sharing is not applicable to this article as no datasets were generated or
analyzed during the current study.

\end{document}